\pdfoutput=1
\RequirePackage{ifpdf}
\ifpdf 
\documentclass[pdftex]{sigma}
\else
\documentclass{sigma}
\fi

\numberwithin{equation}{section}

\newtheorem{Theorem}{Theorem}[section]
 { \theoremstyle{definition}
\newtheorem{Definition}[Theorem]{Definition}
\newtheorem{Example}[Theorem]{Example}
\newtheorem{Remark}[Theorem]{Remark} }

\usepackage[all]{xy}
\usepackage{lscape}

\begin{document}

\allowdisplaybreaks

\newcommand{\arXivNumber}{1703.01379}

\renewcommand{\PaperNumber}{096}

\FirstPageHeading

\ShortArticleName{Four-Dimensional Painlev\'e-Type Equations Associated with Ramif\/ied Linear Equations III}

\ArticleName{Four-Dimensional Painlev\'e-Type Equations\\ Associated with Ramif\/ied Linear Equations III:\\ Garnier Systems and Fuji--Suzuki Systems}

\Author{Hiroshi KAWAKAMI}

\AuthorNameForHeading{H.~Kawakami}

\Address{College of Science and Engineering, Aoyama Gakuin University,\\
5-10-1 Fuchinobe, Chuo-ku, Sagamihara-shi, Kanagawa 252-5258, Japan}
\Email{\href{mailto:kawakami@gem.aoyama.ac.jp}{kawakami@gem.aoyama.ac.jp}}

\ArticleDates{Received April 19, 2017, in f\/inal form December 07, 2017; Published online December 25, 2017}

\Abstract{This is the last part of a series of three papers entitled ``Four-dimensional Painlev\'e-type equations associated with ramif\/ied linear equations''. In this series of papers we aim to construct the complete degeneration scheme of four-dimensional Painlev\'e-type equations. In the present paper, we consider the degeneration of the Garnier system in two variables and the Fuji--Suzuki system.}


\Classification{34M55; 34M56; 33E17}

\section{Introduction}\label{intro}
This is the last part of a series of three papers on four-dimensional Painlev\'e-type equations associated with ramif\/ied linear equations. By the term ``Painlev\'e-type equations'' we mean Hamiltonian systems which describe isomonodromic deformations of linear equations. The isomonodromic deformation is the deformation of a linear dif\/ferential equation which do not change its ``monodromy data'' (see, for example,~\cite{JMU}), and it is known that isomonodromic deformation equations can be written in Hamiltonian form. In this terminology, the classical Painlev\'e equations are Painlev\'e-type equations with two-dimensional phase space.

The classical Painlev\'e equations are non-linear ordinary dif\/ferential equations which were discovered by Painlev\'e~\cite{P} and Gambier~\cite{Gm}.
Originally they were classif\/ied into six equations and are often denoted by $P_{\mathrm{I}}, P_{\mathrm{II}}, \ldots, P_{\mathrm{VI}}$.
However, from a geometric viewpoint, it is natural to classify them into eight equations~\cite{Sak1}. More precisely, the third Painlev\'e equation $P_{\mathrm{III}}$ is divided into three cases $P_{\mathrm{III}(D_6)}$, $P_{\mathrm{III}(D_7)}$, and $P_{\mathrm{III}(D_8)}$. The so-called third Painlev\'e equation is then $P_{\mathrm{III}(D_6)}$.

The eight Hamiltonians associated with the Painlev\'e equations are as follows~\cite{OKSO, Ok++,Ok,Ok+,Ok+++}:
\begin{gather*}
 t(t-1)H_\mathrm{VI}\left({\alpha , \beta \atop \gamma, \delta};t;q,p\right)=q(q-1)(q-t)p^2\\
\qquad{} +\{ \delta q(q-1)-(2\alpha +\beta +\gamma +\delta )q(q-t)+\gamma
 (q-1)(q-t)\} p +\alpha (\alpha +\beta )(q-t),\\
 tH_\mathrm{V}\left({\alpha , \beta \atop \gamma };t;q,p\right)=p(p+t)q(q-1)+\beta pq+\gamma p-(\alpha +\gamma )tq,\\
H_\mathrm{IV}\left(\alpha , \beta;t;q,p\right)=pq(p-q-t)+\beta p+\alpha q,\\
tH_{\mathrm{III}(D_6)}(\alpha , \beta ;t;q,p)=p^2q^2-\big(q^2-\beta q-t\big)p-\alpha q,\\
tH_{\mathrm{III}(D_7)}(\alpha;t;q,p)=p^2q^2+\alpha qp+tp+q,\\
tH_{\mathrm{III}(D_8)}(t;q,p)=p^2q^2+qp-q-\frac{t}{q},\\
H_\mathrm{II}(\alpha;t;q,p)=p^2-\big(q^2+t\big)p-\alpha q,\\
H_\mathrm{I}(t;q,p)=p^2-q^3-tq.
\end{gather*}

The standard linear equations associated with the classical Painlev\'e equations are given by certain second order single linear equations, or equivalently, by f\/irst order $2 \times 2$ systems. Here we review the classif\/ication of the classical Painlev\'e equations in terms of associated linear equations.

It is well-known that the classical Painlev\'e equations admit degeneration. We use the term degeneration in the following sense. Suppose a dif\/ferential equation $E$ has some parameter $\varepsilon$. When the equation $E$ tends to another equation $E'$ as $\varepsilon$ tends to~0,
we say that $E$ degenerates to~$E'$.
The following scheme is the well-known degeneration scheme among the six Painlev\'e equations:
$$
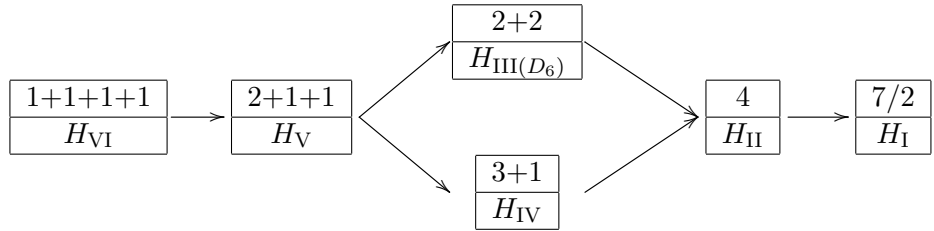
\begin{xy}
{(3,0) *{\begin{tabular}{|c|}
\hline
1+1+1+1\\
\hline
$H_{\rm VI}$\\
\hline
\end{tabular}
}},
{(30,0) *{\begin{tabular}{|c|}
\hline
2+1+1\\
\hline
$H_{\rm V}$\\
\hline
\end{tabular}
}},
{\ar (14,0);(21,0)},
{\ar (39,0);(51,10)},
{\ar (39,0);(51,-10)},
{(60,10) *{\begin{tabular}{|c|}
\hline
2+2\\
\hline
$H_{\mathrm{III}(D_6)}$\\
\hline
\end{tabular}}},
{\ar (69,10);(84,0)},
{\ar (69,-10);(84,0)},
{(60,-10) *{\begin{tabular}{|c|}
\hline
3+1\\
\hline
$H_{\rm IV}$\\
\hline
\end{tabular}}},
{(90,0) *{\begin{tabular}{|c|}
\hline
4\\
\hline
$H_{\rm II}$\\
\hline
\end{tabular}}},
{\ar (96,0);(104,0)},
{(110,0) *{\begin{tabular}{|c|}
\hline
7/2\\
\hline
$H_{\mathrm{I}}$\\
\hline
\end{tabular}}},
\end{xy}
$$
The number in each box is the ``singularity pattern'' of an associated linear equation,
which has information on the Paincar\'e ranks of the singular points of the linear equation.
In this case, the linear equations are well characterized by their singularity pattern.

Here we point out that
\begin{itemize}\itemsep=0pt
\item this scheme lacks the third Painlev\'e equations of type $D_7^{(1)}$ and $D_8^{(1)}$,
\item from the viewpoint of associated linear equations, the degeneration $H_{\mathrm{II}} \to H_{\mathrm{I}}$ is
distinguished from the others. Namely, the other degenerations correspond to the ``conf\/luence of singular points'',
while the degeneration $H_{\mathrm{II}} \to H_{\mathrm{I}}$ corresponds to the ``degeneration of an HTL canonical form''.
\end{itemize}
We note that there are two kinds of degenerations concerning linear equations:
conf\/luence of singular points and degeneration of an HTL canonical form
(where HTL is an abbreviation for Hukuhara--Turrittin--Levelt).

By considering all the possible degenerations of HTL canonical forms (and possible con\-f\/luen\-ces), one can obtain the following complete degeneration scheme of the classical Painlev\'e equations~\cite{Kap,KH, OO}.
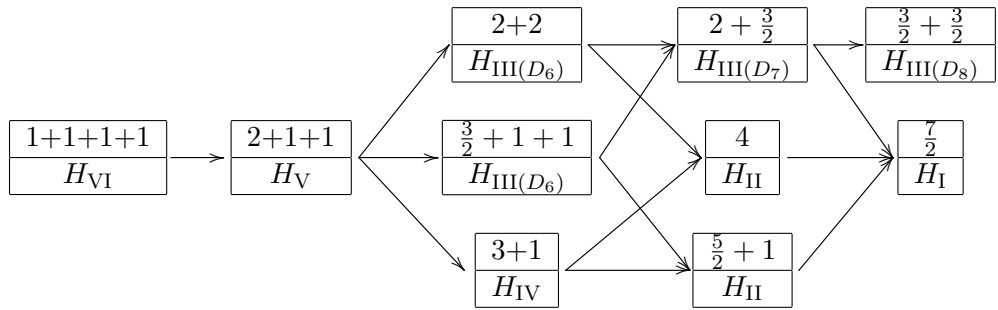
$$
\begin{xy}
{(3,0) *{\begin{tabular}{|c|}
\hline
1+1+1+1\\
\hline
$H_{\mathrm{VI}}$\\
\hline
\end{tabular}
}},
{(30,0) *{\begin{tabular}{|c|}
\hline
2+1+1\\
\hline
$H_{\mathrm{V}}$\\
\hline
\end{tabular}
}},
{\ar (14,0);(21,0)},
{\ar (39,0);(51,15)},
{\ar (39,0);(49.5,0)},
{\ar (39,0);(53,-15)},
{(60,15) *{\begin{tabular}{|c|}
\hline
2+2\\
\hline
$H_{\mathrm{III}(D_6)}$\\
\hline
\end{tabular}}},
{(60,0) *{\begin{tabular}{|c|}
\hline
$\frac{3}{2}+1+1$\\
\hline
$H_{\mathrm{III}(D_6)}$\\
\hline
\end{tabular}}},
{\ar (69.5,15);(81,15)},
{\ar (69.5,15);(84.5,0)},
{\ar (71,0);(81,15)},
{\ar (71,0);(83,-15)},
{\ar (66.5,-15);(84.5,0)},
{\ar (66.5,-15);(83,-15)},
{(60,-15) *{\begin{tabular}{|c|}
\hline
3+1\\
\hline
$H_{\mathrm{IV}}$\\
\hline
\end{tabular}}},
{(90,15) *{\begin{tabular}{|c|}
\hline
$2+\frac{3}{2}$\\
\hline
$H_{\mathrm{III}(D_7)}$\\
\hline
\end{tabular}}},
{(90,0) *{\begin{tabular}{|c|}
\hline
4\\
\hline
$H_{\mathrm{II}}$\\
\hline
\end{tabular}}},
{\ar (99.5,15);(110,0)},
{\ar (99.5,15);(106,15)},
{\ar (96,0);(110,0)},
{\ar (97.5,-15);(110,0)},
{(90,-15) *{\begin{tabular}{|c|}
\hline
$\frac{5}{2}+1$\\
\hline
$H_{\mathrm{II}}$\\
\hline
\end{tabular}}},
{(115,15) *{\begin{tabular}{|c|}
\hline
$\frac{3}{2}+\frac{3}{2}$\\
\hline
$H_{\mathrm{III}(D_8)}$\\
\hline
\end{tabular}}},
{(115,0) *{\begin{tabular}{|c|}
\hline
$\frac{7}{2}$\\
\hline
$H_{\mathrm{I}}$\\
\hline
\end{tabular}}},
\end{xy}
$$
Thus we can say that the sixth Painlev\'e equation is the ``master equation'' from which we can derive
all the other Painlev\'e equations by successive degenerations.

Recently, many generalizations of the Painlev\'e equations have been discovered.
What is important here is that they can be written as compatibility conditions of linear dif\/ferential equations.
Thus they can be regarded as isomonodromic deformation equations of some linear dif\/ferential equations.

The purpose of our study is to understand those of four-dimensional phase space using associated linear equations.
More specif\/ically, we construct the degeneration scheme starting from suitable Painlev\'e-type equations so that we can provide a unif\/ied perspective of four-dimensional analogues of the Painlev\'e equations.

In~\cite{KNS},
the degeneration scheme of four-dimensional Painlev\'e-type equations associated with \textit{unramified} linear equations was obtained.
They considered degenerations from the following four Painlev\'e-type equations:
\begin{gather*}
H_{\mathrm{Gar}}^{1+1+1+1+1}, \qquad
H_{\mathrm{FS}}^{A_5}, \qquad
H_{\mathrm{Ss}}^{D_6}, \qquad
H_{\mathrm{VI}}^{\mathrm{Mat}}.
\end{gather*}
These stand for the Hamiltonians for the four-dimensional Painlev\'e-type equations:
the Garnier system in two variables~\cite{G}, the Fuji--Suzuki system~\cite{FS1, Ts}, the Sasano system~\cite{Ss},
and the sixth matrix Painlev\'e system~\cite{B2, K, KNS}.
They correspond to four Fuchsian systems characterized by the following spectral types:
\begin{gather*}
11,11,11,11,11 \qquad 21,21,111,111 \qquad 31,22,22,1111 \qquad 22,22,22,211
\end{gather*}
respectively.
The above four Painlev\'e-type equations are the master equations in four-di\-men\-sional case.
We note that the degenerations considered in the paper correspond to ``conf\/luence of singular points'' of linear equations.

The aim of the present series of papers is to obtain the ``complete'' degeneration scheme of the Painlev\'e-type equations with four-dimensional phase space.
In order to achieve this goal, we have to consider the degeneration of HTL canonical forms
in four-dimensional case.
Note that, unlike the two-dimensional case, degenerations of HTL canonical forms treated in this study
do not necessarily correspond to degenerations of Jordan canonical forms (see Section~\ref{sec:deg_HTL} and \cite{K2}).

In Part I and Part II \cite{K1, K2} we consider degenerations from the $22,22,22,211$-system and the $31,22,22,1111$-system.
These correspond to the degeneration of the sixth matrix Painlev\'e system and the Sasano system, respectively.

In this Part III, we treat the degeneration of the Garnier system in two variables and what we call the Fuji--Suzuki system.

Note that the degeneration scheme of the Garnier system corresponding to conf\/luences of singular points
(and a degeneration of an HTL canonical form) was already obtained by Kimura~\cite{Km}.
The complete list of degenerate Garnier systems associated with ramif\/ied linear equations was obtained by Kawamuko~\cite{Kw}.
In this paper, we recalculated the Hamiltonians and Lax pairs for the degenerate Garnier systems.

This paper is organized as follows.
In Section~\ref{sec:LDE} we recall the notions of HTL canonical forms, Riemann schemes, singularity patterns, and spectral types.
We also discuss the degeneration of HTL canonical forms.
In Section~\ref{sec:Lax} we present the linear equations of ramif\/ied type which
can be obtained by degeneration from the $11,11,11,11,11$-system and the $21,21,111,111$-system and the corresponding Hamiltonians.
In Section~\ref{sec:Laplace} we discuss correspondences of linear systems through the Laplace transform.
In Section~\ref{sec:conclusion} we give the degeneration scheme of the Hamiltonians of four-dimensional Painlev\'e-type equations.
In the Appendix we give data on degenerations.

Here we present the degeneration scheme of the Garnier system and the Fuji--Suzuki system,
which is the main result of this paper.

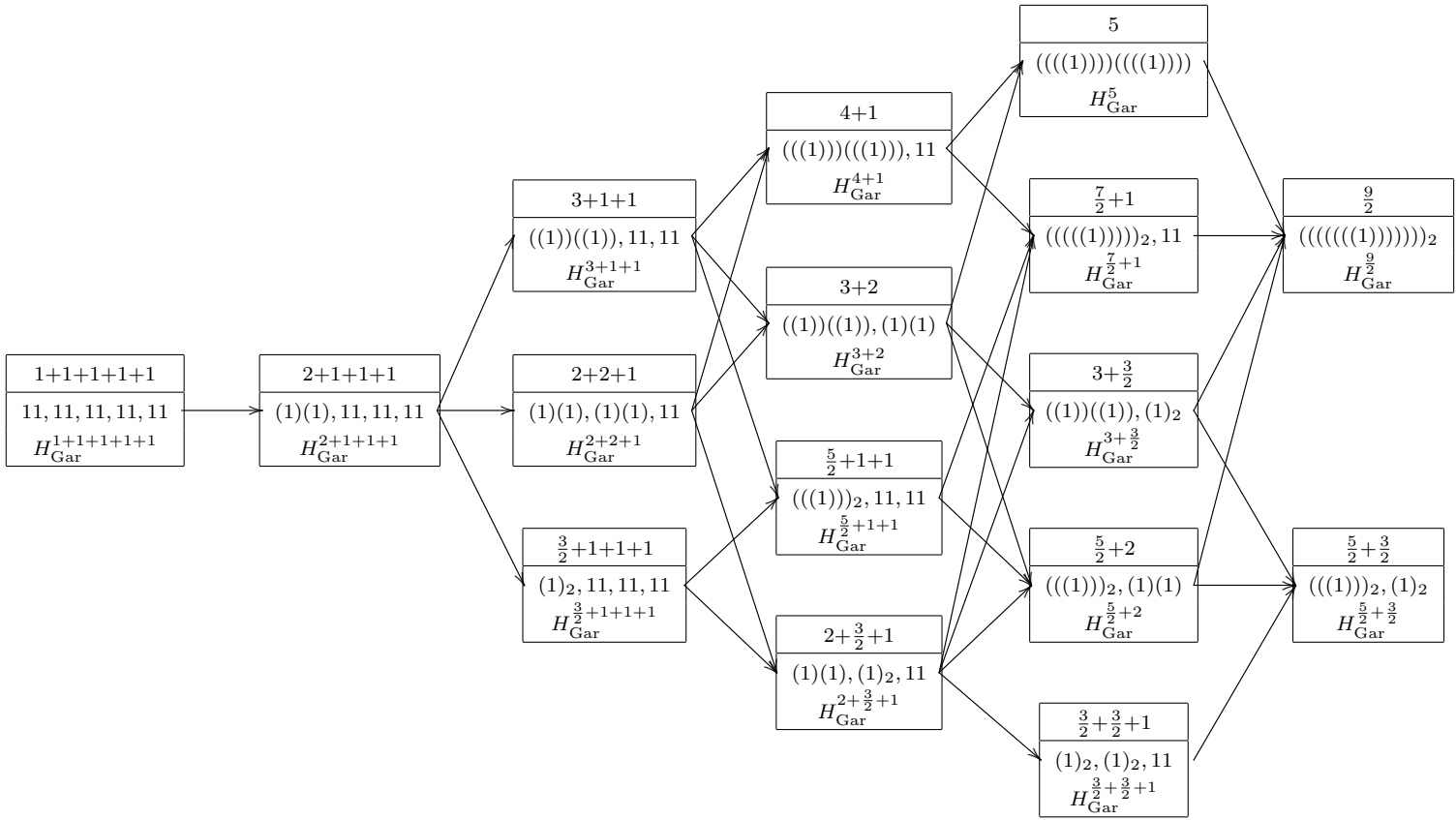
\begin{landscape} \vspace*{20mm}
{\scriptsize
\begingroup
\renewcommand{\arraystretch}{1.5}
\begin{xy}
{(0,-12) *{\begin{tabular}{|c|}
\hline
1+1+1+1+1\\
\hline
$11,11,11,11,11$\\
$H_{\mathrm{Gar}}^{1+1+1+1+1}$\\
\hline
\end{tabular}
}},
{\ar (11.9,-12);(23,-12)},
{(35,-12) *{\begin{tabular}{|c|}
\hline
2+1+1+1\\
\hline
$(1)(1),11,11,11$\\
$H_{\mathrm{Gar}}^{2+1+1+1}$\\
\hline
\end{tabular}
}},
{\ar (47,-12);(57.6,12)},
{\ar (47,-12);(57.6,-12)},
{\ar (47,-12);(59,-36)},
{(70,12) *{\begin{tabular}{|c|}
\hline
3+1+1\\
\hline
$((1))((1)),11,11$\\
$H_{\mathrm{Gar}}^{3+1+1}$\\
\hline
\end{tabular}}},
{(70,-12) *{\begin{tabular}{|c|}
\hline
2+2+1\\
\hline
$(1)(1),(1)(1),11$\\
$H_{\mathrm{Gar}}^{2+2+1}$\\
\hline
\end{tabular}}},
{(70,-36) *{\begin{tabular}{|c|}
\hline
$\frac32$+1+1+1\\
\hline
$(1)_2,11,11,11$\\
$H_{\mathrm{Gar}}^{\frac32+1+1+1}$\\
\hline
\end{tabular}}},
{\ar (82,12);(92.6,24)},
{\ar (82,12);(92.6,0)},
{\ar (82,12);(93.9,-24)},
{\ar (82,-12);(92.6,24)},
{\ar (82,-12);(92.6,0)},
{\ar (82,-12);(93.7,-48)},
{\ar (81,-36);(93.9,-24)},
{\ar (81,-36);(93.7,-48)},
{(105,24) *{\begin{tabular}{|c|}
\hline
4+1\\
\hline
$(((1)))(((1))),11$\\
$H_{\mathrm{Gar}}^{4+1}$\\
\hline
\end{tabular}}},
{(105,0) *{\begin{tabular}{|c|}
\hline
3+2\\
\hline
$((1))((1)),(1)(1)$\\
$H_{\mathrm{Gar}}^{3+2}$\\
\hline
\end{tabular}}},
{(105,-24) *{\begin{tabular}{|c|}
\hline
$\frac52$+1+1\\
\hline
$(((1)))_2,11,11$\\
$H_{\mathrm{Gar}}^{\frac52+1+1}$\\
\hline
\end{tabular}}},
{(105,-48) *{\begin{tabular}{|c|}
\hline
2+$\frac32$+1\\
\hline
$(1)(1),(1)_2,11$\\
$H_{\mathrm{Gar}}^{2+\frac32+1}$\\
\hline
\end{tabular}}},
{\ar (117,24);(127.5,36)},
{\ar (117,24);(129,12)},
{\ar (117,0);(127.5,36)},
{\ar (117,0);(129,-12)},
{\ar (117,0);(128.7,-36)},
{\ar (116,-24);(129,12)},
{\ar (116,-24);(128.7,-36)},
{\ar (116,-48);(129,12)},
{\ar (116,-48);(129,-12)},
{\ar (116,-48);(128.7,-36)},
{\ar (116,-48);(130,-60)},
{(140,36) *{\begin{tabular}{|c|}
\hline
5\\
\hline
$((((1))))((((1))))$\\
$H_{\mathrm{Gar}}^5$\\
\hline
\end{tabular}}},
{(140,12) *{\begin{tabular}{|c|}
\hline
$\frac72$+1\\
\hline
$(((((1)))))_2,11$\\
$H_{\mathrm{Gar}}^{\frac72+1}$\\
\hline
\end{tabular}}},
{(140,-12) *{\begin{tabular}{|c|}
\hline
3+$\frac32$\\
\hline
$((1))((1)),(1)_2$\\
$H_{\mathrm{Gar}}^{3+\frac32}$\\
\hline
\end{tabular}}},
{(140,-36) *{\begin{tabular}{|c|}
\hline
$\frac52$+2\\
\hline
$(((1)))_2,(1)(1)$\\
$H_{\mathrm{Gar}}^{\frac52+2}$\\
\hline
\end{tabular}}},
{(140,-60) *{\begin{tabular}{|c|}
\hline
$\frac32$+$\frac32$+1\\
\hline
$(1)_2,(1)_2,11$\\
$H_{\mathrm{Gar}}^{\frac32+\frac32+1}$\\
\hline
\end{tabular}}},
{\ar (152.4,36);(163.5,12)},
{\ar (151,12);(163.5,12)},
{\ar (151,-12);(163.5,12)},
{\ar (151,-12);(165,-36)},
{\ar (151,-36);(163.5,12)},
{\ar (151,-36);(165,-36)},
{\ar (151,-60);(165,-36)},
{(175,12) *{\begin{tabular}{|c|}
\hline
$\frac92$\\
\hline
$(((((((1)))))))_2$\\
$H_{\mathrm{Gar}}^{\frac92}$\\
\hline
\end{tabular}}},
{(175,-36) *{\begin{tabular}{|c|}
\hline
$\frac52$+$\frac32$\\
\hline
$(((1)))_2,(1)_2$\\
$H_{\mathrm{Gar}}^{\frac52+\frac32}$\\
\hline
\end{tabular}}},
\end{xy}
\endgroup
}
\end{landscape}

\begin{landscape}\vspace*{5mm}
{\scriptsize
\begingroup
\renewcommand{\arraystretch}{1.5}
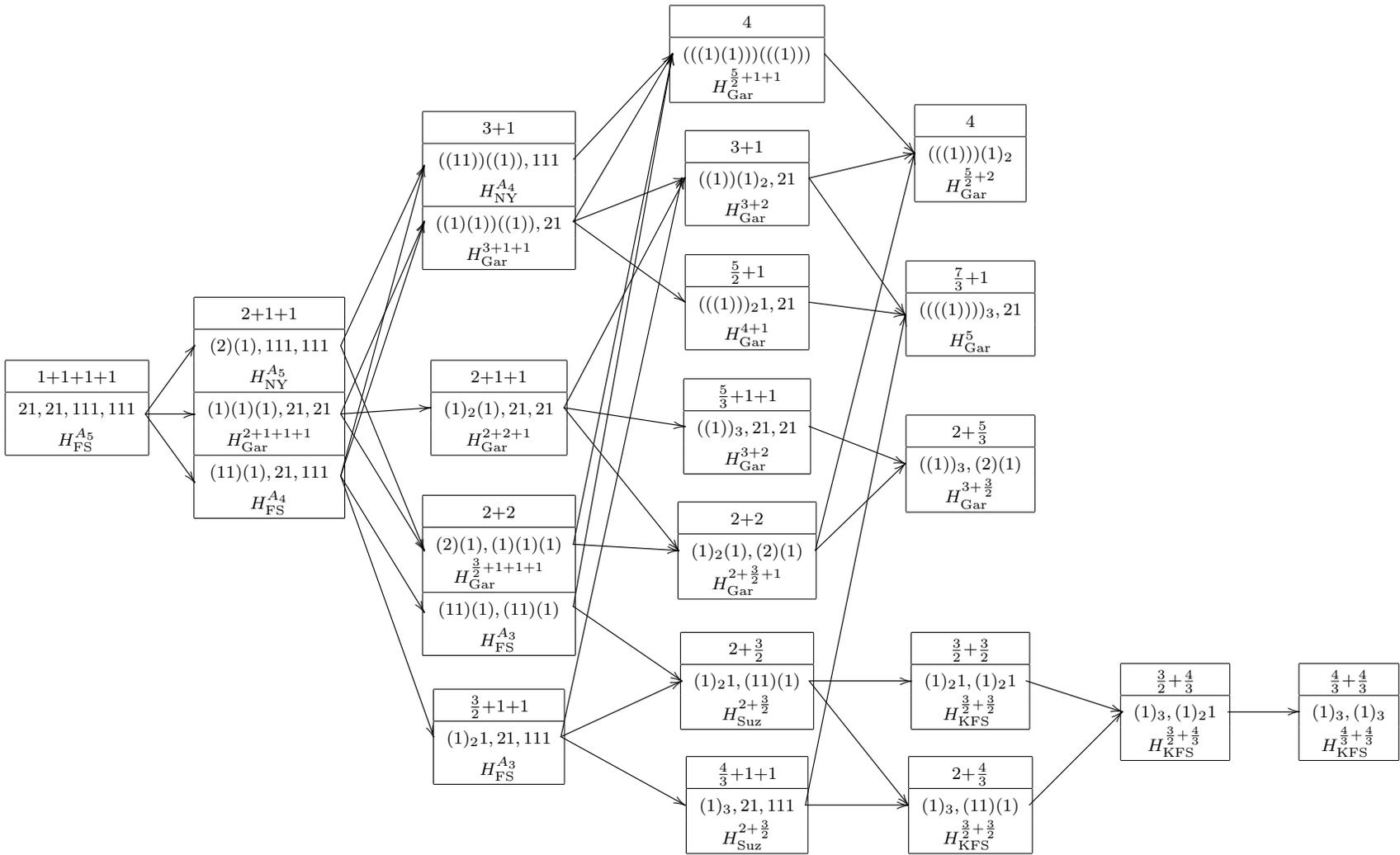
\begin{xy}
{(0,-30) *{\begin{tabular}{|c|}
\hline
1+1+1+1\\
\hline
$21,21,111,111$\\
$H_{\mathrm{FS}}^{A_5}$\\
\hline
\end{tabular}
}},
{\ar (11,-31);(19,-20)},
{\ar (11,-31);(19,-31)},
{\ar (11,-31);(19,-42)},
{(31,-30) *{\begin{tabular}{|c|}
\hline
2+1+1\\
\hline
$(2)(1),111,111$\\
$H_{\mathrm{NY}}^{A_5}$\\
\hline
$(1)(1)(1),21,21$\\
$H_{\mathrm{Gar}}^{2+1+1+1}$\\
\hline
$(11)(1),21,111$\\
$H_{\mathrm{FS}}^{A_4}$\\
\hline
\end{tabular}
}},
{\ar (42.5,-20);(56,9)},
{\ar (42.5,-20);(56,-53)},
{\ar (42.5,-31);(56,0)},
{\ar (42.5,-31);(57,-30)},
{\ar (42.5,-31);(56,-53)},
{\ar (42.5,-41);(56,9)},
{\ar (42.5,-41);(56,0)},
{\ar (42.5,-41);(56,-63)},
{\ar (42.5,-41);(57.5,-83)},
{(68,-30) *{\begin{tabular}{|c|}
\hline
2+1+1\\
\hline
$(1)_2(1), 21, 21$\\
$H_{\mathrm{Gar}}^{2+2+1}$\\
\hline
\end{tabular}}},
{(68,5) *{\begin{tabular}{|c|}
\hline
3+1\\
\hline
$((11))((1)),111$\\
$H_{\mathrm{NY}}^{A_4}$\\
\hline
$((1)(1))((1)),21$\\
$H_{\mathrm{Gar}}^{3+1+1}$\\
\hline
\end{tabular}}},
{(68,-57) *{\begin{tabular}{|c|}
\hline
2+2\\
\hline
$(2)(1),(1)(1)(1)$\\
$H_{\mathrm{Gar}}^{\frac{3}{2}+1+1+1}$\\
\hline
$(11)(1),(11)(1)$\\
$H_{\mathrm{FS}}^{A_3}$\\
\hline
\end{tabular}}},
{(68,-83) *{\begin{tabular}{|c|}
\hline
$\frac32$+1+1\\
\hline
$(1)_21,21,111$\\
$H_{\mathrm{FS}}^{A_3}$\\
\hline
\end{tabular}}},
{\ar (80,10);(96,27)},
{\ar (80,0);(96,27)},
{\ar (80,0);(98,7)},
{\ar (80,0);(98,-13)},
{\ar (78.5,-30);(98,7)},
{\ar (78.5,-30);(98,-33)},
{\ar (78.5,-30);(97,-53)},
{\ar (80,-52);(96,27)},
{\ar (80,-52);(97,-53)},
{\ar (80,-62);(96,27)},
{\ar (80,-62);(97.5,-74)},
{\ar (78,-83);(97.5,7)},
{\ar (78,-83);(97.5,-74)},
{\ar (78,-83);(98,-94)},
{(108,27) *{\begin{tabular}{|c|}
\hline
4\\
\hline
$(((1)(1)))(((1)))$\\
$H_{\mathrm{Gar}}^{\frac{5}{2}+1+1}$\\
\hline
\end{tabular}}},
{(108,7) *{\begin{tabular}{|c|}
\hline
3+1\\
\hline
$((1))(1)_2, 21$\\
$H_{\mathrm{Gar}}^{3+2}$\\
\hline
\end{tabular}}},
{(108,-13) *{\begin{tabular}{|c|}
\hline
$\frac52$+1\\
\hline
$(((1)))_21, 21$\\
$H_{\mathrm{Gar}}^{4+1}$\\
\hline
\end{tabular}}},
{(108,-33) *{\begin{tabular}{|c|}
\hline
$\frac53$+1+1\\
\hline
$((1))_3, 21, 21$\\
$H_{\mathrm{Gar}}^{3+2}$\\
\hline
\end{tabular}}},
{(108,-53) *{\begin{tabular}{|c|}
\hline
2+2\\
\hline
$(1)_2(1), (2)(1)$\\
$H_{\mathrm{Gar}}^{2+\frac32+1}$\\
\hline
\end{tabular}}},
{(108,-74) *{\begin{tabular}{|c|}
\hline
2+$\frac32$\\
\hline
$(1)_21, (11)(1)$\\
$H_{\mathrm{Suz}}^{2+\frac32}$\\
\hline
\end{tabular}}},
{(108,-94) *{\begin{tabular}{|c|}
\hline
$\frac43$+1+1\\
\hline
$(1)_3, 21, 111$\\
$H_{\mathrm{Suz}}^{2+\frac32}$\\
\hline
\end{tabular}}},
{\ar (120.5,27);(135,11)},
{\ar (118,7);(135,11)},
{\ar (118,7);(133.5,-15)},
{\ar (118,-13);(133.5,-15)},
{\ar (118,-33);(133.5,-39)},
{\ar (119,-53);(133.5,-39)},
{\ar (119,-53);(135,11)},
{\ar (118,-74);(134.5,-74)},
{\ar (118,-74);(134,-94)},
{\ar (117.4,-94);(134,-94)},
{\ar (117.4,-94);(133.5,-15)},
{(144,11) *{\begin{tabular}{|c|}
\hline
4\\
\hline
$(((1)))(1)_2$\\
$H_{\mathrm{Gar}}^{\frac52+2}$\\
\hline
\end{tabular}}},
{(144,-14) *{\begin{tabular}{|c|}
\hline
$\frac73$+1\\
\hline
$((((1))))_3, 21$\\
$H_{\mathrm{Gar}}^{5}$\\
\hline
\end{tabular}}},
{(144,-39) *{\begin{tabular}{|c|}
\hline
2+$\frac53$\\
\hline
$((1))_3, (2)(1)$\\
$H_{\mathrm{Gar}}^{3+\frac32}$\\
\hline
\end{tabular}}},
{(144,-74) *{\begin{tabular}{|c|}
\hline
$\frac32$+$\frac32$\\
\hline
$(1)_21, (1)_21$\\
$H_{\mathrm{KFS}}^{\frac32+\frac32}$\\
\hline
\end{tabular}}},
{(144,-94) *{\begin{tabular}{|c|}
\hline
2+$\frac43$\\
\hline
$(1)_3, (11)(1)$\\
$H_{\mathrm{KFS}}^{\frac32+\frac32}$\\
\hline
\end{tabular}}},
{\ar (153,-74);(168.5,-79)},
{\ar (153.5,-94);(168.5,-79)},
{(177,-79) *{\begin{tabular}{|c|}
\hline
$\frac32$+$\frac43$\\
\hline
$(1)_3, (1)_21$\\
$H_{\mathrm{KFS}}^{\frac32+\frac43}$\\
\hline
\end{tabular}}},
{\ar (185.5,-79);(197,-79)},
{(205,-79) *{\begin{tabular}{|c|}
\hline
$\frac43$+$\frac43$\\
\hline
$(1)_3, (1)_3$\\
$H_{\mathrm{KFS}}^{\frac43+\frac43}$\\
\hline
\end{tabular}}},
\end{xy}
\endgroup
}
\end{landscape}

Symbols such as $H^{2+1+1+1}_{\mathrm{Gar}}$ stand for the Hamiltonians for four-dimensional Painlev\'e-type equations.
The explicit expressions for them are given in Section~\ref{sec:Lax}.

\section{Notions on linear dif\/ferential equations}\label{sec:LDE}
In this section we recall some notions on linear dif\/ferential equations.
\subsection{HTL canonical forms}\label{sec:HTLform}
Consider a system of linear dif\/ferential equations
\begin{gather*}
\frac{{\rm d}Y}{{\rm d}x}=A(x)Y.
\end{gather*}
A change of the dependent variable $Y=PZ$ with an invertible matrix $P=P(x)$ yields the following system:
\begin{gather*}
\frac{{\rm d}Z}{{\rm d}x}=\big( P^{-1}A(x)P-P^{-1}P' \big)Z.
\end{gather*}
The transformation
\begin{gather*}
A(x) \mapsto A^P(x):=P^{-1}A(x)P-P^{-1}P'
\end{gather*}
is called the \textit{gauge transformation} by $P$.

A system of linear dif\/ferential equations with rational function coef\/f\/icients
\begin{gather}\label{eq:rational_LDE}
\frac{{\rm d}Y}{{\rm d}x}=
\left(\sum_{\nu=1}^n\sum_{k=0}^{r_{\nu}}\frac{A_{\nu}^{(k)}}{(x-u_{\nu})^{k+1}}
+\sum_{k=1}^{r_{\infty}}A_{\infty}^{(k)}x^{k-1}
\right)Y,\qquad
A_*^{(k)} \in M(m, \mathbb{C})
\end{gather}
can be transformed into the ``canonical form'' at each singular point.

The system (\ref{eq:rational_LDE}) has singularity at $x=u_\nu$ $(\nu=1, \ldots, n)$ and $x=u_0:=\infty$.
Set $z=x-u_{\nu}$ $(\nu=1, \ldots, n)$ or $z=1/x$.
We consider the system around $z=0$:
\begin{gather*}
\frac{{\rm d}Y}{{\rm d}z}=\left( \frac{A_0}{z^{r+1}}+\frac{A_1}{z^{r}}+\cdots+A_{r+1}+A_{r+2} z+\cdots\right)Y.
\end{gather*}

Let $\mathcal{P}_z=\cup_{p > 0}\mathbb{C}\big(\!\big(z^{1/p}\big)\!\big)$ be the f\/ield of Puiseux series where $\mathbb{C}(\!(t)\!)$ is the f\/ield of formal Laurent series in $t$.

\begin{Definition}[HTL canonical form]
An element
\begin{gather*}
\frac{D_0}{z^{l_0}}+\frac{D_1}{z^{l_1}}+\cdots+\frac{D_{s-1}}{z^{l_{s-1}}}+\frac{\Theta}{z^{l_s}}
\end{gather*}
in $M_m(\mathcal{P}_z)$ satisfying the following conditions:
\begin{itemize}\itemsep=0pt
\item $l_j$ is a rational number, $l_0 > l_1 > \cdots > l_{s-1} > l_s=1$,
\item $D_0, \ldots, D_{s-1}$ are commuting diagonalizable matrices,
\item $\Theta$ is a (not necessarily diagonalizable) matrix that commutes with all $D_j$'s
\end{itemize}
is called an \textit{HTL canonical form}, or \textit{HTL form} for short.
\end{Definition}

\begin{Theorem}[Hukuhara~\cite{Huk}, Levelt~\cite{Lev}, Turrittin~\cite{Tur}]
For any
\begin{gather*}
A(z)=\frac{A_0}{z^{r+1}}+\frac{A_1}{z^{r}}+\cdots \in M_m(\mathbb{C}(\!(z)\!)),
\end{gather*}
there exists $P \in \mathrm{GL}(m, \mathcal{P}_z)$ such that $A^P(z)$ is an HTL form
\begin{gather*}
\frac{D_0}{z^{l_0}}+\frac{D_1}{z^{l_1}}+\cdots+\frac{D_{s-1}}{z^{l_{s-1}}}
+\frac{\Theta}{z}.
\end{gather*}
Here $l_0,\ldots,l_{s-1}$ are uniquely determined only by $A(z)$.

If
\begin{gather*}
\frac{\tilde{D}_0}{z^{l_0}}+\frac{\tilde{D}_1}{z^{l_1}}+\cdots+\frac{\tilde{D}_{s-1}}{z^{l_{s-1}}}
+\frac{\tilde{\Theta}}{z}
\end{gather*}
is another HTL form of the same $A(z)$,
then there exist $g \in \mathrm{GL}(m, \mathbb{C})$ and $k \in \mathbb{Z}_{\ge 1}$ such that
\begin{gather*}
\tilde{D}_j=g^{-1}D_j g, \qquad \exp\big(2\pi i k \tilde{\Theta}\big)=g^{-1}\exp(2\pi i k \Theta)g
\end{gather*}
hold.
\end{Theorem}
There is an algorithm for constructing the HTL forms of a given linear system.
We will brief\/ly explain how to construct HTL forms in Section~\ref{sec:ST}.

The number $l_0-1$ is called the \textit{Poincar\'e rank} of the singular point.
If there is a rational number $l_j \in \mathbb{Q} \setminus \mathbb{Z}$, the singular point is called a \textit{ramified} irregular singular point.
A linear equation is said to be of \textit{ramified type} if the equation has a ramif\/ied irregular singular point.

When we roughly express the singularity of a linear dif\/ferential equation, we attach the number ``Poincar\'e rank + 1'' to
each singular point and connect the numbers with ``+''.
We call it the \textit{singularity pattern} of the equation.

\subsection{Riemann schemes}
Consider an HTL form
\begin{gather}\label{eq:HTLform}
\frac{D_0}{z^{l_0}}+\frac{D_1}{z^{l_1}}+\cdots+\frac{D_{s-1}}{z^{l_{s-1}}}+\frac{\Theta}{z}.
\end{gather}
Here we assume $D_j$'s and $\Theta$ are in Jordan canonical form.

Let $d \in \mathbb{Z}_{>0}$ be the minimum element of
\begin{gather*}
\{ k \in \mathbb{Z}_{>0} \,|\, k l_j \in \mathbb{Z} \ (j=0,\ldots, s-1) \}.
\end{gather*}
Notice that $d=1$ is equivalent to the unramif\/iedness of the singular point. Then the HTL form~(\ref{eq:HTLform}) can be rewritten as
\begin{gather*}
\frac{T_0}{z^{\frac{b}{d}+1}}+\frac{T_1}{z^{\frac{b-1}{d}+1}}+\cdots+\frac{T_{b-1}}{z^{\frac{1}{d}+1}}
+\frac{\Theta}{z},
\end{gather*}
where $l_0=\frac{b}{d}+1$.

When $\Theta$ is a diagonal matrix, we denote this HTL form by
\begin{gather*}
\begin{array}{c}
x=u_i \ \left(\frac{1}{d}\right) \\
\overbrace{\begin{array}{ccccc}
 t^0_1 & t^1_1 & \ldots & t^{b-1}_1 & \theta_1\\
 \vdots & \vdots & & \vdots & \vdots \\
 t^0_m & t^1_m & \ldots & t^{b-1}_m & \theta_m
 \end{array}}
\end{array}
\end{gather*}
where $T_j=\operatorname{diag}(t^j_1,\ldots,t^j_m)$, $\Theta=\operatorname{diag}(\theta_1,\ldots,\theta_m)$.
In the case of $d=1$, we omit $\left( \frac{1}{d} \right)$.
\begin{Remark}
In this series of papers, the matrix $\Theta$ is always diagonalizable.
\end{Remark}

The table of the HTL forms (represented by the above formula) at all singular points is called the \textit{Riemann scheme} of a linear equation.

Concerning the computation of HTL forms, the following two theorems are fundamental.

\begin{Theorem}[\cite{Wa}]\label{thm:first_red}
For any
\begin{align*}
A(z)=\frac1z \left( A_0+A_1z+\cdots \right),
\end{align*}
there exists a formal power series $P=P(z)$ such that
\begin{align*}
A^P(z)=\frac{\tilde{A}_0}{z}.
\end{align*}
When no two different eigenvalues of $A_0$ differ by an integer, we can choose the gauge so that $\tilde{A}_0=A_0$.
\end{Theorem}

\begin{Theorem}[block diagonalization \cite{Wa}]\label{thm:block_diag}
Let $A(z)$ be
\begin{gather*}
A(z)=\frac{1}{z^{r+1}} ( A_0+A_1z+\cdots  ), \qquad r \in \mathbb{Z}_{>0},
\end{gather*}
where the eigenvalues of $A_0$ are assumed to be $\lambda_1,\ldots,\lambda_n$.
Without loss of generality, we can assume that $A_0$ is in Jordan canonical form
\begin{gather*}
A_0=J_1(\lambda_1) \oplus \cdots \oplus J_n(\lambda_n)
\end{gather*}
where $J_k(\lambda_k)$ is a direct sum of Jordan blocks with the eigenvalue $\lambda_k$.

Then there exists a formal power series $P=P(z)$ such that
\begin{gather*}
A^P(z)
=
B_1(z) \oplus \cdots \oplus B_n(z),
\qquad
B_k(z)=\frac{1}{z^{r+1}}\big( B^k_0+B^k_1z+\cdots \big),
\end{gather*}
where $B^k_0=J_k(\lambda_k)$.
\end{Theorem}

\subsection{Spectral types}\label{sec:ST}
We have introduced the notion of singularity pattern to represent the singularity of a linear system.
However, if the rank of a linear system is greater than two,
the singularity pattern is too rough to describe the singularity of the linear system suf\/f\/iciently
since it only has information of the Poincar\'e ranks.

In order to describe the singularity of a linear system in more detail,
we need the notion of spectral type of linear equations.

Spectral types are def\/ined through HTL forms.
\begin{itemize}\itemsep=0pt
\item In the unramif\/ied case, the spectral type is def\/ined to be the (tuple of) ``ref\/ining sequence of partitions (RSP)''.
\item In the ramif\/ied case, the spectral type consists of ``copies'' of RSPs.
\end{itemize}
Thus we begin with the unramif\/ied case.

\subsubsection{Spectral types of unramif\/ied irregular singularities}
Consider a system of linear dif\/ferential equations whose coef\/f\/icient matrix is
\begin{gather}\label{eq:irreg_sys}
A(z)=\frac{1}{z^{r+1}}\left( A_0+A_1z+\cdots \right), \qquad r \in \mathbb{Z}_{>0}.
\end{gather}
Now suppose that $A_0$ is diagonalizable.
Applying the block diagonalization (see Theorem~\ref{thm:block_diag}) to (\ref{eq:irreg_sys}),
we see that the leading terms $B^k_0$ $(k=1,\ldots,n)$ are all scalar matrices.

We focus on the next $B^k_1$ in each block.
Suppose that all the $B^k_1$ $(k=1, \ldots, n)$ are also diagonalizable.
Then diagonalize them by constant gauge transformations
and apply Theorem~\ref{thm:block_diag} to each block.
Thus each block decomposes into smaller direct summands again.
We note that the f\/irst two terms of each direct summands thus obtained are scalar matrices.

In general, if $B_{r-j+1}$ is diagonalizable in a direct summand
\begin{gather*}
\frac{c_0 I}{z^{r+1}}+\frac{c_1 I}{z^{r}}+\cdots+\frac{c_{r-j}I}{z^{j+1}}+
\frac{B_{r-j+1}}{z^j}+\cdots,
\end{gather*}
then, by a gauge transformation, we have
\begin{align*}
&\sim
\frac{c_0 I}{z^{r+1}}+\frac{c_1 I}{z^{r}}+\cdots+\frac{c_{r-j}I}{z^{j+1}}+
\frac{D_{r-j+1}}{z^j}+\cdots
\end{align*}
where $D_{r-j+1}=d_1 I_{m_1} \oplus \cdots \oplus d_l I_{m_l}$. After the block diagonalization with respect to the eigenvalues of $D_{r-j+1}$,
the above series decomposes into smaller direct summands
\begin{gather*}
\frac{c_0 I_{m_i}}{z^{r+1}}+\frac{c_1 I_{m_i}}{z^{r}}+\cdots+\frac{c_{r-j}I_{m_i}}{z^{j+1}}+
\frac{d_i I_{m_i}}{z^j}+\frac{*}{z^{j-1}}+\cdots, \qquad i=1\ldots,l.
\end{gather*}

By repeating the above procedure, we can decompose (\ref{eq:irreg_sys}) into a direct sum of the series of the following form
\begin{gather*}
\frac{c_0 I}{z^{r+1}}+\frac{c_1 I}{z^r}+\cdots+\frac{c_{r-1} I}{z^2}+\frac{C_r}{z}+\cdots.
\end{gather*}
By virtue of Theorem~\ref{thm:first_red}, we can eliminate all the terms with non-negative powers of~$z$
by a~suitable gauge transformation. As the result, (\ref{eq:irreg_sys}) transforms into a direct sum of matrices of the following form
\begin{gather*}
\frac{c_0 I}{z^{r+1}}+\frac{c_1 I}{z^r}+\cdots+\frac{c_{r-1} I}{z^2}+\frac{\tilde{C}_r}{z}.
\end{gather*}
The direct sum of the above forms thus obtained is the HTL canonical form.

In this case (i.e., the unramif\/ied case), by construction, the feature of an HTL form can be well represented by a ``ref\/ining sequence of partitions''.

\begin{Definition}[ref\/ining sequence of partitions~\cite{KNS}]
Let $\lambda=\lambda_1\ldots \lambda_p$, $\mu=\mu_1\ldots \mu_q$ be partitions of a natural number~$m$:
\begin{gather*}
\lambda_1+\cdots+\lambda_p=\mu_1+\cdots+\mu_q=m.
\end{gather*}
Here we assume that $\lambda_i$'s and $\mu_i$'s are not necessarily arranged in ascending or descending order.

If there exists a disjoint decomposition $\{ 1,2,\ldots , p\}=I_1\coprod \cdots \coprod I_q$ of the index set of $\lambda$
such that $\mu_k=\sum\limits_{j\in I_k}\lambda_j$ holds, then we call $\lambda$ a \textit{refinement} of $\mu$.

Let $[p_0,\ldots,p_r]$ be an $(r+1)$-tuple of partitions of $m$. When $p_{i+1}$ is a ref\/inement of~$p_i$ for all~$i$ $(i=0,\ldots,r-1)$, we call $[p_0,\ldots,p_r]$ a \textit{refining sequence of partitions},
or \textit{RSP} for short.
\end{Definition}

We denote an RSP in the following way.
\begin{Example}
We consider the following RSP
\begin{gather*}
[321, 2121, 111111]
\end{gather*}
as an example.

First, write the rightmost partition:
\begin{gather*}
111111.
\end{gather*}
Second, put the numbers that are grouped together in the central partition in parentheses:
\begin{gather*}
(11)(1)(11)(1).
\end{gather*}
Finally, put the numbers that are grouped together in the leftmost partition in parentheses:
\begin{gather*}
((11)(1))((11))((1)).
\end{gather*}
Thus we can denote the above RSP by $((11)(1))((11))((1))$.
\end{Example}

\subsubsection{Spectral types of ramif\/ied irregular singularities}\label{sec:STram}

In general, non-semisimple matrices may appear in a sequence of block diagonalizations. Such a~case can be reduced to the above semisimple case by ``shearing transformations''. Here we point out that, usually, the non-semisimplicity implies the ramif\/iedness of a singular point.

A shearing transformation is a gauge transformation by a diagonal matrix, which is typically of the form
\begin{gather*}
S=\operatorname{diag}\big(1,z^s,\ldots,z^{(m-1)s}\big), \end{gather*}
where $s$ is a positive rational number.
The aim of the shearing transformation is to make non-semisimple coef\/f\/icient matrices semisimple by repeating gauge transformations
of the above kind~\cite{Huk, Tur}.
Instead of describing shearing transformations and constructions of HTL forms in the ramif\/ied case in full generality,
we demonstrate a construction of an HTL form using the linear system~(\ref{eq:(1)_3,(11)(1)}).
A general method of constructing HTL forms can be found in \cite{Wa} (see also \cite{K2}).

First we change the dependent and independent variables as $z=1/x, \, Y=\operatorname{diag}(1,-1,1)Z$,
we rewrite (\ref{eq:(1)_3,(11)(1)}) as follows
\begin{gather*}
\frac{{\rm d}Z}{{\rm d}z}=A(z)Z, \\
A(z)=
\left(
\frac{A_0}{z^2}+\frac{A_1}{z}+A_2
\right)
\end{gather*}
where
\begin{gather*}
A_0= \begin{pmatrix}
0 & 1 & 0 \\
0 & 0 & 0 \\
0 & 0 & 0
\end{pmatrix},\qquad
A_1=
\begin{pmatrix}
p_2q_2 & -p_2 & p_1p_2 \\
0 & p_1q_1-p_2q_2-\theta^0_2 & 1 \\
-t & q_1 & -p_1q_1-\theta^0_1
\end{pmatrix},\\
A_2=
\begin{pmatrix}
0 & 0 & 0 \\
q_2 & -1 & p_1 \\
0 & 0 & 0
\end{pmatrix}.
\end{gather*}
We consider the singular point $z=0$.

Now we perform shearing transformations.
Let $S_1$ be the diagonal matrix $\operatorname{diag}(1, z^{1/3}, z^{2/3})$.
Then we have
\begin{gather*}
A^{S_1}(z)=
\frac{1}{z^{5/3}}
\begin{pmatrix} 0 & 1 & 0 \\ 0 & 0 & 0 \\ -t & 0 & 0 \end{pmatrix}
+\frac{1}{z^{4/3}}
 \begin{pmatrix} 0 & 0 & 0 \\ 0 & 0 & 0 \\ 0 & q_1 & 0 \end{pmatrix}\\
\hphantom{A^{S_1}(z)=}{} +\frac{1}{z}
\begin{pmatrix} p_2q_2 & 0 & 0 \\ 0 & p_1q_1-p_2q_2-\theta^0_2-1/3 \\ 0 & 0 & -p_1q_1-\theta^0_1-2/3 \end{pmatrix}+\cdots.
\end{gather*}
How to f\/ind the rational number $s$ of the shearing matrix is described in \cite{Wa}.
In this case, the dimension of the centralizer of the leading coef\/f\/icient matrix of $A^{S_1}(z)$ is less than that of~$A(z)$.
In fact, let $G_1$ be the following matrix
\begin{gather*}
G_1=
\begin{pmatrix}
0 & 1 & 0 \\
0 & 0 & 1 \\
-t & 0 & 0
\end{pmatrix}.
\end{gather*}
Then we have the following Jordan canonical form of the leading matrix of $A^{S_1}(z)$:
\begin{gather*}
{G_1}^{-1} \begin{pmatrix} 0 & 1 & 0 \\ 0 & 0 & 0 \\ -t & 0 & 0 \end{pmatrix}G_1=
\begin{pmatrix}
0 & 1 & 0 \\
0 & 0 & 1 \\
0 & 0 & 0
\end{pmatrix}.
\end{gather*}
Next, let $S_2=\operatorname{diag}(1, z^{1/3}, z^{2/3})$. Then we have
\begin{gather*}
A^{S_1G_1S_2}(z)=
\frac{1}{z^{4/3}}
 \begin{pmatrix} 0 & 1 & 0 \\ 0 & 0 & 1 \\ -t & 0 & 0 \end{pmatrix}\\
 \hphantom{A^{S_1G_1S_2}(z)=}{}
+\frac{1}{z}
 \begin{pmatrix} -p_1q_1-\theta^0_1-\frac{2}{3} & 0 & 0 \\ 0 & p_2q_2-\frac{1}{3} & 0 \\ 0 & 0 & p_1q_1-p_2q_2-\theta^0_2-1 \end{pmatrix}+\cdots.
\end{gather*}
Note that the leading matrix of $A^{S_1G_1S_2}(z)$ is diagonalizable. Indeed, the following matrix
\begin{gather*}
G_2=
\begin{pmatrix}
1 & 1 & 1 \\
-t^{1/3} & -\omega t^{1/3} & -\omega^2 t^{1/3}\\
t^{2/3} & \omega^2 t^{2/3} & \omega t^{2/3}
\end{pmatrix}
\end{gather*}
diagonalizes the leading matrix of $A^{S_1G_1S_2}(z)$:
\begin{gather*}
A^{S_1G_1S_2G_2}(z)=
\frac{1}{z^{4/3}}
\begin{pmatrix} -t^{1/3} & 0 & 0 \\ 0 & -\omega t^{1/3} & 0 \\ 0 & 0 & -\omega^2 t^{1/3} \end{pmatrix}\\
\hphantom{A^{S_1G_1S_2G_2}(z)=}{}+\frac{1}{z}
\begin{pmatrix} \theta^\infty_1/3-2/3 & * & * \\ * & \theta^\infty_1/3-2/3 & * \\ * & * & \theta^\infty_1/3-2/3 \end{pmatrix}+\cdots.
\end{gather*}

Since the leading matrix has three distinct eigenvalues, we can remove the of\/f-diagonal entries by virtue of Theorem~\ref{thm:block_diag}.
That is, there exits a matrix $P=P(z)$ whose entries are formal Laurent series in $z^{1/3}$ such that $A^{S_1G_1S_2G_2P}$ is diagonal:
\begin{gather*}
A^{S_1G_1S_2G_2P}(z)=
\frac{1}{z^{4/3}}
\begin{pmatrix} -t^{1/3} & 0 & 0 \\ 0 & -\omega t^{1/3} & 0 \\ 0 & 0 & -\omega^2 t^{1/3} \end{pmatrix}\\
\hphantom{A^{S_1G_1S_2G_2P}(z)=}{} +\frac{1}{z}
 \begin{pmatrix} \theta^\infty_1/3-2/3 & 0 & 0 \\ 0 & \theta^\infty_1/3-2/3 & 0 \\ 0 & 0 & \theta^\infty_1/3-2/3 \end{pmatrix}+\cdots.
\end{gather*}
Then, by virtue of Theorem~\ref{thm:first_red}, we can truncate $A^{S_1G_1S_2G_2P}(z)$ after the principal part
by a certain diagonal gauge transformation.
Furthermore, by the scalar gauge transformation by~$z^{-2/3}$, we can cancel the term $\frac{-2/3}{z}$.
In this way, we can obtain the HTL form of (\ref{eq:(1)_3,(11)(1)}) at~\smash{$x=\infty$}
\begin{gather}\label{eq:HTL(1)_3}
\frac{1}{z^{4/3}}
\begin{pmatrix} -t^{1/3} & 0 & 0 \\ 0 & -\omega t^{1/3} & 0 \\ 0 & 0 & -\omega^2 t^{1/3} \end{pmatrix}
+\frac{1}{z}
 \begin{pmatrix} \theta^\infty_1/3 & 0 & 0 \\ 0 & \theta^\infty_1/3 & 0 \\ 0 & 0 & \theta^\infty_1/3 \end{pmatrix}.
\end{gather}

Together with the HTL form at $x=0$ (this can be easily seen), we obtain the Riemann scheme of the system~(\ref{eq:(1)_3,(11)(1)}):
\begin{gather*}
\left(
\begin{matrix}
 x=0 & x=\infty \, \left( \frac13 \right) \\
\overbrace{\begin{matrix}
 1 & \theta^0_2 \\
 0 & \theta^0_1 \\
 0 & 0 \end{matrix}}
&
\overbrace{\begin{matrix}
 -t^{\frac13} & \theta^\infty_1/3 \\
 -\omega t^{\frac13} & \theta^\infty_1/3 \\
 -\omega^2 t^{\frac13} & \theta^\infty_1/3
 \end{matrix}}
\end{matrix}
\right).
\end{gather*}

We brief\/ly describe the feature of an HTL form (at a ramif\/ied irregular singularity) here. See~\cite{BV} for a precise statement.
Suppose the HTL form of $A(z) \in M_m(\mathbb{C}(\!(z)\!))$ has
\begin{gather}\label{eq:direct_summand}
\frac{T_{0}}{z^{\frac{b}{d}+1}}+\frac{T_{1}}{z^{\frac{b-1}{d}+1}}+\cdots+\frac{T_{b-1}}{z^{\frac{1}{d}+1}}+\frac{\Theta}{z}
\end{gather}
as its direct summand.
Then the HTL form of $A(z)$ also has
\begin{gather*}
\left\{
\frac{{\zeta_d}^{kb}T_0}{z^{\frac{b}{d}+1}}
+\frac{{\zeta_d}^{k(b-1)}T_1}{z^{\frac{b-1}{d}+1}}+\cdots
+\frac{{\zeta_d}^{k}T_{b-1}}{z^{\frac{1}{d}+1}}+\frac{\Theta}{z} \, \Bigg|\,
k=0,\ldots,d-1 \right\},
\end{gather*}
which is the orbit of (\ref{eq:direct_summand}) under the action $z^{\frac1d} \mapsto \zeta_d z^{\frac1d} \ (\zeta_d=e^{\frac{2\pi i}{d}})$
of a cyclic group, as its direct summands.

Let $S$ be the RSP corresponding to (\ref{eq:direct_summand}).
We denote the collection of $S$ and its $d-1$ copies by $S_d$.
Then an HTL form is generally represented as $S^1_{d_1} \ldots S^k_{d_k}$ where $S^1,\ldots,S^k$ are RSPs;
we call this the \textit{spectral type} of an HTL form.
The \textit{spectral type at a singular point} of a linear system is def\/ined as the spectral type of the HTL form at the point.
The tuple of the spectral types at all singular points is called the \textit{spectral type of the equation}.

\begin{Example}
The spectral types of
\begin{gather*}
\begin{matrix}
x=0 \ \left( \frac12 \right) \\
\overbrace{
\begin{matrix}
a & \alpha\\
a & \alpha \\
-a & \alpha \\
-a & \alpha
\end{matrix}}
\end{matrix}
\quad
\begin{matrix}
x=0 \ \left( \frac12 \right) \\
\overbrace{
\begin{matrix}
a & \alpha\\
a & \beta \\
-a & \alpha \\
-a & \beta
\end{matrix}}
\end{matrix}
\quad
\begin{matrix}
x=0 \ \left( \frac12 \right) \\
\overbrace{
\begin{matrix}
a & \alpha\\
-a & \alpha \\
b & \beta \\
-b & \beta
\end{matrix}}
\end{matrix}
\quad
\begin{matrix}
x=0 \ \left( \frac12 \right) \\
\overbrace{
\begin{matrix}
a & \alpha\\
-a & \alpha \\
0 & \beta \\
0 & \gamma
\end{matrix}}
\end{matrix}
\quad
\begin{matrix}
x=0 \ \left( \frac13 \right) \\
\overbrace{
\begin{matrix}
a & \alpha \\
\omega a & \alpha\\
\omega^2 a & \alpha \\
0 & \beta
\end{matrix}}
\end{matrix}
\quad
\begin{matrix}
 x=0 \ \left( \frac14 \right) \\
\overbrace{
\begin{matrix}
a & \alpha \\
\sqrt{-1}a & \alpha \\
-a & \alpha \\
-\sqrt{-1}a & \alpha
\end{matrix}}
\end{matrix}
\end{gather*}
are $(2)_2$, $(11)_2$, $(1)_2(1)_2$, $(1)_2 11$, $(1)_3 1$, and $(1)_4$ respectively.
\end{Example}

\subsection{Degeneration of HTL canonical forms}\label{sec:deg_HTL}

Suppose a system of linear equations has some parameter, say $\varepsilon$.
When we take the limit $\varepsilon \to 0$,
usually with a gauge transformation by some matrix which depends on $\varepsilon$ and is independent of~$x$,
an HTL form of the linear system at some singular point may change.
We call this situation a degeneration of an HTL form.

In the two dimensional case, the standard linear systems associated with classical Painlev\'e equations are $2 \times 2$,
and degenerations of HTL forms are realized by
the degenerations of the Jordan canonical forms
of the coef\/f\/icient matrices of the leading terms at irregular singular points.
However, when the rank of a linear system is greater than two,
a degeneration of an HTL form does not necessarily correspond to a degeneration of a Jordan canonical form.
We can see such degenerations in the degeneration schemes of the Sasano system and the Fuji--Suzuki system.
Here we take a sequence of degenerations $(11)(1), (11)(1) \to (1)_2 1, (11)(1) \to (1)_3, (11)(1)$ as an example.

Before going into the details, we look at the following simple example.
\begin{Example}\label{eg:2by2}
Let $A(x)$ be a $2 \times 2$ matrix of the form
\begin{gather}\label{eq:eg_2by2}
A(x)=\frac{A_0}{x^2}+\frac{A_1}{x}+\cdots \in M_2( \mathbb{C}(\!(x)\!) ),
\end{gather}
where
\begin{gather*}
A_0=
\begin{pmatrix}
0 & 1 \\
0 & 0
\end{pmatrix}, \qquad
A_k=\big( a^{(k)}_{ij} \big) \in M_2(\mathbb{C}).
\end{gather*}
Applying the gauge transformation by $S=\operatorname{diag}(1, x^{1/2})$ to (\ref{eq:eg_2by2}), we have
\begin{gather*}
A^S(x)=
\frac{1}{x^{3/2}}
\begin{pmatrix}
0 & 1 \\
a^{(1)}_{21} & 0
\end{pmatrix}+
\frac{1}{x}
\begin{pmatrix}
a^{(1)}_{11} & 0 \\
0 & a^{(1)}_{22}-\frac12
\end{pmatrix}+\cdots.
\end{gather*}
When $a^{(1)}_{21} \neq 0$, the leading coef\/f\/icient of $A^S(x)$ is a diagonalizable matrix with eigenvalues $\pm \sqrt{a^{(1)}_{21}}$.
This means that the singular point $x=0$ of the system of linear dif\/ferential equations
corresponding to~(\ref{eq:eg_2by2}) is an irregular singular point of Poincar\'e rank 1/2.
We can see that the HTL form at $x=0$ is
\begin{gather*}
\begin{matrix}
 x=0 \, \left( \frac12 \right) \\
\overbrace{
\begin{matrix}
\sqrt{a^{(1)}_{21}} & \operatorname{tr}(A_1) \\
-\sqrt{a^{(1)}_{21}} & \operatorname{tr}(A_1)
\end{matrix}}
\end{matrix}
\end{gather*}
by diagonalizing the leading matrix. When $a^{(1)}_{21}=0$, using $\tilde{S}=\operatorname{diag}(1, x)$ instead of the above $S$,
we have the dif\/ferent HTL form (we omit the details).
This example implies that if the leading matrix is a Jordan canonical form whose $(1,2)$-entry is non-zero,
then whether the $(2,1)$-entry of the subsequent matrix is zero or not is meaningful.
\end{Example}

Now we consider the degeneration $(11)(1), (11)(1) \to (1)_21, (11)(1)$.
The linear system of the spectral type $(11)(1), (11)(1)$ is given by~(\ref{eq:(11)(1),(11)(1)}).
Note that the leading matrix $A_0$ at the irregular singular point $x=\infty$ is diagonalizable.
The degeneration of the HTL form at $x=\infty$ 
is caused by the degeneration of the Jordan canonical form of $A_0$:
\begin{gather*}
\begin{pmatrix}
-t & & \\
 & 0 & \\
 & & 0
\end{pmatrix}
\to
\begin{pmatrix}
0 & 1 & \\
 & 0 & \\
 & & 0
\end{pmatrix},
\qquad
t \to 0.
\end{gather*}
In fact, changing the variables of~(\ref{eq:(11)(1),(11)(1)}) as in Appendix~\ref{sec:deg_FS}, we have the following new coef\/f\/icient matrices
\begin{gather*}
\tilde{A}_0 :=G^{-1}PUA_0U^{-1}P^{-1}G=
\begin{pmatrix}
0 & \tilde{t} & 0 \\
0 & 0 & 0 \\
0 & 0 & 0
\end{pmatrix}+O(\varepsilon), \\
\tilde{A}_1:=G^{-1}\hat{A}_1G=
\begin{pmatrix}
(\tilde{p}_1+\tilde{p}_2)\tilde{q}_1 & \tilde{q}_1 & \tilde{p}_2\tilde{q}_1 \\
1 & -\tilde{p}_1\tilde{q}_1-\tilde{p}_2\tilde{q}_2+\theta^0_1 & 1 \\
\tilde{p}_2(\tilde{q}_2-\tilde{q}_1)+\theta^0_2+\tilde{\theta}^\infty_2 & \tilde{q}_2-\tilde{q}_1 & \tilde{p}_2(\tilde{q}_2-\tilde{q}_1)+\theta^0_2
\end{pmatrix}
+O(\varepsilon), \\
\tilde{A}_2:=G^{-1}\hat{A}_2G=
\begin{pmatrix}
0 & 0 & 0 \\
\tilde{p}_1+\tilde{p}_2 & 1 & \tilde{p}_2 \\
0 & 0 & 0
\end{pmatrix},
\end{gather*}
where
\begin{gather*}
G=
\begin{pmatrix}
0 & 1 & 0 \\ -1/q_2 & 0 & -1/q_2 \\ \varepsilon q_1 & 0 & 0
\end{pmatrix}.
\end{gather*}
Here $\tilde{A}_0$ is diagonalizable provided that $\varepsilon$ is not equal to $0$,
and it degenerates to a nilpotent matrix when $\varepsilon$ tends to $0$.
We note that the $(2, 1)$-entry of $\lim_{\varepsilon \to 0}\tilde{A}_1$ is not equal to 0,
and thus we have the linear system (\ref{eq:(1)_2,(11)(1)}) of the spectral type $(1)_21, (11)(1)$ by $\varepsilon \to 0$.

\begin{Remark}
It is easy to obtain the HTL form at $x=\infty$, which corresponds to $(1)_21$, of~(\ref{eq:(1)_2,(11)(1)}).
In the same manner as Example~\ref{eg:2by2}, the shearing at $x=\infty$ can be done by the matrix $S=\operatorname{diag}(1, x^{-1/2}, 1)$.
\end{Remark}

On the other hand, the degeneration of a HTL form $(1)_2 1 \to (1)_3$ does not correspond to
the degeneration of a Jordan canonical form.
Let us see the degeneration $(1)_21, (11)(1) \to (1)_3, (11)(1)$.
In the course of the degeneration, the Jordan canonical form of the leading matrix stays unchanged.
Instead, the $(2,1)$-entry of the subsequent matrix goes to zero. 
In fact, changing the variables of (\ref{eq:(1)_2,(11)(1)}) as in Appendix~\ref{sec:deg_FS}, we have
\begin{gather*}
{G}^{-1}A_0G =
\begin{pmatrix}
0 & 1 & 0 \\
0 & 0 & 0 \\
0 & 0 & 0
\end{pmatrix}, \\
{G}^{-1}A_1G =
\begin{pmatrix}
-\tilde{p}_2\tilde{q}_2 & -\tilde{p}_2 & -\tilde{p}_1\tilde{p}_2 \\
0 & -\tilde{p}_1\tilde{q}_1+\tilde{p}_2\tilde{q}_2+\tilde{\theta}^0_2 & 1 \\
\tilde{t} & \tilde{q}_1 & \tilde{p}_1\tilde{q}_1+\tilde{\theta}^0_1
\end{pmatrix}\\
 \hphantom{{G}^{-1}A_1G =}{} +\varepsilon
\begin{pmatrix}
-\tilde{t}\tilde{p}_1\tilde{p}_2 & 0 & 0 \\
\tilde{t} & 0 & 0 \\
\tilde{t}(\tilde{p}_1\tilde{q}_1+\tilde{\theta}^0_1) & \tilde{t}\tilde{p}_2 & \tilde{t}\tilde{p}_1\tilde{p}_2
\end{pmatrix}
+O\big(\varepsilon^2\big),
\end{gather*}
where $G=\operatorname{diag}(t,1,1)$ (we have omitted the expression of $G^{-1}A_2G$).
The limit $\varepsilon \to 0$ causes the degeneration of the HTL form $(1)_2 1 \to (1)_3$.
In this way, we obtain the system~(\ref{eq:(1)_3,(11)(1)}).
The construction of the HTL form at $x=\infty$ of (\ref{eq:(1)_3,(11)(1)}) has been given in Section~\ref{sec:STram}.

We determine the possibility of degeneration as follows.
For example, $(1)_2 1$ is a direct sum of two direct summands (see the Riemann scheme of $(1)_2 1, (11)(1)$)
\begin{gather*}
\frac{1}{z^{3/2}}
\begin{pmatrix}
\sqrt{t} & 0 \\
0 & -\sqrt{t}
\end{pmatrix}
+
\frac{1}{z}
\begin{pmatrix}
\theta^\infty_1/2 & 0 \\
0 & \theta^\infty_1/2
\end{pmatrix}
\end{gather*}
and $\theta^\infty_2/z$.
From this, we expect that it is possible to take a limit $t \to 0$ and
indeed this corresponds to the degeneration $(1)_2 1 \to (1)_3$.
On the other hand, $(1)_3$ itself consists of a single Galois orbit, see (\ref{eq:HTL(1)_3}).
Thus we conclude that $(1)_3$ does not admit degeneration.

Let us make a remark on degeneration of HTL forms.
We do not consider the degeneration of the HTL form $(2)(1)$
since the degeneration of $(2)(1)$ does not preserve the number of accessory parameters.
\begin{Remark}
Accessory parameters of a linear system are free parameters remaining in the linear system when the Riemann scheme is f\/ixed.
The number of accessory parameters of a~linear system coincides with the dimension of the phase space of the corresponding Painlev\'e-type equation.
\end{Remark}
To see this, for example, suppose that the HTL form $(2)(1)$ of the $(2)(1), 111, 111$-system degenerates.
Then the degenerated linear system has the following form:
\begin{gather}\label{eq:deg_(11)(1),111,111}
\frac{{\rm d}Y}{{\rm d}x}=
\left\{
\begin{pmatrix}
0 & 1 & 0 \\
0 & 0 & 0 \\
0 & 0 & 0
\end{pmatrix}
+\frac{A_0}{x}+\frac{A_1}{x-1}
\right\}Y, \qquad
A_* \sim \operatorname{diag}(0, \theta^*_1, \theta^*_2), \qquad *=0,1.
\end{gather}
By a direct calculation, we f\/ind that the number of accessory parameters of the above system is six.
In fact, there is a Fuchsian equation with spectral type $21,111,111,111$, which has six accessory parameters.
The system (\ref{eq:deg_(11)(1),111,111}) turns out to be the degenerated system of the $21,111,111,111$-system.
The same argument applies to the degeneration of $((11))((1))$ of the $((11))((1)), 111$-system.

\section{Lax pairs of degenerate FS and Garnier systems}\label{sec:Lax}
The Garnier systems and the Fuji--Suzuki systems are non-linear dif\/ferential equations, which are regarded as
generalizations of the Painlev\'e equations.

The Garnier system in $N$ variables was derived as the isomonodromic deformation equation of a second order Fuchsian equation
with $N+3$ singular points~\cite{G}.
The Fuji--Suzuki systems were originally derived from the Drinfeld--Sokolov hierarchy by similarity reductions~\cite{FS1}.

In Section~\ref{sec:Gar_Lax}, we present the Riemann schemes, Lax pairs, and corresponding Hamiltonians
of degenerate Garnier systems. In Section~\ref{sec:FS_Lax}, we present similar data of degenerate Fuji--Suzuki systems.

\subsection{degenerate Garnier systems}\label{sec:Gar_Lax}
The Garnier systems were originally derived from a second order single Fuchsian equation with $N+3$ singular points.
However, unlike the original study by Garnier, we adopt f\/irst order systems concerning linear equations.
In~\cite{Sak2}, the Garnier system in $N$ variables was derived from the f\/irst order $2 \times 2$ system of the form
\begin{gather}\label{eq:Lax_Gar}
\frac{{\rm d}Y}{{\rm d}x}=
\left(
\frac{A_0}{x}+\frac{A_1}{x-1}+\sum_{j=1}^N \frac{A_{t_j}}{x-t_j}
\right)Y.
\end{gather}
When $N$ equals 2, the Painlev\'e-type equation corresponding to (\ref{eq:Lax_Gar}) has a four-dimensional phase space. In~\cite{KNS}, conf\/luences from this linear system were considered.

The degeneration of the Garnier system in two variables was considered by Kimura~\cite{Km}. He treated mainly the conf\/luence of singular points of associated linear equations, and he obtained the degenerated Garnier systems with the singularity pattern $2+1+1+1$, $3+1+1$, $2+2+1$, $3+2$, $4+1$, $5$, and $9/2$.
Kawamuko~\cite{Kw} further considered the degeneration of HTL canonical forms and obtained eight degenerate Garnier systems.

In this subsection, we give the Riemann schemes, Lax pairs, and Hamiltonians for degenerate Garnier systems associated with ramif\/ied linear equations.
Although all the Hamiltonians in this subsection are equivalent to those in \cite{Km, Kw},
we recalculated them.
The following is the list of Hamiltonians for the degenerate Garnier systems associated with ramif\/ied linear equations:
\begin{gather*}
t_1H_{\mathrm{Gar},t_1}^{\frac{3}{2}+1+1+1}\left({\alpha,\beta\atop \gamma};{t_1 \atop t_2};{q_1,p_1 \atop q_2,p_2}\right)
=
t_1H_{\mathrm{III}(D_6)}\left({-\alpha,\gamma-\alpha};t_1;q_1,p_1\right)
+q_1(q_1p_1-\alpha)p_2\\
\hphantom{t_1H_{\mathrm{Gar},t_1}^{\frac{3}{2}+1+1+1}\left({\alpha,\beta\atop \gamma};{t_1 \atop t_2};{q_1,p_1 \atop q_2,p_2}\right)=}{}
+\frac{t_1}{t_1-t_2}(p_1(q_1-q_2)-\alpha)(p_2(q_2-q_1)-\beta), \\
t_2H_{\mathrm{Gar},t_2}^{\frac{3}{2}+1+1+1}\left({\alpha,\beta \atop \gamma};{t_1 \atop t_2};{q_1,p_1 \atop q_2,p_2}\right)
=
t_2H_{\mathrm{III}(D_6)}\left({-\beta,\gamma-\beta};t_2;q_2,p_2\right)
+q_2(q_2p_2-\beta)p_1 \\
\hphantom{t_2H_{\mathrm{Gar},t_2}^{\frac{3}{2}+1+1+1}\left({\alpha,\beta \atop \gamma};{t_1 \atop t_2};{q_1,p_1 \atop q_2,p_2}\right)
=}{}+\frac{t_2}{t_2-t_1}(p_1(q_1-q_2)-\alpha)(p_2(q_2-q_1)-\beta),
\\
H_{\mathrm{Gar},t_1}^{\frac{5}{2}+1+1}\left({\alpha,\beta};{t_1 \atop t_2};{q_1,p_1 \atop q_2,p_2}\right)
=
H_{\mathrm{II}}\left(-\alpha;t_1;q_1,p_1\right)
+p_1p_2\\
\hphantom{H_{\mathrm{Gar},t_1}^{\frac{5}{2}+1+1}\left({\alpha,\beta};{t_1 \atop t_2};{q_1,p_1 \atop q_2,p_2}\right)
=}{} +\frac{1}{t_1-t_2}(p_1(q_1-q_2)-\alpha)(p_2(q_2-q_1)-\beta),\\
H_{\mathrm{Gar},t_2}^{\frac{5}{2}+1+1}\left({\alpha,\beta};{t_1 \atop t_2};{q_1,p_1 \atop q_2,p_2}\right)
=
H_{\mathrm{II}}\left(-\beta;t_2;q_2,p_2\right)
+p_1p_2\\
\hphantom{H_{\mathrm{Gar},t_2}^{\frac{5}{2}+1+1}\left({\alpha,\beta};{t_1 \atop t_2};{q_1,p_1 \atop q_2,p_2}\right)
=}{}
+\frac{1}{t_2-t_1}(p_1(q_1-q_2)-\alpha)(p_2(q_2-q_1)-\beta),
\\
t_1H_{\mathrm{Gar},t_1}^{2+\frac{3}{2}+1}\left({\alpha, \beta};{t_1 \atop t_2};{q_1,p_1 \atop q_2,p_2}\right)
=
t_1H_{\mathrm{III}(D_6)}\left({-\alpha, -\beta-\alpha}; t_1; q_1,p_1\right)
+t_2p_1p_2\nonumber\\
\hphantom{t_1H_{\mathrm{Gar},t_1}^{2+\frac{3}{2}+1}\left({\alpha, \beta};{t_1 \atop t_2};{q_1,p_1 \atop q_2,p_2}\right)
=}{} -p_2q_2(2p_1q_1-\alpha)+\frac{t_2}{t_1}p_1q_1-\frac{q_1q_2}{t_1}(p_1q_1-\alpha),\\
t_2H_{\mathrm{Gar},t_2}^{2+\frac{3}{2}+1}\left({\alpha, \beta};{t_1 \atop t_2};{q_1,p_1 \atop q_2,p_2}\right)
=
t_2H_{\mathrm{III}(D_7)}\left(\beta+1; t_2; q_2,p_2\right)
-t_2p_1p_2\nonumber\\
\hphantom{t_2H_{\mathrm{Gar},t_2}^{2+\frac{3}{2}+1}\left({\alpha, \beta};{t_1 \atop t_2};{q_1,p_1 \atop q_2,p_2}\right)
=}{} -\frac{t_2}{t_1}p_1q_1+\frac{q_1q_2}{t_1}(p_1q_1-\alpha),
\\
H^{\frac72+1}_{\mathrm{Gar},t_1}\left({\alpha};{t_1 \atop t_2};{q_1,p_1 \atop q_2,p_2}\right)
=
H_{\mathrm{I}}\left(t_1;q_1,p_1\right)+p_2\big(2q_1-{q_2}^2-t_2\big)+\alpha q_2,\\
H^{\frac72+1}_{\mathrm{Gar},t_2}\left({\alpha};{t_1 \atop t_2};{q_1,p_1 \atop q_2,p_2}\right)
=
{p_2}^2-t_2p_2{q_2}^2-{t_2}^2p_2+\alpha t_2q_2\\
\hphantom{H^{\frac72+1}_{\mathrm{Gar},t_2}\left({\alpha};{t_1 \atop t_2};{q_1,p_1 \atop q_2,p_2}\right)
=}{} +2p_1p_2q_2-q_1q_2(p_2q_2-\alpha)-p_2q_1(q_1-t_2)-t_1p_2-\alpha p_1,\nonumber
\\
t_1H^{3+\frac32}_{\mathrm{Gar},t_1}\left({\alpha};{t_1 \atop t_2};{q_1,p_1 \atop q_2,p_2}\right)
=
t_1H_{\mathrm{III}(D_7)}\left({\alpha};t_1;q_1,p_1\right)
+p_2q_2(p_2q_2+\alpha)-q_2\nonumber\\
\hphantom{t_1H^{3+\frac32}_{\mathrm{Gar},t_1}\left({\alpha};{t_1 \atop t_2};{q_1,p_1 \atop q_2,p_2}\right)
=}{} +q_1(p_2q_1+2p_1p_2q_2-t_2-1),\\
H^{3+\frac32}_{\mathrm{Gar},t_2}\left({\alpha};{t_1 \atop t_2};{q_1,p_1 \atop q_2,p_2}\right)
=
p_2q_1(p_1q_1+2p_2q_2+\alpha-1)+p_1q_2-q_1+t_1p_2-t_2p_2q_2,
\\
H^{\frac52+2}_{\mathrm{Gar},t_1}\left({\alpha};{t_1 \atop t_2};{q_1,p_1 \atop q_2,p_2}\right)
=
H_{\mathrm{II}}\left({\alpha};t_1;q_1,p_1\right)
-2p_2q_2q_1-t_2p_2-q_2,\\
t_2H^{\frac52+2}_{\mathrm{Gar},t_2}\left({\alpha};{t_1 \atop t_2};{q_1,p_1 \atop q_2,p_2}\right)
={p_2}^2{q_2}^2+\alpha p_2q_2+t_2p_2\big(p_1-{q_1}^2-t_1\big)-p_1q_2-t_2q_1,
\\
t_1H_{\mathrm{Gar},t_1}^{\frac{3}{2}+\frac{3}{2}+1}\left({\alpha};{t_1 \atop t_2};{q_1,p_1 \atop q_2,p_2}\right)
=
t_1H_{\mathrm{III}(D_6)}\left({\alpha, \alpha}; t_1; q_1,p_1\right)
+\frac{q_1q_2}{t_1}(p_1q_1+\alpha)-2p_1q_1p_2q_2\nonumber\\
\hphantom{t_1H_{\mathrm{Gar},t_1}^{\frac{3}{2}+\frac{3}{2}+1}\left({\alpha};{t_1 \atop t_2};{q_1,p_1 \atop q_2,p_2}\right)
=}{}-\alpha p_2q_2-\frac{t_2p_1}{q_2},\\
t_2H_{\mathrm{Gar},t_2}^{\frac{3}{2}+\frac{3}{2}+1}\left({\alpha};{t_1 \atop t_2};{q_1,p_1 \atop q_2,p_2}\right)
=
t_2H_{\mathrm{III}(D_8)}\left(t_2; q_2, p_2\right)
-\frac{q_1q_2}{t_1}(p_1q_1+\alpha)+\frac{t_2p_1}{q_2},
\\
H_{\mathrm{Gar},t_1}^{\frac92}\left({t_1 \atop t_2};{q_1,p_1 \atop q_2,p_2}\right)=
-{p_1}^3p_2+{p_1}^2{q_2}^2+t_1{p_1}^3-2p_1q_1q_2-2p_1{p_2}^2+p_2{q_2}^2\nonumber\\
\hphantom{H_{\mathrm{Gar},t_1}^{\frac92}\left({t_1 \atop t_2};{q_1,p_1 \atop q_2,p_2}\right)=}{} +t_1p_1p_2+{q_1}^2-t_1{q_2}^2+{t_1}^2p_1+t_2p_2,\\
H_{\mathrm{Gar},t_2}^{\frac92}\left({t_1 \atop t_2};{q_1,p_1 \atop q_2,p_2}\right)=
{p_1}^4+3{p_1}^2p_2+p_1{q_2}^2-2q_1q_2+{p_2}^2-t_2p_1+t_1p_2,
\\
H^{\frac52+\frac32}_{\mathrm{Gar},t_1}\left({t_1 \atop t_2};{q_1,p_1 \atop q_2,p_2}\right)
=
H_{\mathrm{II}}\left({0};t_1;q_1,p_1\right)
-2p_2q_2q_1-q_2-\frac{t_2}{q_2},\\
t_2H^{\frac52+\frac32}_{\mathrm{Gar},t_2}\left({t_1 \atop t_2};{q_1,p_1 \atop q_2,p_2}\right)
={p_2}^2{q_2}^2-p_1q_2+\frac{t_2}{q_2}\big(p_1-{q_1}^2-t_1\big).
\end{gather*}

\subsubsection*{Singularity pattern $\boldsymbol{\frac32+1+1+1}$}
The Riemann scheme is given by
\begin{gather*}
\left(\begin{matrix}
 x=0 & x=t_1 & x=t_2 & x=\infty \, \left( \frac12 \right)\\
\begin{matrix}0 \\ \theta^0 \end{matrix} &
\begin{matrix}0 \\ \theta^{t_1} \end{matrix} &
\begin{matrix}0 \\ \theta^{t_2} \end{matrix} &
\overbrace{\begin{matrix}
 1 & \theta^\infty_1/2 \\
 -1 & \theta^\infty_1/2
 \end{matrix}}\
 \end{matrix}\right) ,
\end{gather*}
and the Fuchs--Hukuhara relation is written as $\theta^0+\theta^{t_1}+\theta^{t_2}+\theta^\infty_1=0$.
The Lax pair is expressed as
\begin{gather*}
\frac{\partial Y}{\partial x}=\left(\frac{A_0}{x}+\frac{A_{t_1}}{x-t_1}+\frac{A_{t_2}}{x-t_2}+N \right)Y ,\\
\frac{\partial Y}{\partial t_1}=\left( N_1-\frac{A_{t_1}}{x-t_1} \right)Y ,\qquad
\frac{\partial Y}{\partial t_2}=\left( N_2-\frac{A_{t_2}}{x-t_2} \right)Y .
\end{gather*}
Here
\begin{gather*}
A_0=
\begin{pmatrix}
0 \\
1
\end{pmatrix}
\begin{pmatrix}
1-p_1-p_2 & \theta^0
\end{pmatrix},\qquad
A_{t_i}=
\begin{pmatrix}
q_i \\
1
\end{pmatrix}
\begin{pmatrix}
p_i & \theta^{t_i}-p_i q_i
\end{pmatrix},\\
N=
 \begin{pmatrix}
 0 & 1 \\
 0 & 0
 \end{pmatrix},\qquad
N_i=
\frac{q_i(p_i q_i-\theta^{t_i})}{t_i}N, \qquad i=1,2.
\end{gather*}

The Hamiltonians are given by
\begin{gather*}
t_1H_{\mathrm{Gar},t_1}^{\frac{3}{2}+1+1+1}\left({\theta^{t_1}, \theta^{t_2} \atop -\theta^0};{t_1 \atop t_2};{q_1,p_1 \atop q_2,p_2}\right)
=
t_1H_{\mathrm{III}(D_6)}\left({-\theta^{t_1}, -\theta^0-\theta^{t_1}};t_1;q_1,p_1\right)\!+\!q_1(q_1p_1-\theta^{t_1})p_2\!\! \nonumber \\
\hphantom{t_1H_{\mathrm{Gar},t_1}^{\frac{3}{2}+1+1+1}\left({\theta^{t_1}, \theta^{t_2} \atop -\theta^0};{t_1 \atop t_2};{q_1,p_1 \atop q_2,p_2}\right)
=}{}
+\frac{t_1}{t_1-t_2}(p_1(q_1-q_2)-\theta^{t_1})(p_2(q_2-q_1)-\theta^{t_2}),\nonumber\\
t_2H_{\mathrm{Gar},t_2}^{\frac{3}{2}+1+1+1}\left({\theta^{t_1}, \theta^{t_2} \atop -\theta^0};{t_1 \atop t_2};{q_1,p_1 \atop q_2,p_2}\right)
=
t_2H_{\mathrm{III}(D_6)}\left({-\theta^{t_2}, -\theta^0-\theta^{t_2}};t_2;q_2,p_2\right)\!+\!q_2(q_2p_2-\theta^{t_2})p_1\!\! \nonumber \\
\hphantom{t_2H_{\mathrm{Gar},t_2}^{\frac{3}{2}+1+1+1}\left({\theta^{t_1}, \theta^{t_2} \atop -\theta^0};{t_1 \atop t_2};{q_1,p_1 \atop q_2,p_2}\right)
=}{}
+\frac{t_2}{t_2-t_1}(p_1(q_1-q_2)-\theta^{t_1})(p_2(q_2-q_1)-\theta^{t_2}).\nonumber
\end{gather*}

\subsubsection*{Singularity pattern $\boldsymbol{\frac52+1+1}$}
The Riemann scheme is given by
\begin{gather*}
\left(
\begin{matrix}
x=t_1 & x=t_2 & x=\infty \ \left( \frac12 \right)\\
 \begin{matrix}
	0 \\
	\theta^{t_1}
	\end{matrix}
 &\begin{matrix}
	0 \\
	\theta^{t_2}
	\end{matrix}
&\overbrace{
	\begin{matrix}
\sqrt{-1} & 0 & 0 & \theta^\infty_1/2 \\
-\sqrt{-1} & 0 & 0 & \theta^\infty_1/2
	\end{matrix}}
\end{matrix}
\right),
\end{gather*}
and the Fuchs--Hukuhara relation is written as $\theta^{t_1}+\theta^{t_2}+\theta^\infty_1=0$. The Lax pair is given as
\begin{gather*}
\frac{\partial Y}{\partial x}=\left( \frac{A_{t_1}}{x-t_1}+\frac{A_{t_2}}{x-t_2}+A_{\infty\,1}+N x  \right)Y ,\\
\frac{\partial Y}{\partial t_1}=\left(N_1-\frac{A_{t_1}}{x-t_1}\right)Y ,\qquad
\frac{\partial Y}{\partial t_2}=\left(N_2-\frac{A_{t_2}}{x-t_2}\right)Y .
\end{gather*}
Here
\begin{gather*}
A_{t_i} =
\begin{pmatrix}
-q_i \\
1
\end{pmatrix}
\begin{pmatrix}
-p_i & -p_i q_i+\theta^{t_i}
\end{pmatrix},\qquad
N=
\begin{pmatrix}
0 & 1 \\
0 & 0
\end{pmatrix},\\
A_{\infty \, 1} =
\begin{pmatrix}
0 & -p_1-p_2 \\
-1 & 0
\end{pmatrix},\qquad
N_i=-p_i N , \qquad i=1,2.
\end{gather*}

The Hamiltonians are written as
\begin{gather*}
H_{\mathrm{Gar},t_1}^{\frac{5}{2}+1+1}
\left({\theta^{t_1}, \theta^{t_2}};{t_1 \atop t_2};
{q_1,p_1 \atop q_2,p_2}\right)\\
\qquad {} =H_{\mathrm{II}}\left(-\theta^{t_1};t_1;q_1,p_1\right)
+p_1p_2+\frac{1}{t_1-t_2}(p_1(q_1-q_2)-\theta^{t_1})(p_2(q_2-q_1)-\theta^{t_2}),\nonumber\\
H_{\mathrm{Gar},t_1}^{\frac{5}{2}+1+1} \left({\theta^{t_1}, \theta^{t_2}};{t_1 \atop t_2};
{q_1,p_1 \atop q_2,p_2}\right)\\
\qquad {}=
H_{\mathrm{II}}\left(-\theta^{t_2};t_2;q_2,p_2\right)
+p_1p_2+\frac{1}{t_2-t_1}(p_1(q_1-q_2)-\theta^{t_1})(p_2(q_2-q_1)-\theta^{t_2}).\nonumber
\end{gather*}

\subsubsection*{Singularity pattern $\boldsymbol{2+\frac32+1}$}
The Riemann scheme is given by
\begin{gather*}
\left(
\begin{matrix}
x=0 & x=t_1 & x=\infty \, \left( \frac12 \right)\\
\overbrace{
	\begin{matrix}
	0 & 0\\
	-t_2 & \theta^0
	\end{matrix}}
	&\begin{matrix}
	0\\
	\theta^{t_1}
	\end{matrix}
&\overbrace{
	\begin{matrix}
	1 & \theta^\infty_1/2 \\
	-1 & \theta^\infty_1/2
	\end{matrix}}
\end{matrix}
\right),
\end{gather*}
and the Fuchs--Hukuhara relation is written as $\theta^0+\theta^{t_1}+\theta^\infty_1=0$. The Lax pair is expressed as
\begin{gather*}
\frac{\partial Y}{\partial x}=\left( \frac{A_0^{(1)}}{x^2}+\frac{A_0^{(0)}}{x}+\frac{A_{t_1}}{x-t_1}+N \right)Y ,\\
\frac{\partial Y}{\partial t_1}=\left( N_1-\frac{A_{t_1}}{x-t_1} \right)Y ,\qquad
\frac{\partial Y}{\partial t_2}=\left( N_2-\frac{\frac{1}{t_2}A_0^{(1)}}{x} \right)Y .
\end{gather*}
Here
\begin{gather*}
A_0^{(1)}=
\begin{pmatrix}
0 & 0 \\
q_2 & -t_2
\end{pmatrix},\qquad
A_0^{(0)}=
\begin{pmatrix}
-p_2q_2 & t_2p_2 \\
1-p_1 & p_2q_2+\theta^0
\end{pmatrix},\qquad
N=
 \begin{pmatrix}
 0 & 1 \\
 0 & 0
 \end{pmatrix},\\
A_{t_1}=
\begin{pmatrix}
q_1 \\
1
\end{pmatrix}
\begin{pmatrix}
p_1 & \theta^{t_1}-p_1q_1
\end{pmatrix},\qquad
N_1=
\frac{q_1(p_1q_1-\theta^{t_1})}{t_1}N,\qquad
N_2=-p_2N.
\end{gather*}

The Hamiltonians are given by
\begin{gather*}
t_1H_{\mathrm{Gar},t_1}^{2+\frac{3}{2}+1}\left({\theta^{t_1}, \theta^0};{t_1 \atop t_2};{q_1,p_1 \atop q_2,p_2}\right)
=
t_1H_{\mathrm{III}(D_6)}\left({-\theta^{t_1}, -\theta^0-\theta^{t_1}}; t_1; q_1,p_1\right)
+t_2p_1p_2\nonumber\\
\qquad\quad{} -p_2q_2\big(2p_1q_1-\theta^{t_1}\big)+\frac{t_2}{t_1}p_1q_1-\frac{q_1q_2}{t_1}\big(p_1q_1-\theta^{t_1}\big), \nonumber \\
t_2H_{\mathrm{Gar},t_2}^{2+\frac{3}{2}+1}\left({\theta^{t_1}, \theta^0};{t_1 \atop t_2};{q_1,p_1 \atop q_2,p_2}\right)\\
\qquad{} =
t_2H_{\mathrm{III}(D_7)}\left(\theta^0+1; t_2; q_2,p_2\right)
-t_2p_1p_2-\frac{t_2}{t_1}p_1q_1+\frac{q_1q_2}{t_1}(p_1q_1-\theta^{t_1}). \nonumber
\end{gather*}

\subsubsection*{Singularity pattern $\boldsymbol{\frac72+1}$}
The Riemann scheme is given by
\begin{gather*}
\left(
\begin{matrix}
x=0 & x=\infty	\,\left( \frac12 \right)\\
 \begin{matrix}
	0 \\
	\theta^0
	\end{matrix}
&\overbrace{
	\begin{matrix}
1 &	 0 & -\frac{3t_2}{2} & 0 & \frac{t_1}{2}+\frac{3{t_2}^2}{8} & \theta^\infty_1/2 \\
-1 &	 0 & \frac{3t_2}{2} & 0 & -\frac{t_1}{2}-\frac{3{t_2}^2}{8} & \theta^\infty_1/2
	\end{matrix}}
\end{matrix}
\right) ,
\end{gather*}
and the Fuchs--Hukuhara relation is written as $\theta^0+\theta^\infty_1=0$. The Lax pair is expressed as
\begin{gather*}
\frac{\partial Y}{\partial x}=\left( \frac{A_3}{x}+A_2+A_1x+A_0x^2  \right)Y ,\\
\frac{\partial Y}{\partial t_1} =(A_0 x+B_{10})Y ,\qquad \frac{\partial Y}{\partial t_2}=
\left( -A_0 x^2+B_{21}x+B_{20} \right)Y .
\end{gather*}
Here
\begin{gather*}
A_0=
\begin{pmatrix}
0 & 1 \\
0 & 0
\end{pmatrix},\qquad
A_1=
\begin{pmatrix}
0 & q_1-2t_2 \\
1 & 0
\end{pmatrix},\\
A_2=
\begin{pmatrix}
-p_1 & -p_2+{q_1}^2-t_2q_1+t_1+{t_2}^2\\
-q_1-t_2 & p_1
\end{pmatrix},\qquad
A_3=
\begin{pmatrix}
q_2 \\
1
\end{pmatrix}
\begin{pmatrix}
p_2 & -p_2q_2+\theta^0
\end{pmatrix},\\
B_{10}=
\begin{pmatrix}
0 & 2q_1-t_2 \\
1 & 0
\end{pmatrix},\qquad
B_{21}=
\begin{pmatrix}
0 & 2t_2-q_1 \\
-1 & 0
\end{pmatrix},\\
B_{20} =
\begin{pmatrix}
p_1 & 2p_2-{q_1}^2+t_2q_1-t_1-{t_2}^2\\
q_1+t_2 & -p_1
\end{pmatrix}.
\end{gather*}

The Hamiltonians are given by
\begin{gather*}
H^{\frac72+1}_{\mathrm{Gar},t_1} \left({\theta^0};{t_1 \atop t_2}; {q_1,p_1 \atop q_2,p_2}\right)
= H_{\mathrm{I}}\left(t_1;q_1,p_1\right)+p_2\big(2q_1-{q_2}^2-t_2\big)+\theta^0 q_2,\\
H^{\frac72+1}_{\mathrm{Gar},t_2} \left({\theta^0};{t_1 \atop t_2};
{q_1,p_1 \atop q_2,p_2}\right) = {p_2}^2-t_2p_2{q_2}^2-{t_2}^2p_2+\theta^0 t_2q_2+ 2p_1p_2q_2-q_1q_2\big(p_2q_2-\theta^0\big) \nonumber \\
\qquad{}-p_2q_1(q_1-t_2)-t_1p_2-\theta^0 p_1.\nonumber
\end{gather*}

\subsubsection*{Singularity pattern $\boldsymbol{3+\frac32}$}
The Riemann scheme is given by
\begin{gather*}
\left(
\begin{matrix}
x=0\,\left( \frac12 \right) & x=\infty	\\
\overbrace{
 \begin{matrix}
	\sqrt{t_1} & 0 \\
	-\sqrt{t_1} & 0
	\end{matrix}}
&\overbrace{
	\begin{matrix}
	0 & 0 & \theta^\infty_1 \\
	1 & -t_2 & \theta^\infty_2
	\end{matrix}}
\end{matrix}
\right) ,
\end{gather*}
and the Fuchs--Hukuhara relation is written as $\theta^\infty_1+\theta^\infty_2=0$. The Lax pair is expressed as
\begin{gather*}
\frac{\partial Y}{\partial x}=\left(\frac{A_3}{x^2}+\frac{A_2}{x}+A_1+A_0 x  \right)Y ,\qquad
\frac{\partial Y}{\partial t_1}=-\frac{A_3}{t_1x}Y, \qquad \frac{\partial Y}{\partial t_2}=(E_2 x+B_1)Y.
\end{gather*}
Here
\begin{gather*}
A_\xi =
\begin{pmatrix}
1 & \\
 & u
\end{pmatrix}^{-1}
\hat{A}_\xi
\begin{pmatrix}
1 & \\
 & u
\end{pmatrix},\qquad
E_2=\operatorname{diag}(0,\, 1),\\
\hat{A}_3=t_1
\begin{pmatrix}
1 \\
p_2
\end{pmatrix}
\begin{pmatrix}
-p_2 & 1
\end{pmatrix},\qquad
\hat{A}_0=-E_2,\qquad
\hat{A}_1=
\begin{pmatrix}
0 & q_1 \\
p_2(p_2q_1-t_2)+p_1 & t_2
\end{pmatrix},
\\
\hat{A}_2=
\begin{pmatrix}
-p_2q_1(p_2q_1-t_2)-p_1q_1-\theta^\infty_1 & q_1(p_2q_1-t_2)-q_2\\
(\hat{A}_2)_{21} & p_2q_1(p_2q_1-t_2)+p_1q_1+\theta^\infty_1 \\
\end{pmatrix},\\
B_1 =
\begin{pmatrix}
0 & (-A_1)_{12} \\
(-A_1)_{21} & 0
\end{pmatrix},
\end{gather*}
where
\begin{gather*}
(\hat{A}_2)_{21}=-{p_2}^2q_1(p_2q_1-t_2)-p_2\big(2p_1q_1+p_2q_2+2\theta^\infty_1\big)+1.
\end{gather*}

The Hamiltonians are given by
\begin{gather*}
t_1H^{3+\frac32}_{\mathrm{Gar},t_1}\left({2\theta^\infty_1};{t_1 \atop t_2};{q_1,p_1 \atop q_2,p_2}\right)
=
t_1H_{\mathrm{III}(D_7)}\left({2\theta^\infty_1};t_1;q_1,p_1\right)
+p_2q_2\big(p_2q_2+2\theta^\infty_1\big)-q_2\nonumber\\
\qquad{}+q_1 (p_2q_1+2p_1p_2q_2-t_2-1), \nonumber \\
H^{3+\frac32}_{\mathrm{Gar},t_2}\left({2\theta^\infty_1};{t_1 \atop t_2};{q_1,p_1 \atop q_2,p_2}\right)
=
p_2q_1\big(p_1q_1+2p_2q_2+2\theta^\infty_1-1\big)+p_1q_2-q_1+t_1p_2-t_2p_2q_2. \nonumber
\end{gather*}
The gauge parameter $u$ satisf\/ies
\begin{gather*}
\frac{1}{u}\frac{\partial u}{\partial t_1}=-\frac{2}{t_1}\big(p_1q_1+p_2q_2+\theta^\infty_1\big),\qquad
\frac{1}{u}\frac{\partial u}{\partial t_2}=-2p_2q_1.
\end{gather*}

\subsubsection*{Singularity pattern $\boldsymbol{\frac52+2}$}
The Riemann scheme is given by
\begin{gather*}
\left(
\begin{matrix}
x=0 & x=\infty \, \left( \frac12 \right)	\\
\overbrace{
 \begin{matrix}
	0 & 0 \\
	t_2 & \theta^0
	\end{matrix}}
&\overbrace{
	\begin{matrix}
	1 & 0 & -t_1/2 & \theta^\infty_1/2 \\
	-1 & 0 & t_1/2 & \theta^\infty_1/2
	\end{matrix}}
\end{matrix}
\right) ,
\end{gather*}
and the Fuchs--Hukuhara relation is written as $\theta^0+\theta^\infty_1=0$.
The Lax pair is expressed as
\begin{gather*}
\frac{\partial Y}{\partial x}=\left(\frac{A_3}{x^2}+\frac{A_2}{x}+A_1+A_0 x  \right)Y ,\\
\frac{\partial Y}{\partial t_1}=\left( x B_{11}+B_{10}\right)Y ,\qquad
\frac{\partial Y}{\partial t_2}=\left( -\frac{A_3}{t_2x}+B_{20}\right)Y ,
\end{gather*}
where
\begin{gather*}
A_0=
\begin{pmatrix}
0 & 1 \\
0 & 0
\end{pmatrix},\qquad
A_1=
\begin{pmatrix}
q_1 & p_1-{q_1}^2-t_1\\
1 & -q_1 \\
\end{pmatrix},\\
A_2=
\begin{pmatrix}
p_2q_2 & q_2 \\
-p_1 & -p_2q_2-\theta^\infty_1
\end{pmatrix},\qquad
A_3=
\begin{pmatrix}
0 & 0 \\
t_2p_2 & t_2
\end{pmatrix},\\
B_{11}=
\begin{pmatrix}
0 & -1 \\
0 & 0
\end{pmatrix},\qquad
B_{10}=
\begin{pmatrix}
-q_1 & 0 \\
-1 & q_1
\end{pmatrix},\qquad
B_{20}=
\begin{pmatrix}
0 & -q_2/t_2 \\
0 & 0
\end{pmatrix}.
\end{gather*}

The Hamiltonians are given by
\begin{gather*}
H^{\frac52+2}_{\mathrm{Gar},t_1}
\left({\theta^\infty_1};{t_1 \atop t_2};
{q_1,p_1 \atop q_2,p_2}\right)
=
H_{\mathrm{II}}\big({\theta^\infty_1};t_1;q_1,p_1\big)
-2p_2q_2q_1-t_2p_2-q_2,\\
t_2H^{\frac52+2}_{\mathrm{Gar},t_2}
\left({\theta^\infty_1};{t_1 \atop t_2};
{q_1,p_1 \atop q_2,p_2}\right)
={p_2}^2{q_2}^2+\theta^\infty_1p_2q_2+t_2p_2\big(p_1-{q_1}^2-t_1\big)-p_1q_2-t_2q_1.
\end{gather*}

\subsubsection*{Singularity pattern $\boldsymbol{\frac32+\frac32+1}$}
The Riemann scheme is given by
\begin{gather*}
\left(
\begin{matrix}
x=0 \, \left( \frac12 \right) & x=t_1 & x=\infty \, \left( \frac12 \right)\\
\overbrace{
	\begin{matrix}
	\sqrt{t_2} & 0\\
	-\sqrt{t_2} & 0
	\end{matrix}}
	&\begin{matrix}
	0\\
	\theta^{t_1}
	\end{matrix}
&\overbrace{
	\begin{matrix}
	1 & \theta^\infty_1/2 \\
	-1 & \theta^\infty_1/2
	\end{matrix}}
\end{matrix}
\right),
\end{gather*}
and the Fuchs--Hukuhara relation is written as $\theta^{t_1}+\theta^\infty_1=0$.
The Lax pair is expressed as
\begin{gather*}
\frac{\partial Y}{\partial x}=\left( \frac{A_0^{(1)}}{x^2}+\frac{A_0^{(0)}}{x}+\frac{A_{t_1}}{x-t_1}+N \right)Y ,\\
\frac{\partial Y}{\partial t_1}=\left( N_1-\frac{A_{t_1}}{x-t_1} \right)Y ,\qquad
\frac{\partial Y}{\partial t_2}=\left( N_2-\frac{\frac{1}{t_2}A_0^{(1)}}{x} \right)Y .
\end{gather*}
Here
\begin{gather*}
A_0^{(1)}=
\begin{pmatrix}
0 & 0 \\
-q_2 & 0
\end{pmatrix},\qquad
A_0^{(0)}=
\begin{pmatrix}
-p_2q_2 & -t_2/q_2 \\
1-p_1 & p_2q_2
\end{pmatrix},\qquad
N=
 \begin{pmatrix}
 0 & 1 \\
 0 & 0
 \end{pmatrix},\\
A_{t_1}=
\begin{pmatrix}
q_1 \\
1
\end{pmatrix}
\begin{pmatrix}
p_1 & \theta^{t_1}-p_1q_1
\end{pmatrix},\qquad
N_1=
\frac{q_1(p_1q_1-\theta^{t_1})}{t_1}N,\qquad
N_2=\frac{1}{q_2}N.
\end{gather*}

The Hamiltonians are given by
\begin{gather*}
t_1H_{\mathrm{Gar},t_1}^{\frac{3}{2}+\frac{3}{2}+1}\left({\theta^\infty_1};{t_1 \atop t_2};{q_1,p_1 \atop q_2,p_2}\right)\\
\qquad{} =
t_1H_{\mathrm{III}(D_6)}\left({\theta^\infty_1, \theta^\infty_1}; t_1; q_1,p_1\right)
+\frac{q_1q_2}{t_1}\big(p_1q_1+\theta^\infty_1\big)-2p_1q_1p_2q_2
-\theta^\infty_1p_2q_2-\frac{t_2p_1}{q_2}, \nonumber \\
t_2H_{\mathrm{Gar},t_2}^{\frac{3}{2}+\frac{3}{2}+1}\left({\theta^\infty_1};{t_1 \atop t_2};{q_1,p_1 \atop q_2,p_2}\right)
=
t_2H_{\mathrm{III}(D_8)}\left(t_2; q_2, p_2\right)
-\frac{q_1q_2}{t_1}\big(p_1q_1+\theta^\infty_1\big)+\frac{t_2p_1}{q_2}.
\end{gather*}

\subsubsection*{Singularity pattern $\boldsymbol{\frac92}$}
The Riemann scheme is given by
\begin{gather*}
\left(
\begin{matrix}
x=\infty\,\left( \frac12 \right) \\
\overbrace{
	\begin{matrix}
	1 & 0 & 0 & 0 & \frac32 t_1 & 0 & -\frac{t_2}{2} & 0 \\
	-1 & 0 & 0 & 0 & -\frac32 t_1 & 0 & \frac{t_2}{2} & 0
	\end{matrix}}
\end{matrix}
\right).
\end{gather*}
The Lax pair is expressed as
\begin{gather*}
\frac{\partial Y}{\partial x} =\big(
A_0 x^3+A_1 x^2+A_2 x+A_3
 \big)Y ,\\
\frac{\partial Y}{\partial t_1} =\big( A_0 x^2+A_1 x+B_{10} \big)Y ,\qquad \frac{\partial Y}{\partial t_2}= ( -A_0 x+B_{20} )Y ,
\end{gather*}
where
\begin{gather*}
A_0 =
\begin{pmatrix}
0 & 1 \\
0 & 0
\end{pmatrix},\qquad
A_1=
\begin{pmatrix}
0 & p_1 \\
1 & 0
\end{pmatrix},\qquad
A_2=
\begin{pmatrix}
q_2 & {p_1}^2+p_2+2t_1 \\
-p_1 & -q_2 \\
\end{pmatrix},\\
A_3 =
\begin{pmatrix}
q_1-p_1q_2 & {p_1}^3+2p_1p_2-{q_2}^2+t_1p_1-t_2 \\
-p_2+t_1 & -q_1+p_1q_2
\end{pmatrix},\\
B_{10} =
\begin{pmatrix}
q_2 & {p_1}^2+2p_2+t_1 \\
-p_1 & -q_2
\end{pmatrix},\qquad
B_{20}=
\begin{pmatrix}
0 & -2p_1 \\
-1 & 0
\end{pmatrix}.
\end{gather*}

The Hamiltonians are given by
\begin{gather*}
H_{\mathrm{Gar},t_1}^{\frac92}
\left({t_1 \atop t_2};{q_1,p_1 \atop q_2,p_2}\right)
=-{p_1}^3p_2+{p_1}^2{q_2}^2+t_1{p_1}^3-2p_1q_1q_2-2p_1{p_2}^2+p_2{q_2}^2 \nonumber\\
 \qquad{} +t_1p_1p_2+{q_1}^2-t_1{q_2}^2+{t_1}^2p_1+t_2p_2, \nonumber \\
 H_{\mathrm{Gar},t_2}^{\frac92}
\left({t_1 \atop t_2};
{q_1,p_1 \atop q_2,p_2}\right)=
{p_1}^4+3{p_1}^2p_2+p_1{q_2}^2-2q_1q_2+{p_2}^2-t_2p_1+t_1p_2.
\end{gather*}

\subsubsection*{Singularity pattern $\boldsymbol{\frac52+\frac32}$}
The Riemann scheme is given by
\begin{gather*}
\left(
\begin{matrix}
x=0 \,\left( \frac12 \right) & x=\infty \, \left( \frac12 \right)	\\
\overbrace{
 \begin{matrix}
	\sqrt{t_2} & 0 \\
	-\sqrt{t_2} & 0
	\end{matrix}}
&\overbrace{
	\begin{matrix}
	1 & 0 & -t_1/2 & 0 \\
	-1 & 0 & t_1/2 & 0
	\end{matrix}}
\end{matrix}
\right).
\end{gather*}
The Lax pair is expressed as
\begin{gather*}
\frac{\partial Y}{\partial x}=\left( \frac{A_3}{x^2}+\frac{A_2}{x}+A_1+A_0 x \right)Y ,\\
\frac{\partial Y}{\partial t_1}=\left( x B_{11}+B_{10}\right)Y ,\qquad
\frac{\partial Y}{\partial t_2}=\left( -\frac{A_3}{t_2x}+B_{20}\right)Y ,
\end{gather*}
where
\begin{gather*}
A_0 =
\begin{pmatrix}
0 & 1 \\
0 & 0
\end{pmatrix},\qquad
A_1=
\begin{pmatrix}
q_1 & p_1-{q_1}^2-t_1\\
1 & -q_1 \\
\end{pmatrix},\\
A_2 =
\begin{pmatrix}
p_2q_2 & q_2 \\
-p_1 & -p_2q_2
\end{pmatrix},\qquad
A_3=
\begin{pmatrix}
0 & 0 \\
t_2/q_2 & 0
\end{pmatrix},\\
B_{11} =
\begin{pmatrix}
0 & -1 \\
0 & 0
\end{pmatrix},\qquad
B_{10}=
\begin{pmatrix}
-q_1 & 0 \\
-1 & q_1
\end{pmatrix},\qquad
B_{20}=
\begin{pmatrix}
0 & -q_2/t_2 \\
0 & 0
\end{pmatrix}.
\end{gather*}

The Hamiltonians are given by
\begin{gather*}
H^{\frac52+\frac32}_{\mathrm{Gar},t_1} \left({t_1 \atop t_2};{q_1,p_1 \atop q_2,p_2}\right)
= H_{\mathrm{II}}\left({0};t_1;q_1,p_1\right) -2p_2q_2q_1-q_2-\frac{t_2}{q_2},\\
t_2H^{\frac52+\frac32}_{\mathrm{Gar},t_2} \left({t_1 \atop t_2};{q_1,p_1 \atop q_2,p_2}\right)
={p_2}^2{q_2}^2-p_1q_2+\frac{t_2}{q_2}\big(p_1-{q_1}^2-t_1\big).
\end{gather*}

\subsection{Degenerate Fuji--Suzuki systems}\label{sec:FS_Lax}
The Fuji--Suzuki systems were discovered by Fuji and Suzuki~\cite{FS1} in their study on similarity reductions of the Drinfeld--Sokolov hierarchies. A family of Painlev\'e-type equations which includes the Fuji--Suzuki system with $A_5$-symmetry was independently proposed by Tsuda~\cite{Ts}.

Sakai~\cite{Sak2} derived the Fuji--Suzuki system of type $A_5$ from the isomonodromic deformation of the following Fuchsian system:
\begin{gather}\label{eq:21,21,111,111}
\frac{{\rm d}Y}{{\rm d}x}= \left( \frac{A_0}{x}+\frac{A_1}{x-1}+\frac{A_t}{x-t} \right)Y
\end{gather}
where $A_0$, $A_1$, and $A_t$ are $3 \times 3$ matrices satisfying the following conditions
\begin{gather*}
A_0 \sim \operatorname{diag}\big( 0, \theta^0_1, \theta^0_2 \big),\qquad
A_1 \sim \operatorname{diag}\big( 0, 0, \theta^1 \big),\qquad
A_t \sim \operatorname{diag}\big( 0, 0, \theta^t \big),
\end{gather*}
and
\begin{gather}\label{eq:Fuchs_rel}
A_\infty:=-(A_0+A_1+A_t)=\operatorname{diag}\big(\theta^\infty_1, \theta^\infty_2, \theta^\infty_3\big).
\end{gather}
Thus the spectral type of the Fuchsian system~(\ref{eq:21,21,111,111}) is $21,21,111,111$.
Taking the trace of~(\ref{eq:Fuchs_rel}), we have the Fuchs relation
\begin{gather*}
\theta^0_1+\theta^0_2+\theta^1+\theta^t+\theta^\infty_1+\theta^\infty_2+\theta^\infty_3=0.
\end{gather*}

The isomonodromic deformation equation of (\ref{eq:21,21,111,111}) is equivalent to
the Hamiltonian system
\begin{gather*}
\frac{{\rm d}q_i}{{\rm d}t}=\frac{\partial H_{\mathrm{FS}}^{A_5}}{\partial p_i}, \qquad
\frac{{\rm d}p_i}{{\rm d}t}=-\frac{\partial H_{\mathrm{FS}}^{A_5}}{\partial q_i}, \qquad i=1,2,
\end{gather*}
where the Hamiltonian is given by
\begin{gather*}
H_{\mathrm{FS}}^{A_5}\left({\theta^0_2+\theta^\infty_2,\theta^\infty_3,\theta^t
 \atop \theta^1,\theta^0_1,\theta^0_2};t;{q_1, p_1 \atop q_2, p_2}\right)
 =H_{\mathrm{VI}}\left({\theta^0_2+\theta^\infty_2,\theta^1+\theta^\infty_3
\atop \theta^t+\theta^\infty_3,\theta^0_1-\theta^0_2+1};t;q_1,p_1\right) \nonumber \\
\qquad {}+H_{\mathrm{VI}}\left({\theta^\infty_3,\theta^0_2+\theta^1+\theta^\infty_2
\atop \theta^0_2+\theta^t+\theta^\infty_2,\theta^0_1-\theta^0_2-\theta^\infty_2+1};t;q_2,p_2\right)\nonumber\\
\qquad{}+\frac{1}{t(t-1)}(q_1-t)(q_2-1)\big\{\big(p_1q_1-\theta^0_2-\theta^\infty_2\big)p_2+p_1\big(p_2q_2-\theta^\infty_3\big)\big\}.\nonumber
\end{gather*}

The linear system  of the spectral type $21,21,111,111$ has one deformation parameter.
However, linear systems which are degenerated from the $21,21,111,111$-system sometimes admit two-dimensional deformation,
since a degeneration process does not necessarily preserve the number of deformation parameters.
Such degenerations were pointed out in~\cite{KNS}.

In the present paper, the following linear systems
\begin{alignat*}{5}
& (1)_2(1), 21, 21, \qquad && ((1))(1)_2, 21,\qquad & & (((1)))_2 1, 21,\qquad & & ((1))_3, 21, 21,& \\
& (1)_2(1), (2)(1), \qquad && (((1)))(1)_2, \qquad && ((1))_3, (2)(1), \qquad && ((((1))))_3, 21,&
\end{alignat*}
which are degenerated systems of the $21,21,111,111$-system, also admit two-dimensional deformations and the Painlev\'e-type equations associated with these linear systems are degenerate Garnier systems.

The following are the Hamiltonians for the degenerate Fuji--Suzuki systems associated with ramif\/ied linear equations:
\begin{gather*}
tH_{\mathrm{Suz}}^{2+\frac{3}{2}}\left({\alpha, \beta \atop \gamma}; t; {q_1,p_1 \atop q_2,p_2}\right)\\
\qquad{} =
tH_{\mathrm{III}(D_7)}\left(\alpha; t; q_1,p_1\right)
+tH_{\mathrm{III}(D_7)}\left(\beta; t; q_2,p_2\right)+p_2q_1(p_1(q_1+q_2)+\gamma)-q_1,\\
tH_{\mathrm{KFS}}^{\frac32+\frac32}\left(\alpha, \beta; t; {q_1,p_1 \atop q_2,p_2}\right)\\
\qquad{}=tH_{\mathrm{III}(D_7)}\left( \alpha; t; q_1,p_1 \right)
+tH_{\mathrm{III}(D_7)}\left( \beta; t; q_2,p_2\right)
-p_1q_1p_2q_2-t(p_1p_2+p_1+p_2),\\
tH_{\mathrm{KFS}}^{\frac32+\frac43}\left(\alpha ;t; {q_1,p_1 \atop q_2,p_2}\right)\\
\qquad{} =
tH_{\mathrm{III}(D_7)}\left(\alpha; t; q_1,p_1\right)
+tH_{\mathrm{III}(D_7)}\left(1-\alpha; t; q_2,p_2\right)-p_1q_1p_2q_2-t\left(\frac{p_2}{q_1}+p_1+p_2\right),\\
tH_{\mathrm{KFS}}^{\frac43+\frac43}\left(t; {q_1,p_1 \atop q_2,p_2}\right)\\
\qquad{} =
tH_{\mathrm{III}(D_8)}\left(t; q_1,p_1\right)
+tH_{\mathrm{III}(D_8)}\left(t; q_2,p_2\right)-p_1q_1p_2q_2+\left(\frac{q_1q_2}{t}+q_1+q_2\right).
\end{gather*}
Here the Hamiltonian $H_{\mathrm{Suz}}^{2+\frac32}$ can be transformed into the Hamiltonian ${}_0\mathcal{H}_2$ in~\cite{Sz}
via the following canonical transformation:
\begin{gather*}
q_1 \to \frac{t}{q_1}, \qquad p_1 \to -\frac{q_1}{t}(p_1q_1+\gamma), \qquad
q_2 \to t p_2, \qquad p_2 \to -\frac{q_2}{t}.
\end{gather*}

\subsubsection*{Singularity pattern $\boldsymbol{2+1+1}$}

\underline{Spectral type $(1)_2(1), 21, 21$.}
The Riemann scheme is given by
\begin{gather*}
\left(
\begin{matrix}
x=0 & x=1 & x=\infty \, \left( \frac12 \right)\\
\begin{matrix} 0 \\ 0 \\ \theta^0 \end{matrix}&
\begin{matrix} 0 \\ 0 \\ \theta^1 \end{matrix}
&
\overbrace{\begin{matrix}
 0 & \sqrt{t_2} & \theta^\infty_1/2 \\
 0 & -\sqrt{t_2} & \theta^\infty_1/2 \\
-t_1 & 0 & \theta^\infty_2
 \end{matrix}}
\end{matrix}
\right),
\end{gather*}
and the Fuchs--Hukuhara relation is written as $\theta^0+\theta^1+\theta^\infty_1+\theta^\infty_2=0$. The Lax pair is expressed as
\begin{gather*}
\frac{\partial Y}{\partial x}= \left(A_\infty+\frac{A_0}{x}+\frac{A_1}{x-1} \right)Y ,\qquad
\frac{\partial Y}{\partial t_1}=(E_3x+B_1)Y,\qquad \frac{\partial Y}{\partial t_2}=\left( \frac{1}{t_2}Nx+B_2 \right)Y.
\end{gather*}
Here
\begin{gather*}
A_\infty=
\begin{pmatrix}
0 & 1 & 0 \\
0 & 0 & 0 \\
0 & 0 & t_1
\end{pmatrix} ,\qquad
A_0=
\begin{pmatrix}
0 \\
1 \\
q_1
\end{pmatrix}
\begin{pmatrix}
t_2(1-p_2) & a_0 & p_1(1-q_1)-\theta^\infty_2
\end{pmatrix},\\
A_1 =
\begin{pmatrix}
q_2 \\
t_2 \\
t_2(q_1-1)
\end{pmatrix}
\begin{pmatrix}
p_2 & a_1 & \frac{1}{t_2}(p_2q_2+\theta^\infty_2)
\end{pmatrix}, \\
E_3=\operatorname{diag}(0,0,1), \qquad
N=
\begin{pmatrix}
0 & 1 & 0 \\
0 & 0 & 0 \\
0 & 0 & 0
\end{pmatrix},\\
B_1=
\frac{1}{t_1}
\begin{pmatrix}
0 & -\frac{p_1q_1}{t_1}-\frac{q_1q_2}{t_2}(p_1q_1+\theta^\infty_2) & \frac{p_1}{t_1}+\frac{q_2}{t_2}(p_1q_1+\theta^\infty_2) \\
0 & 0 & p_1\\
t_2(q_1-p_2) & (B_1)_{32} & (B_1)_{33}
\end{pmatrix}
,\\
B_2=\frac{1}{t_2}
\begin{pmatrix}
p_2q_2+\theta^\infty_1 & \frac{p_1q_1}{t_1} & -\frac{p_1}{t_1}\\
t_2 & 1-p_2q_2 & 0\\
0 & \frac{t_2}{t_1}(p_2-q_1) & 1-q_2(q_1+p_2-1)
\end{pmatrix},
\end{gather*}
where
\begin{gather*}
a_0 =p_1q_1(q_1-1)+\theta^\infty_2q_1+\theta^0, \\
a_1 =\frac{1}{t_2}\big(\big(p_1q_1+\theta^\infty_2\big)(1-q_1)-p_2q_2+\theta^1\big), \\
(B_1)_{32}=(q_1-1)(p_1q_1-p_2q_2)+\big(\theta^\infty_2-\theta^\infty_1\big)q_1+\frac{t_2}{t_1}(q_1-p_2)-\theta^1-\theta^\infty_2,\\
(B_1)_{33}=(p_1+t_1)(q_1-1)+(p_1+q_2)q_1-p_2q_2.
\end{gather*}

The Hamiltonians are given by
\begin{gather*}
t_1H^{2+2+1}_{\mathrm{Gar},t_1}\left({-\theta^\infty_2,-\theta^\infty_1 \atop -\theta^0};{t_1 \atop t_2};
{q_1,p_1 \atop q_2,p_2}\right) =
t_1H_{\mathrm{V}}
\left({-\theta^0-\theta^\infty_2-\theta^\infty_1, -\theta^\infty_1+\theta^\infty_2
\atop \theta^0+\theta^\infty_1};t_1;q_1,p_1\right)\nonumber\\
\qquad{} +q_1q_2\big(p_1q_1+\theta^\infty_2\big)+p_2q_2\big({-}\theta^\infty_2+p_1-2p_1q_1\big)
-\frac{t_2}{t_1}p_1(p_2-q_1) ,\nonumber\\
 t_2H^{2+2+1}_{\mathrm{Gar},t_2}\left({-\theta^\infty_2,-\theta^\infty_1 \atop -\theta^0};{t_1 \atop t_2};
{q_1,p_1 \atop q_2,p_2}\right) =
t_2H_{\mathrm{III}(D_6)}
\big({\theta^0+\theta^\infty_1+\theta^\infty_2,\theta^\infty_1}
;t_2;q_2,p_2\big)\\
\qquad{} -\big(p_1q_1+\theta^\infty_2\big)q_2(q_1-1)+\frac{t_2}{t_1}p_1(p_2-q_1). \nonumber
\end{gather*}

\subsubsection*{Singularity pattern $\boldsymbol{3+1}$}

\underline{Spectral type $((1))(1)_2, 21$.}
The Riemann scheme is given by
\begin{gather*}
\left(
\begin{matrix}
 x=0 & x=\infty \,\left( \frac12 \right)\\
\begin{matrix} 0 \\ 0 \\ \theta^0 \end{matrix}
&\overbrace{\begin{matrix}
1 & 0 & t_2 & 0 & \theta^\infty_1 \\
0 & 0 & 0 & \sqrt{t_1} & \theta^\infty_2/2 \\
0 & 0 & 0 & -\sqrt{t_1} & \theta^\infty_2/2
 	\end{matrix}}
\end{matrix}
\right) ,
\end{gather*}
and the Fuchs--Hukuhara relation is written as $\theta^0+\theta^\infty_1+\theta^\infty_2=0$.
The Lax pair is expressed as
\begin{gather*}
\frac{\partial Y}{\partial x}=
\left( A_0 x+A_1+\frac{A_2}{x} \right)Y ,\qquad
\frac{\partial Y}{\partial t_1} =(B_{11}x+B_{10})Y ,\qquad
\frac{\partial Y}{\partial t_2}=(B_{21}x+B_{20})Y ,
\end{gather*}
where
\begin{gather*}
A_0=
\begin{pmatrix}
-1 & 0 & 0 \\
0 & 0 & 0 \\
0 & 0 & 0
\end{pmatrix},\qquad
A_1=
\begin{pmatrix}
-t_2 & 1 & 0 \\
-p_1q_1+p_2q_2-\theta^\infty_2 & 0 & 1 \\
q_1q_2+t_1 & 0 & 0
\end{pmatrix},\\
A_2=
\begin{pmatrix}
1 \\
p_2 \\
q_1
\end{pmatrix}
\begin{pmatrix}
-p_1q_1+p_2q_2+\theta^0 & -q_2 & p_1
\end{pmatrix},\\
B_{11} =
\begin{pmatrix}
0 & 0 & 0 \\
0 & 0 & 1/t_1 \\
0 & 0 & 0
\end{pmatrix}, \qquad
B_{10}=\frac{1}{t_1}
\begin{pmatrix}
0 & 0 & 1 \\
q_1q_2+t_1 & 0 & q_2+t_2 \\
0 & t_1 & 2p_1q_1+\theta^\infty_2+1
\end{pmatrix},\\
B_{21}=
\begin{pmatrix}
-1 & 0 & 0 \\
0 & 0 & 0 \\
0 & 0 & 0
\end{pmatrix},\qquad
B_{20}=
\begin{pmatrix}
0 & 1 & 0 \\
-p_1q_1+p_2q_2-\theta^\infty_2 & q_2+t_2 & 1-p_1 \\
q_1q_2+t_1 & 0 & q_2+t_2
\end{pmatrix}.
\end{gather*}

The Hamiltonians are given by
\begin{gather*}
H^{3+2}_{\mathrm{Gar},t_1}
\left({\theta^\infty_2, \theta^0};{t_1 \atop t_2};
{q_1,p_1 \atop q_2,p_2}\right)\\
\qquad =
H_{\mathrm{III}(D_6)}\left({-\theta^0,\theta^\infty_2+1};t_1;q_1,p_1\right)
-p_1-\frac{q_1q_2}{t_1}(q_2-p_2+t_2)+p_1p_2-q_2, \nonumber \\
H^{3+2}_{\mathrm{Gar},t_2}
\left({\theta^\infty_2, \theta^0};{t_1 \atop t_2};
{q_1,p_1 \atop q_2,p_2}\right) =
H_{\mathrm{IV}}\left({\theta^\infty_2,\theta^0};t_2;q_2,p_2\right)
-p_1q_1(p_2-2q_2-t_2)-q_1q_2+t_1p_1. \nonumber
\end{gather*}

\subsubsection*{Singularity pattern $\boldsymbol{\frac53+1+1}$}
\underline{Spectral type $((1))_3, 21, 21$.}
The Riemann scheme is given by
\begin{gather*}
\left(
\begin{matrix}
x=0 & x=1 & x=\infty \, \left( \frac13 \right)\\
\begin{matrix} 0 \\ 0 \\ \theta^0 \end{matrix}&
\begin{matrix} 0 \\ 0 \\ \theta^1 \end{matrix}
&
\overbrace{\begin{matrix}
-{t_1}^{\frac23} & \frac13{t_1}^{\frac13}t_2 & \theta^\infty_1/3 \\
-{t_1}^{\frac23}\omega & \frac13{t_1}^{\frac13}t_2\omega^2 & \theta^\infty_1/3 \\
-{t_1}^{\frac23}\omega^2 & \frac13{t_1}^{\frac13}t_2\omega & \theta^\infty_1/3
 \end{matrix}}
\end{matrix}
\right),
\end{gather*}
where $\omega=e^{2\pi i/3}$ is a cube root of unity.
The Fuchs--Hukuhara relation is written as $\theta^0+\theta^1+\theta^\infty_1$ $=0$.
The Lax pair is expressed as
\begin{gather*}
\frac{\partial Y}{\partial x} = \left(A_\infty+\frac{A_0}{x}+\frac{A_1}{x-1} \right)Y ,\qquad
\frac{\partial Y}{\partial t_1} =(B_{11}x+B_{10})Y,\qquad \frac{\partial Y}{\partial t_2}=(B_{21}x+B_{20})Y.
\end{gather*}
Here
\begin{gather*}
A_\infty =
\begin{pmatrix}
0 & 1 & 0 \\
0 & 0 & 1 \\
0 & 0 & 0
\end{pmatrix} ,\qquad
A_0=
\begin{pmatrix}
0 \\
0 \\
1
\end{pmatrix}
\begin{pmatrix}
-t_1q_1q_2 & a_0 & \theta^0
\end{pmatrix},\\
A_1 =
\begin{pmatrix}
\frac{q_2-p_2+t_2}{t_1q_1}+\frac{1}{{q_1}^2} \\
1/q_1 \\
1
\end{pmatrix}
\begin{pmatrix}
t_1(q_1q_2+t_1) & a_1 & p_1q_1+\frac{t_1p_2}{q_1}-\theta^0+\theta^1
\end{pmatrix}, \\
B_{11} =
\frac{1}{t_1}
\begin{pmatrix}
0 & 1 & 0 \\
0 & 0 & 1 \\
0 & 0 & 0
\end{pmatrix}, \qquad
B_{21}=-\frac{1}{t_1}
\begin{pmatrix}
0 & 0 & 1 \\
0 & 0 & 0 \\
0 & 0 & 0
\end{pmatrix},\\
B_{10}=
\frac{1}{t_1}
\begin{pmatrix}
(B_{10})_{11} & (B_{10})_{12} & 0 \\
q_2+\frac{t_1}{q_1} & (B_{10})_{22} & 0\\
t_1 & -q_2-t_2-\frac{t_1}{q_1} & \frac{1}{t_1}(p_1q_1-\theta^0+2)+\frac{p_2}{q_1}
\end{pmatrix},\\
B_{20} =
\begin{pmatrix}
-q_2-t_2-\frac{t_1}{q_1} & (B_{20})_{12} & 0\\
-t_1 & 0 & 0\\
0 & -t_1 & q_2+\frac{t_1}{q_1}
\end{pmatrix},
\end{gather*}
where
\begin{gather*}
a_0 =q_1q_2(q_2-p_2+t_2)+q_1\big(p_1q_1-\theta^0\big)+t_1q_2, \\
a_1 =q_1q_2(p_2-q_2-t_2)-q_1\big(p_1q_1-\theta^0\big)-t_1\left( 2q_2+\frac{t_1}{q_1}+t_2 \right), \\
(B_{10})_{11} =\frac{1}{t_1}\left( \left(q_2+\frac{t_1}{q_1}\right)(q_2-p_2+t_2)+\theta^\infty_1 \right)+\frac{q_1q_2+t_1}{{q_1}^2},\\
(B_{10})_{12} =\frac{1}{{t_1}^2}\left(\left(p_1q_1-q_2(p_2-q_2-t_2)
+\frac{t_1(3q_2+t_2)}{q_1}-\theta^0\right)(p_2-q_2-t_2)-2t_1p_1\right)\\
\hphantom{(B_{10})_{12} =}{} +\frac{2\theta^0-\theta^1}{t_1q_1}-\frac{3q_2+2t_2}{{q_1}^2}-\frac{t_1}{{q_1}^3},\\
(B_{10})_{22} =\frac{1}{t_1}\big(q_2(p_2-q_2-t_2)-p_1q_1-\theta^1+1\big)-\frac{2q_2+t_2}{q_1}-\frac{t_1}{{q_1}^2}, \\
(B_{20})_{12} =\frac{1}{t_1}\left( \left( q_2+\frac{2t_1}{q_1} \right)(p_2-q_2-t_2)+p_1q_1+\theta^1-1 \right)+\frac{t_2}{q_1}-\frac{t_1}{{q_1}^2}.
\end{gather*}

The Hamiltonians are given by
\begin{gather*}
 H^{3+2}_{\mathrm{Gar},t_1}
\left({\theta^1-\theta^0,\theta^0};{t_1 \atop t_2};
{q_1,p_1 \atop q_2,p_2}\right)\\
 \qquad{} =
H_{\mathrm{III}(D_6)}\left({-\theta^0,\theta^1-\theta^0+1};t_1;q_1,p_1\right)
-p_1-\frac{q_1q_2}{t_1}(q_2-p_2+t_2)+p_1p_2-q_2, \nonumber \\
H^{3+2}_{\mathrm{Gar},t_2}
\left({\theta^1-\theta^0,\theta^0};{t_1 \atop t_2};
{q_1,p_1 \atop q_2,p_2}\right)\\
 \qquad{} =
H_{\mathrm{IV}}\left({\theta^1-\theta^0,\theta^0};t_2;q_2,p_2\right)
-p_1q_1(p_2-2q_2-t_2)-q_1q_2+t_1p_1. \nonumber
\end{gather*}

\subsubsection*{Singularity pattern $\boldsymbol{\frac52+1}$}

\underline{Spectral type $(((1)))_2 1, 21$.} The Riemann scheme is given by
\begin{gather*}
\left(
\begin{matrix}
 x=0 & x=\infty \, \left( \frac12 \right) \\
\begin{matrix} 0 \\ 0 \\ \theta^0 \end{matrix}
&\overbrace{\begin{matrix}
 1 & t_2 & -t_1/2 & \theta^\infty_1/2 \\
 -1 & t_2 & t_1/2 & \theta^\infty_1/2 \\
 0 & 0 & 0 & \theta^\infty_2
 	\end{matrix}}
\end{matrix}
\right) ,
\end{gather*}
and the Fuchs--Hukuhara relation is written as
$\theta^0+\theta^\infty_1+\theta^\infty_2=0$.
The Lax pair is expressed as
\begin{gather*}
\frac{\partial Y}{\partial x}= \left( A_0 x+A_1+\frac{A_2}{x} \right)Y ,\qquad
\frac{\partial Y}{\partial t_1}=(B_{11}x+B_{10})Y ,\qquad
\frac{\partial Y}{\partial t_2}=(B_{21}x+B_{20})Y.
\end{gather*}
Here
\begin{gather*}
A_0 =
\begin{pmatrix}
0 & 1 & 0 \\
0 & 0 & 0 \\
0 & 0 & 0
\end{pmatrix},\qquad
A_1=
\begin{pmatrix}
q_1-t_2 & p_1-{q_1}^2-t_1 & -1 \\
1 & -q_1-t_2 & 0 \\
0 & p_2q_2-\theta^\infty_2 & 0
\end{pmatrix},\\
A_2 =
\begin{pmatrix}
0 \\
1 \\
-p_2
\end{pmatrix}
\begin{pmatrix}
-p_1 & p_2q_2+\theta^0 & q_2
\end{pmatrix},\\
B_{11} =
\begin{pmatrix}
0 & -1 & 0 \\
0 & 0 & 0 \\
0 & 0 & 0
\end{pmatrix}, \qquad
B_{10}=
\begin{pmatrix}
-q_1 & p_2 & 1 \\
-1 & q_1 & 0 \\
0 & -p_2q_2+\theta^\infty_2 & q_2-t_2
\end{pmatrix},\\
B_{21}=
\begin{pmatrix}
-1 & 0 & 0 \\
0 & -1 & 0 \\
0 & 0 & 0
\end{pmatrix},\qquad
B_{20}=
\begin{pmatrix}
0 & p_2(q_1-q_2+t_2) & q_1-q_2+t_2 \\
0 & 0 & 1 \\
-p_2q_2+\theta^\infty_2 & (B_{20})_{32} & (B_{20})_{33}
\end{pmatrix},
\end{gather*}
where
\begin{gather*}
(B_{20})_{32} =-p_1p_2+\big(p_2q_2-\theta^\infty_2\big)(q_1-t_2), \qquad
(B_{20})_{33} =p_1-t_1-{t_2}^2-q_2(q_1-t_2).
\end{gather*}

The Hamiltonians are given by
\begin{gather*}
 H^{4+1}_{\mathrm{Gar},t_1}\left({\theta^\infty_2, \theta^0};{t_1 \atop t_2};{q_1,p_1 \atop q_2,p_2}\right)
 =
H_{\rm{II}}\left({-\theta^0};t_1;q_1,p_1\right)
+p_2q_2(q_1-q_2+t_2)+p_1p_2+\theta^\infty_2 q_2, \nonumber \\
 H^{4+1}_{\mathrm{Gar},t_2}\left({\theta^\infty_2, \theta^0};{t_1 \atop t_2};{q_1,p_1 \atop q_2,p_2}\right)
 =
-{p_2}^2q_2-t_2p_2{q_2}^2+{t_2}^2p_2q_2+\theta^\infty_2 t_2q_2-\theta^0 p_2\\
\qquad{} +p_1p_2(q_1-2q_2+t_2)
 +q_1q_2\big(p_2q_2-\theta^\infty_2\big)+\theta^\infty_2 p_1+t_1p_2q_2.\nonumber
\end{gather*}

\subsubsection*{Singularity pattern $\boldsymbol{2+2}$}

\underline{Spectral type $(1)_2(1), (2)(1)$.} The Riemann scheme is given by
\begin{gather*}
\left(
\begin{matrix}
 x=0 & x=\infty \, \left( \frac12 \right)\\
\overbrace{\begin{matrix}
 0 & 0 \\
 0 & 0 \\
 1 & \theta^0
 \end{matrix}}
&
\overbrace{\begin{matrix}
 0 & \sqrt{-t_2} & \theta^\infty_1/2 \\
 0 & -\sqrt{-t_2} & \theta^\infty_1/2 \\
 -t_1 & 0 & \theta^\infty_2
 \end{matrix}}
\end{matrix}
\right) ,
\end{gather*}
and the Fuchs--Hukuhara relation is written as
$\theta^0+\theta^\infty_1+\theta^\infty_2=0$.
The Lax pair is expressed as
\begin{gather*}
\frac{\partial Y}{\partial x}=
\left(
\frac{A_2}{x^2}+\frac{A_1}{x}+A_0
\right)Y ,\qquad
\frac{\partial Y}{\partial t_1}=(E_3x+B_1)Y ,\qquad
\frac{\partial Y}{\partial t_2}=\left( \frac{1}{t_2}Nx+B_2 \right)Y.
\end{gather*}
Here
\begin{gather*}
A_0 =
\begin{pmatrix}
0 & 1 & 0 \\
0 & 0 & 0 \\
0 & 0 & t_1
\end{pmatrix},\qquad
A_2=
\begin{pmatrix}
0 \\
1 \\
1
\end{pmatrix}
\begin{pmatrix}
q_2 & 1-p_1 & p_1
\end{pmatrix},\\
A_1 =
\begin{pmatrix}
-p_2q_2 & (p_1-1)p_2 & -p_1p_2 \\
-t_2 & p_2q_2-\theta^\infty_1 & -p_1q_1-\theta^\infty_2 \\
q_1q_2-t_2 & (1-p_1)q_1+p_2q_2-\theta^\infty_1 & -\theta^\infty_2
\end{pmatrix}, \\
E_3 =\operatorname{diag}(0 \ 0 \ 1),\qquad
N=
\begin{pmatrix}
0 & 1 & 0 \\
0 & 0 & 0 \\
0 & 0 & 0
\end{pmatrix},\\
B_1 =
\frac{1}{t_1}
\begin{pmatrix}
0 & p_1\left( p_2+\frac{q_1}{t_1} \right)+\frac{\theta^\infty_2}{t_1} & -p_1(p_2+\frac{q_1}{t_1})-\frac{\theta^\infty_2}{t_1} \\
0 & 0 & -p_1q_1-\theta^\infty_2\\
q_1q_2-t_2 & (B_1)_{32} & (B_1)_{33}
\end{pmatrix},\\
B_2=\frac{1}{t_1t_2}
\begin{pmatrix}
-t_1(p_2q_2-\theta^\infty_1) & -p_1q_1-\theta^\infty_2 & p_1q_1+\theta^\infty_2\\
-t_1t_2 & t_1(p_2q_2+1) & 0\\
0 & -q_1q_2+t_2 & q_1q_2+t_1(p_2q_2+1)-t_2,
\end{pmatrix}
\end{gather*}
where
\begin{gather*}
(B_1)_{32} =q_1\left( -p_1+\frac{q_2}{t_1}+1 \right)+\left( p_2q_2-\theta^\infty_1-\frac{t_2}{t_1} \right), \\
(B_1)_{33} =-q_2\left( p_2+\frac{q_1}{t_1} \right)-q_1+\theta^\infty_1-\theta^\infty_2+\frac{t_2}{t_1}.
\end{gather*}

The Hamiltonians are given by
\begin{gather*}
t_1H_{\mathrm{Gar},t_1}^{2+\frac{3}{2}+1}
\left({-\theta^\infty_2, -\theta^\infty_1};{t_1 \atop t_2};
{q_1,p_1 \atop q_2,p_2}\right)
 =
t_1H_{\mathrm{III}(D_6)}\big({\theta^\infty_2, \theta^\infty_1+\theta^\infty_2}; t_1; q_1,p_1\big)
+t_2p_1p_2\\
 \qquad{}-p_2q_2\big(2p_1q_1+\theta^\infty_2\big)+\frac{t_2}{t_1}p_1q_1-\frac{q_1q_2}{t_1}\big(p_1q_1+\theta^\infty_2\big), \nonumber \\
t_2H_{\mathrm{Gar},t_2}^{2+\frac{3}{2}+1}
\left({-\theta^\infty_2, -\theta^\infty_1};{t_1 \atop t_2};
{q_1,p_1 \atop q_2,p_2}\right)
 =
t_2H_{\mathrm{III}(D_7)}\big({-}\theta^\infty_1+1; t_2; q_2,p_2\big)
-t_2p_1p_2\\
 \qquad{}-\frac{t_2}{t_1}p_1q_1+\frac{q_1q_2}{t_1}\big(p_1q_1+\theta^\infty_2\big).\nonumber
\end{gather*}

\subsubsection*{Singularity pattern $\boldsymbol{2+\frac53}$}

\underline{Spectral type $((1))_3, (2)(1)$.}
The Riemann scheme is given by
\begin{gather*}
\left(
\begin{matrix}
 x=0 & x=\infty \, \left( \frac13 \right)\\
\overbrace{\begin{matrix}
 0 & 0 \\
 0 & 0 \\
 t_ 1 & \theta^0
 \end{matrix}}
&
\overbrace{\begin{matrix}
-1 & \frac{t_2}{3} & \theta^\infty_1/3 \\
-\omega & \frac{\omega^2 t_2}{3} & \theta^\infty_1/3 \\
-\omega^2 & \frac{\omega t_2}{3} & \theta^\infty_1/3
 \end{matrix}}
\end{matrix}
\right) ,
\end{gather*}
and the Fuchs--Hukuhara relation is written as $\theta^0+\theta^\infty_1=0$. The Lax pair is expressed as
\begin{gather*}
\frac{\partial Y}{\partial x} =
\left(
\frac{A_2}{x^2}+\frac{A_1}{x}+A_0
\right)Y ,\qquad
\frac{\partial Y}{\partial t_1} =\left( \frac{B_{11}}{x}+B_{10} \right)Y ,\qquad
\frac{\partial Y}{\partial t_2}=\left( B_{21}x+B_{20} \right)Y ,
\end{gather*}
where
\begin{gather*}
A_0 =N,\qquad
A_2=t_1
\begin{pmatrix}
0 \\
0 \\
1
\end{pmatrix}
\begin{pmatrix}
-p_2 & p_1 & 1
\end{pmatrix},\qquad
N=
\begin{pmatrix}
0 & 1 & 0 \\
0 & 0 & 1 \\
0 & 0 & 0
\end{pmatrix},\\
A_1 =
\begin{pmatrix}
p_2q_2 & -p_1q_2 & -q_2 \\
-p_2q_1 & p_1q_1 & q_1 \\
1 & p_2q_1-t_2 & -p_1q_1-p_2q_2 -\theta^\infty_1
\end{pmatrix}, \\
B_{11} =
-\frac{1}{t_1}A_2, \qquad
B_{21}=-N^2,\\
B_{10} =
\frac{1}{t_1}
\begin{pmatrix}
0 & -q_1 & q_2 \\
0 & 0 & -q_1\\
0 & 0 & 0
\end{pmatrix},\qquad
B_{20}=
\begin{pmatrix}
p_2q_1-t_2 & -p_1q_1-2p_2q_2-\theta^\infty_1+1 & 0\\
-1 & 0 & 0\\
0 & -1 & -p_2q_1
\end{pmatrix}.
\end{gather*}

The Hamiltonians are given by
\begin{gather*}
t_1H^{3+\frac32}_{\mathrm{Gar},t_1}
\left({\theta^\infty_1};{t_1 \atop t_2};
{q_1,p_1 \atop q_2,p_2}\right)
 =
t_1H_{\mathrm{III}(D_7)}\big({\theta^\infty_1};t_1;q_1,p_1\big)
+p_2q_2\big(p_2q_2+\theta^\infty_1\big)-q_2\\
 \qquad{}+q_1(p_2q_1+2p_1p_2q_2-t_2-1), \nonumber \\
H^{3+\frac32}_{\mathrm{Gar},t_2}
\left({\theta^\infty_1};{t_1 \atop t_2};
{q_1,p_1 \atop q_2,p_2}\right)
 =
p_2q_1\big(p_1q_1+2p_2q_2+\theta^\infty_1-1\big)+p_1q_2-q_1+t_1p_2-t_2p_2q_2.
\end{gather*}

\subsubsection*{Singularity pattern 4}

\underline{Spectral type $(((1)))(1)_2$.}
The Riemann scheme is given by
\begin{gather*}
\left(
\begin{matrix}
 x=\infty \, \left( \frac12 \right)\\
\overbrace{\begin{matrix}
1 & 0 & 0 & 0 & t_1 & 0 & \theta^\infty_1\\
0 & 0 & 0 & 0 & 0 & \sqrt{t_2} & \theta^\infty_2/2\\
0 & 0 & 0 & 0 & 0 & -\sqrt{t_2} & \theta^\infty_2/2
 	\end{matrix}}
\end{matrix}
\right) ,
\end{gather*}
and the Fuchs--Hukuhara relation is written as $\theta^\infty_1+\theta^\infty_2=0$.
The Lax pair is expressed as
\begin{gather*}
\frac{\partial Y}{\partial x} =
\big( A_0x^2+A_1x+A_2 \big)Y ,\qquad
\frac{\partial Y}{\partial t_1} =(B_{11}x+B_{10})Y ,\qquad
\frac{\partial Y}{\partial t_2}=(B_{21}x+B_{20})Y ,
\end{gather*}
where
\begin{gather*}
A_0 =
\begin{pmatrix}
-1 & 0 & 0 \\
 0 & 0 & 0 \\
 0 & 0 & 0
\end{pmatrix},\qquad
A_1=
\begin{pmatrix}
0 & 1 & 0 \\
p_1 & 0 & 0 \\
p_2 & 0 & 0
\end{pmatrix},\qquad
A_2=
\begin{pmatrix}
p_1-t_1 & -q_1 & -q_2 \\
p_1q_1+p_2q_2+\theta^\infty_1 & -p_1 & t_2 \\
p_2q_1+1 & -p_2 & 0
\end{pmatrix}, \\
B_{11} =
\begin{pmatrix}
-1 & 0 & 0 \\
0 & 0 & 0 \\
0 & 0 & 0
\end{pmatrix},\qquad
B_{10}=
\begin{pmatrix}
0 & 1 & 0 \\
p_1 & q_1 & q_2 \\
p_2 & 0 & q_1
\end{pmatrix},\\
B_{21} =
\begin{pmatrix}
0 & 0 & 0 \\
0 & 0 & 1 \\
0 & 0 & 0
\end{pmatrix},\qquad
B_{20}=
\begin{pmatrix}
0 & 0 & 1 \\
p_2 & 0 & q_1 \\
0 & 1/t_2 & -\frac{1}{t_2}(2p_2q_2+\theta^\infty_1)
\end{pmatrix}.
\end{gather*}

The Hamiltonians are given by
\begin{gather*}
H^{\frac52+2}_{\mathrm{Gar},t_1}
\left({\theta^\infty_1};{t_1 \atop t_2};
{q_1,p_1 \atop q_2,p_2}\right)
 =
H_{\mathrm{II}}\big({\theta^\infty_1};t_1;q_1,p_1\big)
-2p_2q_2q_1-t_2p_2-q_2,\\
t_2H^{\frac52+2}_{\mathrm{Gar},t_2}
\left({\theta^\infty_1};{t_1 \atop t_2};
{q_1,p_1 \atop q_2,p_2}\right)
 ={p_2}^2{q_2}^2+\theta^\infty_1 p_2q_2+t_2p_2\big(p_1-{q_1}^2-t_1\big)-p_1q_2-t_2q_1.
\end{gather*}

\subsubsection*{Singularity pattern $\boldsymbol{\frac73+1}$}

\underline{Spectral type $((((1))))_3, 21$.}
The Riemann scheme is given by
\begin{gather*}
\left(
\begin{matrix}
 x=0 & x=\infty \,\left( \frac13 \right)\\
\begin{matrix} 0 \\ 0 \\ \theta^0 \end{matrix}
&\overbrace{\begin{matrix}
1 & 0 & -\frac23 t_2 & -t_1/3 & \theta^\infty_1/3 \\
\omega & 0 & -\frac23 \omega^2 t_2 & -\omega t_1/3 & \theta^\infty_1/3 \\
\omega^2 & 0 & -\frac23 \omega t_2 & -\omega^2 t_1/3 & \theta^\infty_1/3
 	\end{matrix}}
\end{matrix}
\right) ,
\end{gather*}
and the Fuchs--Hukuhara relation is written as $\theta^0+\theta^\infty_1=0$.
The Lax pair is expressed as
\begin{gather*}
\frac{\partial Y}{\partial x} =
\left( A_0 x+A_1+\frac{A_2}{x} \right)Y ,\qquad
\frac{\partial Y}{\partial t_1} =(B_{11}x+B_{10})Y ,\qquad
\frac{\partial Y}{\partial t_2}=(B_{21}x+B_{20})Y ,
\end{gather*}
where
\begin{gather*}
A_0 =
\begin{pmatrix}
0 & 1 & 0 \\
0 & 0 & 0 \\
0 & 0 & 0
\end{pmatrix},\qquad
A_1=
\begin{pmatrix}
0 & p_1-q_1q_2-t_1 & -p_2 \\
0 & 0 & -1 \\
1 & q_2-p_2+2t_2 & 0
\end{pmatrix},\\
A_2 =
\begin{pmatrix}
0 \\
1 \\
q_1
\end{pmatrix}
\begin{pmatrix}
q_2 & -p_1q_1+\theta^0 & p_1
\end{pmatrix},\\
B_{11} =
\begin{pmatrix}
0 & -1 & 0 \\
0 & 0 & 0 \\
0 & 0 & 0
\end{pmatrix}, \qquad
B_{10}=
\begin{pmatrix}
0 & -q_1(p_2+q_2) & p_2+q_2 \\
0 & 0 & 1 \\
-1 & p_2-2(q_2+t_2) & 0
\end{pmatrix},\\
B_{21} =
\begin{pmatrix}
0 & 0 & -1 \\
0 & 0 & 0 \\
0 & 1 & 0
\end{pmatrix},\qquad
B_{20}=
\begin{pmatrix}
q_2-p_2+t_2 & q_1(p_1-2q_1q_2-t_1) & -p_1+2q_1q_2+t_1 \\
-1 & p_2-t_2 & 0 \\
0 & 2p_1-q_1q_2-t_1 & -q_2-t_2
\end{pmatrix}.
\end{gather*}

The Hamiltonians are given by
\begin{gather*}
H^{5}_{\mathrm{Gar},t_1}\left({\theta^0};{t_1 \atop t_2};{q_1,p_1 \atop q_2,p_2}\right)
=
-q_1\big(p_1q_1-\theta^0\big)+q_2(q_1(p_2+q_2)-2p_1+t_1)+p_1(p_2-2t_2), \nonumber \\
H^{5}_{\mathrm{Gar},t_2}\left({\theta^0};{t_1 \atop t_2};{q_1,p_1 \atop q_2,p_2}\right)
=
H_{\mathrm{IV}}\big({-}1,\theta^0 ; 2t_2; q_2, p_2\big)\\
\qquad{}
+q_1q_2(q_1q_2-2p_1+t_1)+p_1(p_1-p_2q_1-t_1). \nonumber
\end{gather*}

\subsubsection*{Singularity pattern $\boldsymbol{\frac32+1+1}$}
\underline{Spectral type $(1)_2 1, 21, 111$.}
The Riemann scheme is given by
\begin{gather*}
\left(
\begin{matrix}
 x=0 & x=1 & x=\infty \, \left( \frac12 \right) \\
\begin{matrix} 0 \\ \theta^0_1 \\ \theta^0_2 \end{matrix}
& \begin{matrix} 0 \\ 0 \\ \theta^1 \end{matrix}
&\overbrace{\begin{matrix}
 \sqrt{-t} & \theta^\infty_1/2 \\
 -\sqrt{-t} & \theta^\infty_1/2 \\
 0 & \theta^\infty_2 \\
 	\end{matrix}}
\end{matrix}
\right) ,
\end{gather*}
and the Fuchs--Hukuhara relation is written as
$\theta^0_1+\theta^0_2+\theta^1+\theta^\infty_1+\theta^\infty_2=0$.

The Lax pair is expressed as
\begin{gather*}
\frac{\partial Y}{\partial x} =\left(\frac{A_0}{x}+\frac{A_1}{x-1}+A_\infty\right)Y ,\qquad
\frac{\partial Y}{\partial t} =\left(\frac{1}{t}A_\infty x+B_0\right)Y .
\end{gather*}
Here
\begin{gather*}
A_{\xi} =
P^{-1}\hat{A}_{\xi}P \ \ (\xi=0,1),\qquad
P=
\begin{pmatrix}
0 & 1 & 0 \\
1 & 0 & -t \\
\frac{1}{t}(p_1q_1-\theta^\infty_2) & 0 & 0
\end{pmatrix},\\
\hat{A}_0 =
\begin{pmatrix}
\theta^0_1 & -p_2 & q_1 \\
0 & \theta^0_2 & q_1q_2-t \\
0 & 0 & 0
\end{pmatrix}, \qquad
\hat{A}_1=
\begin{pmatrix}
1 \\
q_2 \\
p_1
\end{pmatrix}
\begin{pmatrix}
p_1q_1+(1-p_2)q_2+\theta^1 & p_2-1 & -q_1
\end{pmatrix},\\
A_\infty =
\begin{pmatrix}
0 & t & 0 \\
0 & 0 & 0 \\
0 & 0 & 0
\end{pmatrix}, \qquad
B_0=\frac{1}{t}
\begin{pmatrix}
(A_0+A_1)_{11} & 0 & 0 \\
-1 & (A_0+A_1)_{22} & t \\
(A_0+A_1)_{31} & 0 & (p_2-1)q_2+\theta^0_2-1
\end{pmatrix}.
\end{gather*}
The Hamiltonian is given by
\begin{gather*}
 tH_{\mathrm{FS}}^{A_3}
\left({-\theta^\infty_2,\, 1+\theta^0_1+\theta^1 \atop \theta^0_2-\theta^0_1,\, \theta^0_2-\theta^0_1-\theta^1};
t; {q_1,p_1\atop q_2,p_2}\right)  =tH_{\mathrm{III}(D_6)}\big( {-}\theta^\infty_2,\, 1+\theta^0_1+\theta^1 ; t; q_1,p_1 \big)\\
\qquad{} +tH_{\mathrm{III}(D_6)}\big( {\theta^0_2-\theta^0_1,\, \theta^0_2-\theta^0_1-\theta^1}; t; q_2,p_2\big)
 -p_1p_2(q_1q_2+t).\nonumber
\end{gather*}
\begin{Remark}\label{rem:FSA3}
Note that the expression of the Hamiltonian for the degenerate Fuji--Suzuki system of type $A_3$ is slightly dif\/ferent from that in~\cite{KNS}.
Here we adopted the following expression:
\begin{gather*}
 tH_{\mathrm{FS}}^{A_3}\left( {\alpha_,\beta \atop \gamma,\delta}; t; {q_1, p_1 \atop q_2, p_2} \right)
 =
tH_{\mathrm{III}(D_6)} (\alpha, \beta; t; q_1,p_1 )
+tH_{\mathrm{III}(D_6)} (\gamma, \delta; t; q_2,p_2 )-p_1p_2(q_1q_2+t).
\end{gather*}
This Hamiltonian can be transformed into the expression of $H_{\mathrm{FS}}^{A_3}$ in \cite{KNS} via the following
canonical transformation:
\begin{gather*}
q_2 \to \frac{t}{q_2}, \qquad p_2 \to -\frac{q_2}{t}(p_2q_2+\gamma).
\end{gather*}
\end{Remark}

\subsubsection*{Singularity pattern $\boldsymbol{2+\frac32}$}

\underline{Spectral type $(1)_2 1, (11)(1)$.}
The Riemann scheme is given as
\begin{gather*}
\left(
\begin{matrix}
 x=0 & x=\infty \, \left( \frac12 \right) \\
\overbrace{\begin{matrix}
 1 & \theta^0_1 \\
 0 & \theta^0_2 \\
 0 & 0
 \end{matrix}}
&
\overbrace{\begin{matrix}
 \sqrt{t} & \theta^\infty_1/2 \\
 -\sqrt{t} & \theta^\infty_1/2 \\
 0 & \theta^\infty_2
 \end{matrix}}
\end{matrix}
\right),
\end{gather*}
and the Fuchs--Hukuhara relation is written as
$\theta^0_1+\theta^0_2+\theta^\infty_1+\theta^\infty_2=0$.
The Lax pair is expressed as
\begin{gather}
\frac{\partial Y}{\partial x}= \left(\frac{A_2}{x^2}+\frac{A_1}{x}+A_0 \right)Y ,\qquad
\frac{\partial Y}{\partial t}=\left( \frac{1}{t}A_0 x+B_0 \right)Y,\label{eq:(1)_2,(11)(1)}
\end{gather}
where
\begin{gather*}
A_0=
\begin{pmatrix}
0 & t & 0 \\
0 & 0 & 0 \\
0 & 0 & 0
\end{pmatrix},\\
A_1=
\begin{pmatrix}
(p_1+p_2)q_1 & q_1 & p_2q_1 \\
1 & -p_1q_1-p_2q_2+\theta^0_1 & 1 \\
p_2(q_2-q_1)+\theta^0_2+\theta^\infty_2 & q_2-q_1 & p_2(q_2-q_1)+\theta^0_2
\end{pmatrix},\\
A_2 =
\begin{pmatrix}
0 & 0 & 0 \\
p_1+p_2 & 1 & p_2 \\
0 & 0 & 0
\end{pmatrix},\\
B_0=\frac{1}{t}
\begin{pmatrix}
(p_1+p_2)q_1 & 0 & 0 \\
1 & -p_1q_1-p_2q_2+\theta^0_1 & 1 \\
p_2(q_2-q_1)+\theta^0_2+\theta^\infty_2 & 0 & p_2q_2+\theta^0_2
\end{pmatrix}.
\end{gather*}

The Hamiltonian is given by
\begin{gather*}
tH_{\mathrm{Suz}}^{2+\frac{3}{2}}
\left({-\theta^0_1,\, \theta^0_2-\theta^0_1 \atop \theta^\infty_2}; t; {q_1,p_1 \atop q_2,p_2}\right)=tH_{\mathrm{III}(D_7)}\big({-}\theta^0_1; t; q_1,p_1\big)\\
\qquad{}
+tH_{\mathrm{III}(D_7)}\big(\theta^0_2-\theta^0_1; t; q_2,p_2\big)+p_2q_1\big(p_1(q_1+q_2)+\theta^\infty_2\big)-q_1.\nonumber
\end{gather*}

\subsubsection*{Singularity pattern $\boldsymbol{\frac43+1+1}$}

\underline{Spectral type $(1)_3, 21, 111$.}
The Riemann scheme is given by
\begin{gather*}
\left(
\begin{matrix}
 x=0 & x=1 & x=\infty \, \left( \frac13 \right) \\
\begin{matrix} 0 \\ \theta^0_1 \\ \theta^0_2 \end{matrix}
& \begin{matrix} 0 \\ 0 \\ \theta^1 \end{matrix}
&\overbrace{\begin{matrix}
 t^{\frac13} & \theta^\infty_1/3 \\
 \omega t^{\frac13} & \theta^\infty_1/3 \\
 \omega^2 t^{\frac13} & \theta^\infty_1/3 \\
 	\end{matrix}}
\end{matrix}
\right) ,
\end{gather*}
and the Fuchs--Hukuhara relation is written as
$\theta^0_1+\theta^0_2+\theta^1+\theta^\infty_1=0$.
The Lax pair is expressed as
\begin{gather*}
\frac{\partial Y}{\partial x} =
\left(
\frac{A_0}{x}+\frac{A_1}{x-1}+A_\infty
\right)Y ,\qquad
\frac{\partial Y}{\partial t} =\left(\frac{1}{t}A_\infty x+B_0\right)Y .
\end{gather*}
Here
\begin{gather*}
A_0 =
\begin{pmatrix}
0 & 0 & 0 \\
-q_2 & \theta^0_1 & p_1q_1-\theta^0_1+\theta^0_2 \\
q_1-q_2 & 0 & \theta^0_2
\end{pmatrix}, \qquad
A_\infty=
\begin{pmatrix}
0 & 1 & 0 \\
0 & 0 & 0 \\
0 & 0 & 0
\end{pmatrix}, \\
A_1 =
\begin{pmatrix}
-p_2 \\
1 \\
1
\end{pmatrix}
\begin{pmatrix}
q_2 & p_1q_1+p_2q_2+\frac{t}{q_1}-\theta^0_1+\theta^0_2+\theta^1 & -p_1q_1-\frac{t}{q_1}+\theta^0_1-\theta^0_2
\end{pmatrix},\\
B_0 =\frac{1}{t}
\begin{pmatrix}
-p_2q_2-1 & 0 & 0 \\
0 & p_1q_1+p_2q_2+\frac{t}{q_1}+\theta^0_2+\theta^1 & -t/q_1 \\
q_1 & 0 & p_1q_1+p_2(q_1+q_2)+\theta^0_2+\theta^1
\end{pmatrix}.
\end{gather*}

The Hamiltonian is given by
\begin{gather*}
 tH_{\mathrm{Suz}}^{2+\frac{3}{2}}\left({\theta^0_2-\theta^0_1+\theta^1+1,\, \theta^0_2+\theta^1+1 \atop \theta^0_2-\theta^0_1};
t; {q_1,p_1\atop q_2,p_2}\right) =tH_{\mathrm{III}(D_7)}\big( \theta^0_2-\theta^0_1+\theta^1+1; t; q_1,p_1 \big)\\
\qquad{} +tH_{\mathrm{III}(D_7)}\big( \theta^0_2+\theta^1+1; t; q_2,p_2\big) +p_2q_1(p_1(q_1+q_2)+\theta^0_2-\theta^0_1)-q_1.
\end{gather*}

\subsubsection*{Singularity pattern $\boldsymbol{\frac32+\frac32}$}

\underline{Spectral type $(1)_2 1, (1)_2 1$.}
The Riemann scheme is given by
\begin{gather*}
\left(
\begin{matrix}
 x=0 \, \left( \frac12 \right)& x=\infty \, \left( \frac12 \right) \\
\overbrace{\begin{matrix}
 0 & 0 \\
 1 & \theta^0/2 \\
 -1 & \theta^0/2
 \end{matrix}}
&
\overbrace{\begin{matrix}
 \sqrt{t} & \theta^\infty_1/2 \\
 -\sqrt{t} & \theta^\infty_1/2 \\
 0 & \theta^\infty_2
 \end{matrix}}
\end{matrix}
\right),
\end{gather*}
and the Fuchs--Hukuhara relation is written as
$\theta^0+\theta^\infty_1+\theta^\infty_2=0$.
The Lax pair is expressed as
\begin{gather*}
\frac{\partial Y}{\partial x} =
\left(
\frac{A_2}{x^2}+\frac{A_1}{x}+A_0
\right)Y ,\qquad
\frac{\partial Y}{\partial t}  =\left( \frac{1}{t}A_0 x+B_0 \right)Y,
\end{gather*}
where
\begin{gather*}
A_0 =
\begin{pmatrix}
0 & t & 0 \\
0 & 0 & 0 \\
0 & 0 & 0
\end{pmatrix},\qquad
A_1=
\begin{pmatrix}
-p_2q_2 & 0 & -tp_2 \\
1 & -p_1q_1+p_2q_2-\theta^\infty_1 & q_1 \\
p_1 & 1 & p_1q_1-\theta^\infty_2
\end{pmatrix},\\
A_2 =
\begin{pmatrix}
0 & 0 & 0 \\
q_2/t & 0 & 1 \\
0 & 0 & 0
\end{pmatrix}, \qquad
B_0=\frac{1}{t}
\begin{pmatrix}
-p_2q_2+\theta^\infty_1 & 0 & 0 \\
1 & -p_1q_1+p_2q_2 & q_1 \\
p_1 & 0 & -p_1q_1+p_2q_2
\end{pmatrix}.
\end{gather*}

The Hamiltonian is given by
\begin{gather*}
 tH_{\mathrm{KFS}}^{\frac32+\frac32}\left(\theta^\infty_1-\theta^\infty_2,\, 1-\theta^\infty_1 ;t; {q_1,p_1 \atop q_2,p_2}\right)
 =tH_{\mathrm{III}(D_7)}\big( \theta^\infty_1-\theta^\infty_2; t; q_1,p_1 \big)\\
\qquad{} +tH_{\mathrm{III}(D_7)}\big( 1-\theta^\infty_1; t; q_2,p_2\big)
-p_1q_1p_2q_2-t(p_1p_2+p_1+p_2).\nonumber
\end{gather*}

\subsubsection*{Singularity pattern $\boldsymbol{2+\frac43}$}
\underline{Spectral type $(1)_3, (11)(1)$.}
The Riemann scheme is given by
\begin{gather*}
\left(
\begin{matrix}
 x=0 & x=\infty \, \left( \frac13 \right) \\
\overbrace{\begin{matrix}
 1 & \theta^0_2 \\
 0 & \theta^0_1 \\
 0 & 0 \end{matrix}}
&
\overbrace{\begin{matrix}
 -t^{\frac13} & \theta^\infty_1/3 \\
 -\omega t^{\frac13} & \theta^\infty_1/3 \\
 -\omega^2 t^{\frac13} & \theta^\infty_1/3
 \end{matrix}}
\end{matrix}
\right),
\end{gather*}
and the Fuchs--Hukuhara relation is written as
$\theta^0_1+\theta^0_2+\theta^\infty_1=0$.
The Lax pair is expressed as
\begin{gather}
\frac{\partial Y}{\partial x} =
\left(
\frac{A_2}{x^2}+\frac{A_1}{x}+A_0
\right)Y ,\qquad
\frac{\partial Y}{\partial t} =\left( \frac{1}{t}A_0 x+B_0 \right)Y ,\label{eq:(1)_3,(11)(1)}
\end{gather}
where
\begin{gather*}
A_0=
\begin{pmatrix}
0 & 1 & 0 \\
0 & 0 & 0 \\
0 & 0 & 0
\end{pmatrix},\qquad
A_1=
\begin{pmatrix}
-p_2q_2 & -p_2 & -p_1p_2 \\
0 & -p_1q_1+p_2q_2+\theta^0_2 & 1 \\
t & q_1 & p_1q_1+\theta^0_1
\end{pmatrix},\\
A_2=
\begin{pmatrix}
0 & 0 & 0 \\
q_2 & 1 & p_1 \\
0 & 0 & 0
\end{pmatrix},\qquad
B_0=\frac1t
\begin{pmatrix}
-p_2q_2-1 & 0 & 0 \\
0 & -p_1q_1+p_2q_2+\theta^0_2 & 1 \\
t & 0 & p_1q_1+\theta^0_1
\end{pmatrix}.
\end{gather*}

The Hamiltonian is given by
\begin{gather*}
 tH_{\mathrm{KFS}}^{\frac32+\frac32}\left({\theta^0_1-\theta^0_2,\, \theta^0_2+1}; t; {q_1,p_1 \atop q_2,p_2}\right)
 =tH_{\mathrm{III}(D_7)}\big(\theta^0_1-\theta^0_2; t; q_1,p_1\big)\\
 \qquad{}
+tH_{\mathrm{III}(D_7)}\big(\theta^0_2+1; t; q_2,p_2\big)-p_1q_1p_2q_2-t(p_1p_2+p_1+p_2).\nonumber
\end{gather*}

\subsubsection*{Singularity pattern $\boldsymbol{\frac32+\frac43}$}

\underline{Spectral type $(1)_3, (1)_2 1$.} The Riemann scheme is given by
\begin{gather*}
\left(
\begin{matrix}
 x=0 \, \left( \frac12 \right)& x=\infty \, \left( \frac13 \right) \\
\overbrace{\begin{matrix}
 0 & 0 \\
 1 & \theta^0/2 \\
 -1 & \theta^0/2
 \end{matrix}}
&
\overbrace{\begin{matrix}
 -t^{\frac13} & \theta^\infty_1/3 \\
-\omega t^{\frac13} & \theta^\infty_1/3 \\
-\omega^2 t^{\frac13} & \theta^\infty_1/3
 \end{matrix}}
\end{matrix}
\right),
\end{gather*}
and the Fuchs--Hukuhara relation is written as $\theta^0+\theta^\infty_1=0$. The Lax pair is expressed as
\begin{gather*}
\frac{\partial Y}{\partial x}=
\left(\frac{A_2}{x^2}+\frac{A_1}{x}+A_0\right)Y ,\qquad
\frac{\partial Y}{\partial t}=\left( \frac{1}{t}A_0x+B_0 \right)Y,
\end{gather*}
where
\begin{gather*}
A_0=
\begin{pmatrix}
0 & 1 & 0 \\
0 & 0 & 0 \\
0 & 0 & 0
\end{pmatrix},\qquad
A_1=
\begin{pmatrix}
-p_2q_2 & 0 & -p_2 \\
0 & -p_1q_1+p_2q_2-\theta^\infty_1 & q_1 \\
t/q_1 & 1 & p_1q_1
\end{pmatrix},\\
A_2=
\begin{pmatrix}
0 & 0 & 0 \\
q_2 & 0 & 1 \\
0 & 0 & 0
\end{pmatrix}, \qquad
B_0 =\frac{1}{t}
\begin{pmatrix}
-p_2q_2-1 & 0 & 0 \\
0 & -p_1q_1+p_2q_2-\theta^\infty_1 & q_1 \\
t/q_1 & 0 & -p_1q_1+p_2q_2-\theta^\infty_1
\end{pmatrix}.
\end{gather*}

The Hamiltonian is given by
\begin{gather*}
 tH_{\mathrm{KFS}}^{\frac32+\frac43}\left(\theta^\infty_1 ;t; {q_1,p_1 \atop q_2,p_2}\right) =
tH_{\mathrm{III}(D_7)}\big(\theta^\infty_1; t; q_1,p_1\big)\\
\qquad{} +tH_{\mathrm{III}(D_7)}\big(1-\theta^\infty_1; t; q_2,p_2\big)-p_1q_1p_2q_2-t\left(\frac{p_2}{q_1}+p_1+p_2\right)\nonumber.
\end{gather*}

\subsubsection*{Singularity pattern $\boldsymbol{\frac43+\frac43}$}
\underline{Spectral type $(1)_3, (1)_3$.}
The Riemann scheme is given by
\begin{gather*}
\left(
\begin{matrix}
 x=0 \, \left( \frac13 \right)& x=\infty \, \left( \frac13 \right) \\
\overbrace{\begin{matrix}
 1 & 0 \\
 \omega & 0 \\
 \omega^2 & 0
 \end{matrix}}
&
\overbrace{\begin{matrix}
 -t^{\frac13} & 0 \\
-\omega t^{\frac13} & 0 \\
-\omega^2t^{\frac13} & 0
 \end{matrix}}
\end{matrix}
\right).
\end{gather*}

The Lax pair is expressed as
\begin{gather*}
\frac{\partial Y}{\partial x} =
\left(
\frac{A_2}{x^2}+\frac{A_1}{x}+A_0
\right)Y ,\qquad
\frac{\partial Y}{\partial t} =\left( \frac{1}{t}A_0x+B_0 \right)Y,
\end{gather*}
where
\begin{gather*}
A_0 =
\begin{pmatrix}
0 & 1 & 0 \\
0 & 0 & 0 \\
0 & 0 & 0
\end{pmatrix},\qquad
A_1=
\begin{pmatrix}
p_2q_2+1 & 0 & \frac{q_1q_2}{t} \\
0 & p_1q_1-p_2q_2 & t \\
1 & -1/q_1 & -p_1q_1
\end{pmatrix},\\
A_2 =
\begin{pmatrix}
0 & 0 & 0 \\
-t/q_2 & 0 & 0 \\
0 & 0 & 0
\end{pmatrix},\qquad
B_0=\frac{1}{t}
\begin{pmatrix}
p_2q_2 & 0 & 0 \\
0 & p_1q_1-p_2q_2 & t \\
1 & 0 & -p_1q_1-1
\end{pmatrix}.
\end{gather*}

The Hamiltonian is given by
\begin{gather*}
tH_{\mathrm{KFS}}^{\frac43+\frac43}\left(t; {q_1,p_1 \atop q_2,p_2}\right) =
tH_{\mathrm{III}(D_8)} (t; q_1,p_1 ) \\
\qquad{}
  +tH_{\mathrm{III}(D_8)} (t; q_2,p_2 )-p_1q_1p_2q_2+\left(\frac{q_1q_2}{t}+q_1+q_2\right)\nonumber.
\end{gather*}

\section{Laplace transform}\label{sec:Laplace}
In the degeneration scheme of the Garnier system and the Fuji--Suzuki system, we can see that the same Hamiltonian appears in several places. The linear systems corresponding to the same Hamiltonian can be transformed into one another by the Laplace transform. In this section, we present the correspondences through the Laplace transform.

It is known that a linear system (\ref{eq:rational_LDE}) with $r_\infty \le 1$ can canonically be written in the following form~\cite{Y1}:
\begin{gather}\label{eq:QP}
\frac{{\rm d}Y}{{\rm d}x}=\big(Q(x-T)^{-1}P+S\big)Y.
\end{gather}
This system transforms into
\begin{gather}\label{eq:PQ}
\frac{{\rm d}Y}{{\rm d}x}=\big({-}P(x-S)^{-1}Q-T\big)Y
\end{gather}
by the Laplace transform $x \mapsto -{\rm d}/{\rm d}x$, ${\rm d}/{\rm d}x \mapsto x$. The correspondence between (\ref{eq:QP}) and (\ref{eq:PQ}) is known as the Harnad duality when both $T$ and $S$ are semisimple~\cite{Har}, while here we do not impose the semisimplicity of $T$ or $S$.

Using this, we have the following correspondences of linear systems:
\begin{alignat*}{3}
& H^{2+2+1}_{\mathrm{Gar}}\colon \quad &&  \stackrel{\infty}{(1)(1)}, (1)(1), 11 \leftrightarrow \,
\stackrel{\infty}{(1)_2(1)}, 21, 21, &\\
& H^{3+2}_{\mathrm{Gar}}\colon \quad && ((1))((1)), \stackrel{\infty}{(1)(1)} \, \leftrightarrow \,
\stackrel{\infty}{((1))_3}, 21, 21, &\\
& H^{2+\frac32+1}_{\mathrm{Gar}}\colon \quad && \stackrel{\infty}{(1)_2}, (1)(1), 11 \leftrightarrow \,
\stackrel{\infty}{(1)_2(1)}, (2)(1), &\\
& H^{3+\frac{3}{2}}_{\mathrm{Gar}}\colon \quad && ((1))((1)), \stackrel{\infty}{(1)_2} \,\leftrightarrow \,
((1))_3, (2)(1),& \\
& H^{A_3}_{\mathrm{FS}}\colon \quad && \stackrel{\infty}{(11)(1)}, (11)(1) \leftrightarrow \,
\stackrel{\infty}{(1)_2 1}, 21, 111,& \\
& H^{2+\frac32}_{\mathrm{Suz}}\colon \quad && \stackrel{\infty}{(1)_2 1}, (11)(1) \leftrightarrow \,
\stackrel{\infty}{(1)_3}, 21, 111, &\\
& H_\mathrm{KFS}^{\frac32+\frac32}\colon \quad && \stackrel{\infty}{(1)_2 1}, (1)_2 1 \leftrightarrow \,
\stackrel{\infty}{(1)_3}, (11)(1).&
\end{alignat*}
Here we indicated which spectral types correspond to $\infty$.

When $r_\infty$ is greater than one, the correspondence through the Laplace transform is somewhat complicated. However, by using the method described in~\cite{KNS}, we obtain the following correspondences:
\begin{alignat*}{3}
&H^{3+2}_{\mathrm{Gar}}\colon \quad && \stackrel{\infty}{((1))((1))}, (1)(1) \leftrightarrow \,
\stackrel{\infty}{((1))(1)}_2, 21, &\\
& H^{4+1}_{\mathrm{Gar}}\colon \quad && \stackrel{\infty}{(((1)))(((1)))}, 11 \leftrightarrow \,
\stackrel{\infty}{(((1)))_2 1}, 21, &\\
& H^{\frac{5}{2}+2}_{\mathrm{Gar}}\colon \quad && \stackrel{\infty}{(((1)))_2}, (1)(1) \leftrightarrow \,
\stackrel{\infty}{(((1)))(1)_2}, &\\
& H^{5}_{\mathrm{Gar}}\colon \quad && \stackrel{\infty}{((((1))))_2}, 21 \leftrightarrow \,
\stackrel{\infty}{((((1))))((((1))))}.&
\end{alignat*}

\section{Conclusion} \label{sec:conclusion}

The degeneration scheme presented in this series of papers focus on linear equations. Thus the scheme is redundant in terms of Hamiltonians since it happens that a certain Hamiltonian appears in several places.

When we focus on Hamiltonians, the degeneration scheme can be reduced to the following scheme. This scheme shows the relationship among 40 Painlev\'e-type equations including those already known.

\begin{table}[t!]

\rotatebox{90}{\begin{minipage}{\textheight}\vspace*{20mm}
\centering
\begingroup
\renewcommand{\arraystretch}{1.5}
\begin{xy}
{(0,-12) *{H_{\mathrm{Gar}}^{1+1+1+1+1}}},
{\ar (11,-12);(20,-12)},
{(30,-12) *{H_{\mathrm{Gar}}^{2+1+1+1}}},
{\ar (40,-12);(53,-2)},
{\ar (40,-12);(53,-12)},
{\ar (40,-12);(51,-22)},
{(60,-2) *{H_{\mathrm{Gar}}^{3+1+1}}},
{(60,-12) *{H_{\mathrm{Gar}}^{2+2+1}}},
{(60,-22) *{H_{\mathrm{Gar}}^{\frac32+1+1+1}}},
{\ar (67,-2);(85,3)},
{\ar (67,-2);(85,-7)},
{\ar (67,-2);(83,-17)},
{\ar (67,-12);(85,3)},
{\ar (67,-12);(85,-7)},
{\ar (67,-12);(83,-27)},
{\ar (68,-22);(83,-17)},
{\ar (68,-22);(83,-27)},
{(90,3) *{H_{\mathrm{Gar}}^{4+1}}},
{(90,-7) *{H_{\mathrm{Gar}}^{3+2}}},
{(90,-17) *{H_{\mathrm{Gar}}^{\frac52+1+1}}},
{(90,-27) *{H_{\mathrm{Gar}}^{2+\frac32+1}}},
{\ar (96,3);(115,8)},
{\ar (96,3);(115,-2)},
{\ar (96,-7);(115,8)},
{\ar (96,-7);(115,-12)},
{\ar (96,-7);(115,-22)},
{\ar (97,-17);(115,-2)},
{\ar (97,-17);(115,-22)},
{\ar (97,-27);(115,-2)},
{\ar (97,-27);(115,-12)},
{\ar (97,-27);(115,-22)},
{\ar (97,-27);(113,-32)},
{(120,8) *{H_{\mathrm{Gar}}^5}},
{(120,-2) *{H_{\mathrm{Gar}}^{\frac72+1}}},
{(120,-12) *{H_{\mathrm{Gar}}^{3+\frac32}}},
{(120,-22) *{H_{\mathrm{Gar}}^{\frac52+2}}},
{(120,-32) *{H_{\mathrm{Gar}}^{\frac32+\frac32+1}}},
{\ar (126,8);(145,-2)},
{\ar (126,-2);(145,-2)},
{\ar (126,-12);(145,-2)},
{\ar (126,-12);(145,-22)},
{\ar (126,-22);(145,-2)},
{\ar (126,-22);(145,-22)},
{\ar (126,-32);(145,-22)},
{(150,-2) *{H_{\mathrm{Gar}}^{\frac92}}},
{(150,-22) *{H_{\mathrm{Gar}}^{\frac52+\frac32}}},
{(0,-42) *{H_{\mathrm{FS}}^{A_5}}},
{\ar (5,-42);(25,-42)},
{\ar (5,-42);(25,-52)},
{\ar (5,-42);(20,-12)},
{(30,-42) *{H_{\mathrm{FS}}^{A_4}}},
{\ar (35,-42);(55,-42)},
{\ar (35,-42);(55,-52)},
{\ar (35,-42);(53,-2)},
{(60,-42) *{H_{\mathrm{FS}}^{A_3}}},
{\ar (65,-42);(85,-42)},
{\ar (65,-42);(83,-17)},
{\ar (65,-42);(85,-7)},
{(90,-42) *{H_{\mathrm{Suz}}^{2+\frac32}}},
{\ar (95,-42);(115,-42)},
{\ar (95,-42);(115,8)},
{(120,-42) *{H_{\mathrm{KFS}}^{\frac32+\frac32}}},
{\ar (125,-42);(145,-42)},
{(150,-42) *{H_{\mathrm{KFS}}^{\frac32+\frac43}}},
{\ar (155,-42);(175,-42)},
{(180,-42) *{H_{\mathrm{KFS}}^{\frac43+\frac43}}},
{(30,-52) *{H_{\mathrm{NY}}^{A_5}}},
{\ar (35,-52);(55,-52)},
{\ar (35,-52);(51,-22)},
{(60,-52) *{H_{\mathrm{NY}}^{A_4}}},
{\ar (65,-52);(83,-17)},
{(0,-62) *{H_{\mathrm{Ss}}^{D_6}}},
{\ar (5,-62);(25,-62)},
{\ar (5,-62);(25,-52)},
{(30,-62) *{H_{\mathrm{Ss}}^{D_5}}},
{\ar (35,-62);(55,-62)},
{\ar (35,-62);(55,-52)},
{(60,-62) *{H_{\mathrm{Ss}}^{D_4}}},
{\ar (65,-62);(85,-62)},
{\ar (65,-62);(83,-17)},
{(90,-62) *{H_{\mathrm{KSs}}^{2+\frac32}}},
{\ar (95,-62);(115,-62)},
{(120,-62) *{H_{\mathrm{KSs}}^{2+\frac43}}},
{\ar (125,-62);(145,-62)},
{(150,-62) *{H_{\mathrm{KSs}}^{2+\frac54}}},
{\ar (155,-62);(175,-62)},
{(180,-62) *{H_{\mathrm{KSs}}^{\frac32+\frac54}}},
{(30,-100) *{H^{\mathrm{Mat}}_{\mathrm{VI}}}},
{\ar (35,-100);(55,-100)},
{(60,-100) *{H^{\mathrm{Mat}}_{\mathrm{V}}}},
{\ar (65,-100);(83,-90)},
{\ar (65,-100);(85,-110)},
{(90,-90) *{H^{\mathrm{Mat}}_{\mathrm{III}(D_6)}}},
{\ar (97,-90);(113,-90)},
{\ar (97,-90);(115,-110)},
{(90,-110) *{H^{\mathrm{Mat}}_{\mathrm{IV}}}},
{\ar (95,-110);(115,-110)},
{(120,-90) *{H^{\mathrm{Mat}}_{\mathrm{III}(D_7)}}},
{\ar (127,-90);(143,-90)},
{\ar (127,-90);(145,-110)},
{(120,-110) *{H^{\mathrm{Mat}}_{\mathrm{II}}}},
{\ar (125,-110);(145,-110)},
{(150,-90) *{H^{\mathrm{Mat}}_{\mathrm{III}(D_8)}}},
{(150,-110) *{H^{\mathrm{Mat}}_{\mathrm{I}}}}
\end{xy}
\endgroup
\end{minipage}}
\end{table}

\begin{Remark}
To determine whether the Hamiltonians given dif\/ferent names in this series of papers are actually dif\/ferent or not requires further  consideration. Concerning this problem, we refer to~\cite{N}, which is an attempt to characterize the four-dimensional Painlev\'e-type equations from an algebro-geometric point of view.
\end{Remark}

\begin{Remark}
We think that there are no degenerations of linear systems other than what we considered in this series of papers, and hence we believe that there are no other four-dimensional Painlev\'e-type equations. However, further research is needed to show that the four-dimensional Painlev\'e-type equations obtained in this series of papers actually constitute a complete list. We think that a way to classify the unramif\/ied linear equations with four accessory parameters gives a hint (see~\cite{HO}).
\end{Remark}

\appendix

\section{Data on degenerations}
In this appendix, we give explicit coordinate transformations used in degenerations. For examp\-le, for the linear system associated with $H_{\mathrm{Gar}}^{2+1+1+1}$, changing the variables and parameters $q_i$, $p_i$, $t_i$, $\theta^*_i$, $x$, $Y$ as shown in the table and taking the limit $\varepsilon \to 0$, we obtain the linear system associated with $H_{\mathrm{Gar}}^{\frac32+1+1+1}$. Also at this time, $\lim_{\varepsilon \to 0}\tilde{H}_i$ gives the Hamiltonian $H_{\mathrm{Gar}, t_i}^{\frac32+1+1+1}$.

The data that do not appear in the table below do not need to be changed.

\subsection{degenerations of the Garnier system}
$\boldsymbol{2+1+1+1 \to 3/2+1+1+1}$
\begin{gather*}
\theta^\infty_1=\tilde{\theta}^\infty_1+\varepsilon^{-1}, \qquad \theta^\infty_2=-\varepsilon^{-1}, \qquad
t_i=\varepsilon\tilde{t}_i,\qquad H_{t_i}=\varepsilon^{-1}\tilde{H}_{i}, \qquad i=1,2,\\
q_1=-\frac{1}{\varepsilon\tilde{q}_1},\qquad p_1=\varepsilon\tilde{q}_1\big(\tilde{p}_1\tilde{q}_1-\theta^1\big),\qquad
q_2=-\frac{1}{\varepsilon\tilde{q}_2},\qquad p_2=\varepsilon\tilde{q}_2\big(\tilde{p}_2\tilde{q}_2-\theta^t\big),\\
x=\frac{\tilde{x}}{\tilde{t}_1},\qquad
Y={t_1}^{\tilde{\theta}^\infty_1+\varepsilon^{-1}}
\begin{pmatrix}
1 & \\
 & u
\end{pmatrix}^{-1}
\begin{pmatrix}
0 & 1 \\
\varepsilon & 1
\end{pmatrix}
\tilde{Y}.
\end{gather*}

\noindent
$\boldsymbol{3+1+1 \to 5/2+1+1}$
\begin{gather*}
 \theta^\infty_1=\tilde{\theta}^\infty_1-\varepsilon^{-6},\qquad
\theta^\infty_2=\varepsilon^{-6},\qquad
t_i=\varepsilon^{-3}(\varepsilon^4 \tilde{t}_i-2),\qquad H_{t_i}=\varepsilon^{-1}\tilde{H}_{i},\\
 q_i=\varepsilon^{-3}\big(\varepsilon^2 \tilde{q}_i+1\big),\qquad p_i=\varepsilon \tilde{p}_i, \qquad i=1,2,\\
 x=-\varepsilon(\tilde{x}-\tilde{t}_2),\qquad
Y=e^{\varepsilon^{-2}(x-\tilde{t}_2)}
\begin{pmatrix}
\varepsilon^{-1} & -\varepsilon^{-3} \\
0 & -1
\end{pmatrix}\tilde{Y}.
\end{gather*}

\noindent
$\boldsymbol{2+2+1 \to 2+3/2+1}$
\begin{gather*}
\theta^\infty_1=\tilde{\theta}^\infty_1+\varepsilon^{-1},\qquad \theta^\infty_2=-\varepsilon^{-1},\qquad
t_1=\varepsilon\tilde{t}_1,\qquad t_2=-\varepsilon\tilde{t}_2,\\
 H_{t_1}=\varepsilon^{-1}\tilde{H}_{1},\qquad
H_{t_2}=-\varepsilon^{-1}\left( \tilde{H}_{2}-\frac{\tilde{p}_2\tilde{q}_2}{\tilde{t}_2} \right),\\
 q_1=-\frac{1}{\varepsilon \tilde{q}_1},\qquad p_1=\varepsilon \tilde{q}_1\big(\tilde{p}_1\tilde{q}_1-\theta^1\big),\qquad
q_2=\varepsilon\tilde{t}_2\tilde{p}_2,\qquad p_2=-\frac{\tilde{q}_2}{\varepsilon\tilde{t}_2},\\
 x=\frac{\tilde{x}}{\tilde{t}_1},\qquad
Y={t_1}^{\tilde{\theta}^\infty_1+\varepsilon^{-1}}
\begin{pmatrix}
1 & \\
 & u
\end{pmatrix}^{-1}
\begin{pmatrix}
1 & -\varepsilon^{-1} \\
0 & -\varepsilon^{-1}
\end{pmatrix}
\tilde{Y}.
\end{gather*}

\noindent
$\boldsymbol{3/2+1+1+1 \to 5/2+1+1}$
\begin{gather*}
 \theta^0=-2\varepsilon^{-3},\qquad \theta^\infty_1=\tilde{\theta}^\infty_1+2\varepsilon^{-3},\qquad
t_i=-\varepsilon^{-6}\big(\varepsilon^2 \tilde{t}_i+1\big),\qquad H_{t_i}=-\varepsilon^{4}\tilde{H}_{i},\\
 q_i=\varepsilon^{-2}\tilde{q}_i+\varepsilon^{-3},\qquad p_i=\varepsilon^2 \tilde{p}_i, \qquad i=1,2,\\
 x=-\varepsilon^{-6}\big(\varepsilon^2 \tilde{x}+1\big),\qquad
Y=e^{-\varepsilon^{-1}x}
\begin{pmatrix}
1 & -\varepsilon^{-1} \\
0 & \varepsilon^2
\end{pmatrix}
\tilde{Y}.
\end{gather*}

\noindent
$\boldsymbol{3/2+1+1+1 \to 2+3/2+1}$
\begin{gather*}
 \theta^0=-\varepsilon^{-1},\qquad \theta^t=\tilde{\theta}^0+\varepsilon^{-1},\qquad t_2=-\varepsilon\tilde{t}_2, \\
 q_2=\varepsilon\tilde{t}_2\tilde{p}_2,\qquad p_2=-\frac{\tilde{q}_2}{\varepsilon\tilde{t}_2},\qquad
H_{t_2}=-\varepsilon^{-1}\left( \tilde{H}_{2}-\frac{\tilde{p}_2\tilde{q}_2}{\tilde{t}_2} \right).
\end{gather*}

\noindent
$\boldsymbol{4+1 \to 7/2+1}$
\begin{gather*}
 \theta^\infty_1=\tilde{\theta}^\infty_1-2\varepsilon^{-15},\qquad \theta^\infty_2=2\varepsilon^{-15},\qquad
t_1=\varepsilon^2\tilde{t}_1-3\varepsilon^{-10},\qquad t_2=-\varepsilon\tilde{t}_2+\varepsilon^{-5},\\
 H_{t_1}=\varepsilon^{-2}\left( \tilde{H}_{1}-\frac{\varepsilon^3}{2}\tilde{q}_1 \right),\qquad
H_{t_2}=-\varepsilon^{-1}\tilde{H}_{2},\\
 q_1=\varepsilon \tilde{q}_1+\varepsilon^{-5},\qquad
p_1=\frac{\varepsilon^2}{2}\big(\tilde{q}_1{}^2+\tilde{t}_1\big)+\varepsilon^{-1}\tilde{p}_1+\varepsilon^{-4}\tilde{q}_1-\varepsilon^{-10},\\
 q_2=\varepsilon^{-2}\tilde{q}_2+\varepsilon^{-5},\qquad p_2=\varepsilon^2\tilde{p}_2,\\
 x=\varepsilon\tilde{x},\qquad
Y=\mathrm{exp}\left( -\frac{1}{2\varepsilon^3}(\tilde{x}-\tilde{t}_2)^2+\varepsilon^{-9}(\tilde{x}-\tilde{t}_2)-\varepsilon^{-3}\tilde{t}_1 \right)
\begin{pmatrix}
1 & \varepsilon^{-3} \\
0 & \varepsilon^2
\end{pmatrix}
\tilde{Y}.
\end{gather*}

\noindent
$\boldsymbol{3+2 \to 3+3/2}$
\begin{gather*}
 \theta^0=2\varepsilon^{-1},\qquad
\theta^\infty_1=\tilde{\theta}^\infty_1-\varepsilon^{-1},\qquad \theta^\infty_2=\tilde{\theta}^\infty_2-\varepsilon^{-1},\qquad
t_1=\sqrt{-1}\varepsilon\tilde{t}_1,\qquad t_2=\sqrt{-1}\tilde{t}_2,\\
 H_{t_1}=\frac{1}{\sqrt{-1}\varepsilon}\left( \tilde{H}_1
-\frac{\tilde{p}_1\tilde{q}_1+\tilde{p}_2\tilde{q}_2}{\tilde{t}_1}+\frac{1}{\tilde{t}_1\tilde{p}_2} \right),\qquad
H_{t_2}=\frac{1}{\sqrt{-1}}\big(\tilde{H}_2+\tilde{p}_2\tilde{q}_1\big),
\\
q_1=\frac{\tilde{t}_1}{f},\qquad
p_1=-\frac{f}{\tilde{t}_1}\left( \tilde{p}_1\tilde{q}_1+\tilde{p}_2\tilde{q}_2+\tilde{\theta}^\infty_1+\frac{1}{\varepsilon}-\frac{1}{\tilde{p}_2} \right),\\
q_2=\sqrt{-1}f\tilde{p}_2,\qquad
p_2=\sqrt{-1}\left( \frac{1}{f}\left(\tilde{q}_2+\frac{\tilde{p}_1\tilde{q}_1+2\tilde{\theta}^\infty_1}{\tilde{p}_2}-\frac{1}{\tilde{p}_2{}^2} \right)+\tilde{p}_2\tilde{q}_1 \right),\\
x=\sqrt{-1}\tilde{x},\qquad
Y=x^{\varepsilon^{-1}}
\begin{pmatrix}
1 & 0 \\
0 & u
\end{pmatrix}^{-1}\,
\begin{pmatrix}
1 & 0 \\
0 & \frac{1}{\sqrt{-1}\tilde{p}_2{}^2f}
\end{pmatrix}
\begin{pmatrix}
1 & 0 \\
0 & \tilde{u}
\end{pmatrix}
\tilde{Y},
\end{gather*}
where
\begin{gather*}
f=\tilde{q}_1-\frac{\tilde{t}_2}{\tilde{p}_2}+\frac{\tilde{p}_1}{\tilde{p}_2{}^2}.
\end{gather*}
The new gauge parameter $\tilde{u}$ satisf\/ies the following equations:
\begin{gather*}
\frac{1}{\tilde{u}}\frac{\partial \tilde{u}}{\partial \tilde{t}_1}=-\frac{2}{\tilde{t}_1}\big(\tilde{p}_1\tilde{q}_1+\tilde{p}_2\tilde{q}_2+\tilde{\theta}^\infty_1\big),\qquad
\frac{1}{\tilde{u}}\frac{\partial \tilde{u}}{\partial \tilde{t}_2}=-2\tilde{p}_2\tilde{q}_1.
\end{gather*}

\noindent
$\boldsymbol{3+2 \to 5/2+2}$
\begin{gather*}
 \theta^\infty_1=\tilde{\theta}^\infty_1-\varepsilon^{-6},\qquad \theta^\infty_2=\varepsilon^{-6},\qquad
t_1=\varepsilon \tilde{t}_2,\ t_2=\varepsilon \tilde{t}_1-2\varepsilon^{-3},\\
 H_{t_1}=\varepsilon^{-1}\tilde{H}_{2}-\frac{\tilde{p}_2\tilde{q}_2}{\varepsilon\tilde{t}_2},\qquad
H_{t_2}=\varepsilon^{-1}\tilde{H}_{1}-\varepsilon\tilde{p}_1,\qquad
 q_1=\varepsilon^2\tilde{t}_2\tilde{p}_2,\qquad p_1=-\frac{\tilde{q}_2}{\varepsilon^2\tilde{t}_2},\\
q_2=\varepsilon^{-3}+\varepsilon^{-1}\tilde{q}_1+\varepsilon\big(\tilde{p}_1-\tilde{t}_1\big),\qquad p_2=\varepsilon\tilde{p}_1,\\
x=\varepsilon\tilde{x},\qquad Y=
e^{\varepsilon^{-2}(\tilde{x}-\tilde{t}_1)}
\begin{pmatrix}
0 & 1 \\
\varepsilon^5 & 1/q_2
\end{pmatrix}
\tilde{Y}.
\end{gather*}

\noindent
$\boldsymbol{5/2+1+1 \to 7/2+1}$
\begin{gather*}
\theta^{t_1}=\tilde{\theta}^0,\qquad \theta^{t_2}=-2\varepsilon^{-15},\qquad \theta^\infty_1=\tilde{\theta}^\infty_1+2\varepsilon^{-15},\\
t_1=\varepsilon^{-4}\tilde{t}_2-2\varepsilon^{-10},\qquad t_2=\varepsilon^2\tilde{t}_1-3\varepsilon^{-10},\qquad
H_{t_1}=\varepsilon^4\tilde{H}_{2},\qquad H_{t_2}=\varepsilon^{-2}\left( \tilde{H}_{1}-\frac{\varepsilon^3}{2}\tilde{q}_1 \right),\\
 q_1=\varepsilon^{-2}\tilde{q}_2,\qquad p_1=\varepsilon^2\tilde{p}_2,\qquad
q_2=\varepsilon \tilde{q}_1+\varepsilon^{-5},\\
p_2=\frac{\varepsilon^2}{2}\big(\tilde{q}_1{}^2+\tilde{t}_1\big)+\varepsilon^{-1}\tilde{p}_1+\varepsilon^{-4}\tilde{q}_1-\varepsilon^{-10},\\
 x=-\varepsilon^{-4}\tilde{x}+\varepsilon^{-4}\tilde{t}_2-2\varepsilon^{-10},\qquad
Y=\mathrm{exp}\left( \varepsilon^{-5}x+\frac{\varepsilon^5}{2}x^2+\varepsilon^{-3}\tilde{t}_1 \right)
\begin{pmatrix}
1 & 0 \\
0 & -\varepsilon^2
\end{pmatrix}
\tilde{Y}.
\end{gather*}

\noindent
$\boldsymbol{5/2+1+1 \to 5/2+2}$
\begin{gather*}
\theta^{t_1}=\varepsilon^{-1},\!\qquad \theta^{t_2}=\tilde{\theta}^0-\varepsilon^{-1},\!\qquad
t_2=\tilde{t}_1+\varepsilon \tilde{t}_2,\!\qquad
H_{t_1}=\tilde{H}_{1}-\varepsilon^{-1}\tilde{H}_{2},\!\qquad
H_{t_2}=\varepsilon^{-1}\tilde{H}_{2},\\
q_1=\tilde{q}_1,\qquad p_1=\tilde{p}_1-\varepsilon^{-1}\tilde{p}_2,\qquad
q_2=\tilde{q}_1+\varepsilon\tilde{q}_2,\qquad p_2=\varepsilon^{-1}\tilde{p}_2,\\
x=-\tilde{x}+t_1,\qquad
Y=
\begin{pmatrix}
1 & -q_1 \\
0 & 1
\end{pmatrix}
\tilde{Y}.
\end{gather*}

\noindent
$\boldsymbol{2+3/2+1 \to 7/2+1}$
\begin{gather*}
 \theta^0=-3\varepsilon^{-5},\qquad \theta^\infty_1=\tilde{\theta}^\infty_1+3\varepsilon^{-5},\qquad
  \theta^{t_1}=\tilde{\theta}^0,\\
t_1=-\varepsilon^{-8}\tilde{t}_2-\varepsilon^{-10},\qquad t_2=2\varepsilon^{-11}\tilde{t}_1+2\varepsilon^{-15},\\
H_{t_1}=-\varepsilon^8\left( \tilde{H}_{2}
+\frac{\varepsilon^2 \tilde{p}_2\tilde{q}_2+\varepsilon\tilde{p}_2}{1+\varepsilon^2 \tilde{t}_2} \right),\qquad
H_{t_2}=\frac{\varepsilon^{11}}{2}\tilde{H}_{1},\\
q_1=\varepsilon^{-5}\big(1+\varepsilon^2 \tilde{t}_2\big)(1+\varepsilon \tilde{q}_2),\qquad
p_1=\frac{\varepsilon^4\tilde{p}_2}{1+\varepsilon^2\tilde{t}_2},\\
q_2=\varepsilon^{-10}\big(1-\varepsilon^2 \tilde{q}_1\big),\qquad
p_2=\frac{\varepsilon^5}{2}-\frac{\varepsilon^7\tilde{q}_1}{2}-\varepsilon^8\tilde{p}_1,\\
x=t_1\big(1-\varepsilon^2\tilde{x}\big),\qquad
Y=\mathrm{exp}\left[ \frac{\tilde{x}-\tilde{t}_2}{2\varepsilon^3}-\frac{\tilde{t}_1}{\varepsilon}-\frac{1}{\varepsilon}
\left( \frac{(\tilde{x}-\tilde{t}_2)^2}{4}-\frac{\tilde{t}_2\tilde{x}}{2} \right) \right]
\begin{pmatrix}
1 & \varepsilon^{-1} \\
0 & -\frac{1}{\varepsilon^6 t_1}
\end{pmatrix}
\tilde{Y}.
\end{gather*}

\noindent
$\boldsymbol{2+3/2+1 \to 3+3/2}$
\begin{gather*}
\theta^0=\tilde{\theta}^\infty_2-\tilde{\theta}^\infty_1-\frac{1}{\varepsilon^2},\qquad \theta^{t_1}=\frac{1}{\varepsilon^2},\qquad
\theta^\infty_1=2\tilde{\theta}^\infty_1,\qquad
t_1=\varepsilon \tilde{t}_1,\qquad t_2=\tilde{t}_1\tilde{t}_2+\frac{\tilde{t}_1}{\varepsilon},\\
 H_{t_1}=\frac{1}{\varepsilon}\left( \tilde{H}_1-\frac{1}{\tilde{t}_1}\left( \tilde{t}_2+\frac{1}{\varepsilon} \right)\tilde{H}_2
-\frac{\tilde{p}_2\tilde{q}_2}{\tilde{t}_1}-\frac{\tilde{p}_2\tilde{q}_1}{\varepsilon\tilde{t}_1} \right),\qquad
H_{t_2}=\frac{1}{\tilde{t}_1}\tilde{H}_2,\\
 q_1=\varepsilon\tilde{q}_1,\ p_1=\frac{1}{\varepsilon}\left( \tilde{p}_1-\frac{\tilde{p}_2}{\varepsilon} \right),\qquad
q_2=\frac{\tilde{t}_1\tilde{p}_2}{\varepsilon},\qquad
p_2=-\frac{1}{\tilde{t}_1}(\varepsilon \tilde{q}_2+\tilde{q}_1),\\
 x=\frac{\tilde{t}_1}{\tilde{x}},\qquad
Y=\tilde{x}^{\tilde{\theta}^\infty_1}v
\begin{pmatrix}
-1/\tilde{p}_2 & 0 \\
-1 & 1/\tilde{p}_2
\end{pmatrix}
\begin{pmatrix}
1 & 0 \\
0 & \tilde{u}
\end{pmatrix}
\tilde{Y}.
\end{gather*}
Here $\tilde{u}$ and $v$ satisf\/ies
\begin{gather*}
 \frac{1}{\tilde{u}}\frac{\partial \tilde{u}}{\partial \tilde{t}_1}=-\frac{2}{\tilde{t}_1}\big(\tilde{p}_1\tilde{q}_1+\tilde{p}_2\tilde{q}_2+\tilde{\theta}^\infty_1\big),\qquad
\frac{1}{\tilde{u}}\frac{\partial \tilde{u}}{\partial \tilde{t}_2}=-2\tilde{p}_2\tilde{q}_1,\\
 \frac{1}{v}\frac{\partial v}{\partial \tilde{t}_1}=-\frac{1}{\tilde{t}_1}\left(\tilde{p}_1\tilde{q}_1+\tilde{p}_2\tilde{q}_2+2\tilde{\theta}^\infty_1-\frac{1}{\tilde{p}_2}\right),\qquad
\frac{1}{v}\frac{\partial v}{\partial \tilde{t}_2}=-\frac{\tilde{p}_1}{\tilde{p}_2}-\tilde{p}_2\tilde{q}_1+\tilde{t}_2.
\end{gather*}

\noindent
$\boldsymbol{2+3/2+1 \to 5/2+2}$
\begin{gather*}
 \theta^{t_1}=2\varepsilon^{-3},\qquad \theta^\infty_1=\tilde{\theta}^\infty_1-2\varepsilon^{-3},\qquad
t_1=\varepsilon^{-4}\tilde{t}_1+\varepsilon^{-6},\qquad t_2=-\varepsilon^{-4}\tilde{t}_2,\\
 H_{t_1}=\varepsilon^4\tilde{H}_{1},\qquad
H_{t_2}=-\varepsilon^4\left( \tilde{H}_{2}-\frac{\tilde{p}_2\tilde{q}_2}{\tilde{t}_2} \right),\\
 q_1=\varepsilon^{-3}-\varepsilon^{-2}\tilde{q}_1,\qquad
p_1=1-\varepsilon^2\tilde{p}_1,\qquad
q_2=-\varepsilon^{-2}\tilde{t}_2\tilde{p}_2,\qquad p_2=\varepsilon^2\frac{\tilde{q}_2}{\tilde{t}_2},\\
x=\varepsilon^{-4}\tilde{x},\qquad
Y=e^{\varepsilon^{-1}(\tilde{t}_1-\tilde{x})}
\begin{pmatrix}
1 & 0 \\
0 & -\varepsilon^8 t_1
\end{pmatrix}
\tilde{Y}.
\end{gather*}

\noindent
$\boldsymbol{2+3/2+1 \to 3/2+3/2+1}$
\begin{gather*}
 \theta^0=-2\varepsilon^{-1},\qquad \theta^\infty_1=\tilde{\theta}^\infty_1+2\varepsilon^{-1},\qquad
t_2=\varepsilon \tilde{t}_2,\qquad
H_{t_2}=\varepsilon^{-1}\tilde{H}_{2},\\
q_2=-\tilde{q}_2,\qquad p_2=-\tilde{p}_2-\frac{1}{\varepsilon \tilde{q}_2},\qquad
 Y=x^{-\varepsilon^{-1}}\tilde{Y}.
\end{gather*}

\noindent
$\boldsymbol{5 \to 9/2}$
\begin{gather*}
 \theta^\infty_1=-\frac{5}{\varepsilon^{28}},\qquad \theta^\infty_2=\frac{5}{\varepsilon^{28}},\qquad
t_1=\frac{2\tilde{t}_1}{3\varepsilon^5}-\varepsilon^3\tilde{t}_2+\frac{160}{27\varepsilon^{21}},\qquad
t_2=\varepsilon^2\tilde{t}_1-\frac{5}{3\varepsilon^{14}},\\
 H_{t_1}=-\frac{1}{\varepsilon^3}\tilde{H}_2,\qquad
H_{t_2}=\frac{1}{\varepsilon^2}
\left(
\tilde{H}_1+\frac{2}{3\varepsilon^8}\tilde{H}_2+\frac{\tilde{p}_1}{\varepsilon^4}
\right),\\
 q_1=\varepsilon^5\tilde{q}_2-\varepsilon\tilde{p}_1-\frac{4}{3\varepsilon^7},\qquad
p_1=\frac{1}{\varepsilon^5}\big(\tilde{p}_1{}^2+\tilde{p}_2\big)+\frac{\tilde{q}_2}{\varepsilon^9}-\frac{\tilde{p}_1}{3\varepsilon^{13}}
+\frac{4}{3\varepsilon^{21}},\\
 q_2=-\frac{\tilde{p}_1}{\varepsilon^6}+\frac{1}{\varepsilon^{14}},\qquad
p_2=(\tilde{q}_1-2\tilde{p}_1\tilde{q}_2)\varepsilon^6+(-\tilde{p}_2+\tilde{t}_1)\varepsilon^2
-\frac{2\tilde{q}_2}{3\varepsilon^2}+\frac{2\tilde{p}_1}{3\varepsilon^6}+\frac{4}{9\varepsilon^{14}},\\
 x=\varepsilon\tilde{x}+\frac{2}{3\varepsilon^7},\qquad
Y=\exp\left[
\frac{\tilde{t}_1\tilde{x}+\frac{\tilde{x}{}^3}{3}-\tilde{t}_2}{\varepsilon^4}-\frac{\frac{\tilde{x}{}^2}{2}+\tilde{t}_1}{\varepsilon^{12}}+\frac{2\tilde{x}}{\varepsilon^{20}}
\right]
\begin{pmatrix}
1 & 0 \\
0 & u
\end{pmatrix}^{-1}
\begin{pmatrix}
1 & -\varepsilon^{-4} \\
0 & -\varepsilon^3
\end{pmatrix}
\tilde{Y}.
\end{gather*}

\noindent
$\boldsymbol{7/2+1 \to 9/2}$
\begin{gather*}
 \theta^0=-\frac{2}{\varepsilon^{35}},\qquad \theta^\infty_1=\frac{2}{\varepsilon^{35}},\qquad
t_1=\frac{\tilde{t}_1}{\varepsilon^8}+\frac{5}{3\varepsilon^{28}},\qquad
t_2=\varepsilon^6\tilde{t}_1-\varepsilon^{16}\tilde{t}_2-\frac{5}{3\varepsilon^{14}},\\
 H_{t_1}=\varepsilon^8\tilde{H}_1+\frac{\tilde{H}_2}{\varepsilon^2},\qquad
H_{t_2}=-\frac{\tilde{H}_2}{\varepsilon^{16}},\\
 q_1=\frac{\tilde{p}_1}{\varepsilon^4}-\frac{1}{3\varepsilon^{14}},\qquad
p_1=-\varepsilon^4\tilde{q}_1-\frac{\tilde{p}_1{}^2+\tilde{p}_2}{\varepsilon}+\frac{\tilde{q}_2}{\varepsilon^6},\\
 q_2=-\varepsilon^8\tilde{q}_2-\varepsilon^3\tilde{p}_1-\frac{1}{\varepsilon^7},\qquad
p_2=-\frac{\tilde{p}_2}{\varepsilon^8}-\frac{\tilde{p}_1}{\varepsilon^{18}}+\frac{1}{\varepsilon^{28}},\\
 x=\frac{1}{\varepsilon^{14}}\big(\varepsilon^{10}\tilde{x}-1\big),\qquad
Y=\exp\left[
\frac{\tilde{x}{}^3}{3\varepsilon^5}+\frac{\tilde{x}{}^2}{2\varepsilon^{15}}+\frac{\tilde{x}}{\varepsilon^{25}}\right]
\begin{pmatrix}
1 & 0 \\
0 & -\varepsilon^2
\end{pmatrix}\tilde{Y}.
\end{gather*}

\noindent
$\boldsymbol{3+3/2 \to 9/2}$
\begin{gather*}
 \theta^\infty_1=-\frac{15}{8\varepsilon^{14}},\qquad \theta^\infty_2=\frac{15}{8\varepsilon^{14}},\qquad
t_1=-\frac{\sqrt{2}}{\varepsilon^{35}}\big(1+\varepsilon^8 \tilde{t}_1+\varepsilon^{12}\tilde{t}_2\big),\qquad
t_2=-\sqrt{2}\varepsilon \tilde{t}_1+\frac{5}{\sqrt{2}\varepsilon^7},\\
 H_{t_1}=-\frac{\varepsilon^{23}}{\sqrt{2}}
\left(\tilde{H}_2
+\varepsilon^{12}\frac{\tilde{p}_1Q_1+(\tilde{p}_2-\tilde{t}_1)Q_2}{1+\varepsilon^8 \tilde{t}_1+\varepsilon^{12}\tilde{t}_2}\right),\qquad
H_{t_2}=-\frac{1}{\sqrt{2}\varepsilon}
\left( \tilde{H}_1-\frac{\tilde{H}_2}{\varepsilon^4}-\frac{5\tilde{p}_1}{8\varepsilon^2}+\tilde{q}_2 \right),\\
 q_1=\sqrt{2}\big(1+\varepsilon^8 \tilde{t}_1+\varepsilon^{12}\tilde{t}_2\big)
\left(\frac{Q_1}{\varepsilon^{11}}-\frac{Q_2}{\varepsilon^{15}}\right),\qquad
p_1=\frac{\varepsilon^{11}\tilde{p}_1}{\sqrt{2}(1+\varepsilon^8 \tilde{t}_1+\varepsilon^{12}\tilde{t}_2)}-\sqrt{2}\varepsilon^7,\\
 q_2=\frac{2(1+\varepsilon^8 \tilde{t}_1+\varepsilon^{12}\tilde{t}_2)Q_2}{\varepsilon^{22}},\qquad
p_2=\frac{\varepsilon^{18}\tilde{p}_1+\varepsilon^{22}(\tilde{p}_2-\tilde{t}_1)}{2(1+\varepsilon^8 \tilde{t}_1+\varepsilon^{12}\tilde{t}_2)}-\frac{\varepsilon^{14}}{2},\\
 x=\frac{\sqrt{2}}{\varepsilon^7(1-\varepsilon^4 \tilde{x})},\\
 Y=\exp\left[\frac{3\tilde{x}}{2\varepsilon^{10}}+\frac{\tilde{x}{}^2+\frac{\tilde{t}_1}{2}}{\varepsilon^6}
+\frac{-\tilde{t}_1\tilde{x}+\frac{\tilde{x}{}^3}{2}+\tilde{t}_2}{\varepsilon^2}\right]
\begin{pmatrix}
1 & 0 \\
0 & u
\end{pmatrix}^{-1}
\begin{pmatrix}
1 & \varepsilon^{-2} \\
p_2q_1 & q_1\left( \frac{p_2}{\varepsilon^2}-\frac{\sqrt{2}}{\varepsilon^{23}t_1} \right)
\end{pmatrix}\tilde{Y},
\end{gather*}
where
\begin{gather*}
Q_1=\tilde{q}_1+\frac{1}{8}
\left( \frac{\tilde{p}_1}{\varepsilon^6}-\frac{3\tilde{p}_1{}^2}{\varepsilon^2}
-\frac{3\tilde{p}_2}{\varepsilon^2}-\frac{2\tilde{t}_1}{\varepsilon^2}-\frac{3}{\varepsilon^{10}} \right),\qquad
Q_2=\tilde{q}_2-\frac{3\tilde{p}_1}{8\varepsilon^2}+\frac{1}{8\varepsilon^6}.
\end{gather*}

\noindent
$\boldsymbol{3+3/2 \to 5/2+3/2}$
\begin{gather*}
 \theta^\infty_1=\varepsilon^{-6},\qquad \theta^\infty_2=-\varepsilon^{-6},\qquad
t_1=\varepsilon\tilde{t}_2,\qquad t_2=\varepsilon\tilde{t}_1-2\varepsilon^{-3},\\
H_{t_1}=\varepsilon^{-1}\tilde{H}_2,\qquad
H_{t_2}=\varepsilon^{-1}\big(\tilde{H}_1+\tilde{q}_1\big),\\
 q_1=\varepsilon^{-1}\tilde{q}_2,\qquad
p_1=\frac{1}{\varepsilon \tilde{q}_2}\big(\tilde{p}_1-\tilde{q}_1{}^2-\tilde{t}_1+\varepsilon^2(\tilde{p}_2\tilde{q}_2-\tilde{p}_1\tilde{q}_1+\tilde{t}_1\tilde{q}_1)\big),\\
 q_2=\frac{\tilde{q}_2}{\varepsilon^4}\big(1-\varepsilon^2 \tilde{q}_1\big),\qquad
p_2=-\frac{1}{\varepsilon^2\tilde{q}_2}\big(1+\varepsilon^2\tilde{q}_1+\varepsilon^4(\tilde{p}_1-\tilde{t}_1)\big),\\
 x=\varepsilon\tilde{x},\qquad
Y=\exp\big(\varepsilon^{-2}(\tilde{t}_1-\tilde{x})\big)
\begin{pmatrix}
1 & 0 \\
0 & u
\end{pmatrix}^{-1}
\begin{pmatrix}
0 & \varepsilon q_1\\
1 & \varepsilon p_2q_1
\end{pmatrix}\tilde{Y}.
\end{gather*}

\noindent
$\boldsymbol{5/2+2 \to 9/2}$
\begin{gather*}
 \theta^0=-\frac{5}{\varepsilon^{21}},\qquad \theta^\infty_1=\frac{5}{\varepsilon^{21}},\qquad
t_1=\frac{\tilde{t}_1}{\varepsilon^2}+\varepsilon^4\tilde{t}_2-\frac{5}{\varepsilon^{14}},\qquad
t_2=-\frac{2\tilde{t}_1}{\varepsilon^{23}}-\frac{2}{\varepsilon^{35}},\\
 H_{t_1}=\frac{\tilde{H}_2}{\varepsilon^4},\\
 H_{t_2}=\frac{\varepsilon^{23}}{2}\left[
-\tilde{H}_1+\frac{\tilde{H}_2}{\varepsilon^6}+
\frac{\varepsilon^3}{2(1+\varepsilon^{12}\tilde{t}_1)}
\big\{
\tilde{p}_1{}^2-\tilde{p}_2-2\varepsilon^3\tilde{p}_1\tilde{q}_2+\varepsilon^6(\tilde{p}_2+\tilde{t}_1)\tilde{p}_1
-2\varepsilon^9\tilde{p}_2\tilde{q}_2
\big\}\right],\\
 q_1=\varepsilon^8\tilde{q}_1-\frac{\varepsilon^5}{2}\big(\tilde{p}_1{}^2+\tilde{p}_2+\tilde{t}_1\big)
-\varepsilon^2\tilde{q}_2+\frac{\tilde{p}_1}{\varepsilon}+\frac{1}{\varepsilon^7},\qquad
p_1=\frac{\tilde{p}_1}{\varepsilon^8}-\frac{2}{\varepsilon^{14}},\\
 q_2=\frac{2(1+\varepsilon^{12}\tilde{t}_1)}{\varepsilon^{19}}
\left(\tilde{q}_2-\frac{\tilde{p}_1}{2\varepsilon^3}+\frac{1}{2\varepsilon^9}\right),\qquad
p_2=\frac{\varepsilon^{19}}{2(1+\varepsilon^{12}\tilde{t}_1)}
\left(\tilde{p}_2-\tilde{t}_1+\frac{\tilde{p}_1}{\varepsilon^6}-\frac{1}{\varepsilon^{12}}\right),\\
 x=\frac{1}{\varepsilon^{14}}\big(\varepsilon^6\tilde{x}-1\big),\qquad
Y=\exp\left[
-\frac{\frac{\tilde{x}{}^3}{6}+\tilde{t}_1\tilde{x}}{\varepsilon^3}
+\frac{\frac{\tilde{x}{}^2}{4}-\tilde{t}_1}{\varepsilon^9}
+\frac{3\tilde{x}}{2\varepsilon^{15}}\right]
\begin{pmatrix}
1 & \varepsilon^{-3} \\
0 & \varepsilon^4
\end{pmatrix}
\tilde{Y}.
\end{gather*}

\noindent
$\boldsymbol{5/2+2 \to 5/2+3/2}$
\begin{gather*}
\theta^0=2\varepsilon^{-1},\qquad \theta^\infty_1=-2\varepsilon^{-1},\qquad
t_2=\varepsilon\tilde{t}_2,\qquad
H_{t_2}=\varepsilon^{-1}\tilde{H}_2,\\
q_2=\tilde{q}_2,\qquad p_2=\tilde{p}_2+\frac{1}{\varepsilon\tilde{q}_2},\qquad
Y=\tilde{x}^{\varepsilon^{-1}}\tilde{Y}.
\end{gather*}

\noindent
$\boldsymbol{3/2+3/2+1 \to 5/2+3/2}$
\begin{gather*}
\theta^{t_1}=2\varepsilon^{-3},\qquad \theta^\infty_1=-2\varepsilon^{-3},\qquad
t_1=\varepsilon^{-4}\tilde{t}_1+\varepsilon^{-6},\qquad t_2=\varepsilon^{-4}\tilde{t}_2,\\
H_{t_1}=\varepsilon^4\tilde{H}_1,\qquad H_{t_2}=\varepsilon^4\left( \tilde{H}_2-\frac{\tilde{p}_2\tilde{q}_2}{\tilde{t}_2} \right),\\
q_1=-\varepsilon^{-2}\tilde{q}_1+\varepsilon^{-3},\qquad p_1=1-\varepsilon^2\tilde{p}_1,\qquad
q_2=\frac{\tilde{t}_2}{\varepsilon^2\tilde{q}_2},\qquad p_2=-\frac{\varepsilon^2\tilde{p}_2\tilde{q}_2{}^2}{\tilde{t}_2},\\
x=\varepsilon^{-4}\tilde{x},\qquad
Y=\exp\big(\varepsilon^{-1}(\tilde{t}_1-\tilde{x})\big)
\begin{pmatrix}
1 & 0 \\
0 & -\varepsilon^2
\end{pmatrix}\tilde{Y}.
\end{gather*}

\subsection{degenerations of the Fuji--Suzuki system}\label{sec:deg_FS}

We f\/irst note that, concerning the following two linear systems, we adopt slightly dif\/ferent parametrizations from those in \cite{KNS} for convenience of calculation.

As for the linear system of the spectral type $(11)(1),21,111$, the irregular singular point is moved to $x=\infty$, so that the Riemann scheme is given by
\begin{gather*}
\left(
\begin{matrix}
 x=0 & x=1 & x=\infty \\
\begin{matrix} 0 \\ \theta^0_1 \\ \theta^0_2 \end{matrix}
& \begin{matrix} 0 \\ 0 \\ \theta^1 \end{matrix}
& \overbrace{\begin{matrix}
 t & \theta^\infty_1 \\
 0 & \theta^\infty_2 \\
 0 & \theta^\infty_3
 	\end{matrix}}
\end{matrix}
\right) ,
\end{gather*}
and the Fuchs--Hukuhara relation is written as $\theta^0_1+\theta^0_2+\theta^1+\theta^\infty_1+\theta^\infty_2+\theta^\infty_3=0$. The Lax pair is expressed as
\begin{gather*}
\frac{\partial Y}{\partial x} = \left( \frac{A_0}{x}+\frac{A_1}{x-1}+A_\infty \right)Y ,\qquad
\frac{\partial Y}{\partial  t} =\left( \frac{1}{t}A_\infty x+B_0\right)Y .
\end{gather*}
Here
\begin{gather*}
A_{\xi} =
U^{-1}P^{-1}\hat{A}_{\xi}PU, \qquad  \xi=0,1,\qquad U=\operatorname{diag}(1,u,v), \\
P=
\begin{pmatrix}
1 & 0 & 0 \\
p_2q_2-\theta^0_2-\theta^\infty_2 & 1 & 0 \\
p_1q_1-\theta^\infty_3 & \frac{p_1(q_1-q_2)-\theta^\infty_3}{\theta^\infty_3-\theta^\infty_2} & 1
\end{pmatrix}, \\
\hat{A}_0 =
\begin{pmatrix}
1 & 0 \\
0 & 1 \\
0 & 0
\end{pmatrix}
\begin{pmatrix}
\theta^0_1 & q_2-1 & q_1-1 \\
0 & \theta^0_2 & p_2(q_1-q_2)+\theta^0_2+\theta^\infty_2
\end{pmatrix},\\
\hat{A}_1 =
\begin{pmatrix}
1 \\
p_2 \\
p_1
\end{pmatrix}
\begin{pmatrix}
p_1q_1+p_2q_2+\theta^1 & -q_2 & -q_1
\end{pmatrix},\\
A_\infty=
\begin{pmatrix}
-t & 0 & 0 \\
0 & 0 & 0 \\
0 & 0 & 0
\end{pmatrix},\qquad
B_0=\frac{1}{t}
\begin{pmatrix}
0 & (A_0+A_1)_{12} & (A_0+A_1)_{13} \\
(A_0+A_1)_{21} & 0 & 0 \\
(A_0+A_1)_{31} & 0 & 0
\end{pmatrix}.
\end{gather*}

The Hamiltonian is given by
\begin{gather*}
 tH_{\mathrm{FS}}^{A_4}
 \left({\theta^1, \theta^\infty_3-\theta^\infty_2, \theta^0_1+\theta^\infty_2+1\atop
 \theta^\infty_2, \theta^0_2+\theta^\infty_2};t;
{q_1,p_1\atop q_2,p_2}\right)
 =tH_\mathrm{V}\left(
{-\theta^0_1-\theta^\infty_1, \theta^0_1+\theta^1 \atop -\theta^0_2-\theta^1-\theta^\infty_2}; t; q_1, p_1 \right)\\
\qquad{} +tH_\mathrm{V}\left(
{\theta^0_2+\theta^1+\theta^\infty_2, \theta^0_1-\theta^0_2+\theta^1 \atop -\theta^1}; t; q_2, p_2 \right)
 +p_1(q_2-1)\big\{p_2(q_1+q_2)-\theta^0_2-\theta^\infty_2\big\}.\nonumber
\end{gather*}
The gauge parameters satisfy
\begin{gather*}
\frac{t}{u}\frac{{\rm d}u}{{\rm d}t}=-p_1q_2-tq_2+\theta^\infty_1-\theta^\infty_2,\qquad
\frac{t}{v}\frac{{\rm d}v}{{\rm d}t}=-tq_1+\theta^\infty_1-\theta^\infty_3.
\end{gather*}

Concerning the linear system of the spectral type $(11)(1),(11)(1)$, as the expression of the Hamiltonian $H_{\mathrm{FS}}^{A_3}$ changed (see Remark~\ref{rem:FSA3}), the parametrization of the linear system also changed slightly. The variables~$x$ and~$Y$ also changed in the obvious way.

As the result, the Riemann scheme is given by
\begin{gather*}
\left(
\begin{matrix}
 x=0 & x=\infty \\
\overbrace{\begin{matrix}
 1 & \theta^0_1 \\
 0 & \theta^0_2 \\
 0 & 0
 \end{matrix}}
&
\overbrace{\begin{matrix}
 t & \theta^\infty_1 \\
 0 & \theta^\infty_2 \\
 0 & \theta^\infty_3
 \end{matrix}}
\end{matrix}
\right) ,
\end{gather*}
and the Fuchs--Hukuhara relation is written as
$\theta_1^0 +\theta_2^0 +\theta_1^\infty +\theta_2^\infty+\theta_3^\infty =0$.
The Lax pair is expressed as
\begin{gather}
\frac{\partial Y}{\partial x}= \left( \frac{A_2}{x^2}+\frac{A_1}{x}+A_0 \right)Y ,\qquad
\frac{\partial Y}{\partial t}=\left( \frac{1}{t}A_0 x+B_0 \right)Y ,\label{eq:(11)(1),(11)(1)}
\end{gather}
where
\begin{gather*}
A_{\xi} =
U^{-1}P^{-1}\hat{A}_{\xi}PU , \qquad \xi=1,2,\qquad
P=
\begin{pmatrix}
1 & 0 & 0 \\
-p_2 & 1 & 0 \\
-q_1 & \frac{q_1q_2-t}{\theta^\infty_2-\theta^\infty_3} & 1
\end{pmatrix},\qquad U={\rm diag}(1,u,v),\\
A_0 =
\begin{pmatrix}
-t & 0 & 0 \\
0 & 0 & 0 \\
0 & 0 & 0
\end{pmatrix},\qquad
\hat{A}_1=
\begin{pmatrix}
p_1q_1-p_2q_2-\theta^\infty_1 & -q_2 & p_1 \\
p_2q_2-\theta^0_2-\theta^\infty_2 & p_2q_2-\theta^\infty_2 & -p_1p_2 \\
t & t & -p_1q_1-\theta^\infty_3
\end{pmatrix},\\
\hat{A}_2=
\begin{pmatrix}
1 \\
0 \\
0
\end{pmatrix}
\begin{pmatrix}
1 & 1 & -\frac{1}{t}(p_1q_1+\theta^\infty_3)
\end{pmatrix},\qquad
B_0=\frac{1}{t}
\begin{pmatrix}
0 & (A_1)_{12} & (A_1)_{13} \\
(A_1)_{21} & 0 & 0 \\
(A_1)_{31} & 0 & 0
\end{pmatrix}.
\end{gather*}

The Hamiltonian is given by
\begin{gather*}
 tH_{\mathrm{FS}}^{A_3}
\left({ \alpha, \beta \atop \gamma, \delta};t;
 {q_1,p_1 \atop q_2,p_2}\right) =tH_{\mathrm{III}(D_6)}\big( {\theta^\infty_3, \theta^\infty_3-\theta^\infty_1+1}; t; q_1, p_1 \big)\\
 \qquad{} +tH_{\mathrm{III}(D_6)}\big( {-}\theta^0_2-\theta^\infty_2, \theta^\infty_1-\theta^\infty_2; t; q_2, p_2 \big)
 -p_1p_2(q_1q_2+t).\nonumber
\end{gather*}
The gauge parameters satisfy
\begin{gather*}
\frac{t}{u}\frac{{\rm d}u}{{\rm d}t}=p_1q_1-2p_2q_2+q_2,\qquad
\frac{t}{v}\frac{{\rm d}v}{{\rm d}t}=2p_1q_1-p_2q_2-q_1+1.
\end{gather*}

Below are the data for degenerations.

\noindent
$\boldsymbol{2+1+1 \to 2+1+1}$

\noindent
$(1)(1)(1), 21, 21 \to (1)_2(1), 21, 21$
\begin{gather*}
\theta^\infty_1=\tilde{\theta}^\infty_1+\varepsilon^{-1},\qquad \theta^\infty_3=-\varepsilon^{-1},\qquad
t_2=\varepsilon \tilde{t}_2,\ H_{t_2}=\varepsilon^{-1}\tilde{H}_2,\\
q_1=\tilde{q}_1+\frac{\theta^\infty_2}{\tilde{p}_1},\qquad p_1=\tilde{p}_1,\qquad
q_2=\tilde{p}_2,\qquad p_2=-\tilde{q}_2,\\
Y={t_2}^{\varepsilon^{-1}}
\operatorname{diag}(1,\, u,\, v)^{-1}
\begin{pmatrix}
0 & \varepsilon & 0 \\
0 & 0 & \varepsilon p_1q_1 \\
t_2q_2 & \varepsilon q_2(p_2q_2-\theta^\infty_3) & 0
\end{pmatrix}
\tilde{Y}.
\end{gather*}

\noindent
$\boldsymbol{3+1 \to 3+1}$

\noindent
$((1)(1))((1)), 21 \to ((1))(1)_2, 21$
\begin{gather*}
\theta^\infty_2=\tilde{\theta}^\infty_2+\varepsilon^{-1},\qquad \theta^\infty_3=-\varepsilon^{-1},\qquad
t_1=\tilde{t}_2-\varepsilon\tilde{t}_1,\qquad t_2=\tilde{t}_2,\\
H_{t_1}=-\frac{1}{\varepsilon}\left( \tilde{H}_1-\frac{\tilde{p}_1\tilde{q}_1}{\tilde{t}_1} \right),\qquad
H_{t_2}=\varepsilon^{-1}\tilde{H}_1+\tilde{H}_2-\frac{\tilde{p}_1\tilde{q}_1}{\varepsilon\tilde{t}_1},\\
q_1=\varepsilon\tilde{t}_1\tilde{p}_1+\tilde{q}_2,\qquad p_1=-\frac{\tilde{q}_1}{\varepsilon\tilde{t}_1},\qquad
q_2=\tilde{q}_2,\qquad p_2=\tilde{p}_2+\frac{\tilde{q}_1}{\varepsilon\tilde{t}_1},\\
x=-\tilde{x},\qquad
Y=\exp\left( \frac12 x^2-t_2 x \right)
\begin{pmatrix}
1 & & \\
 & u & \\
 & & v
\end{pmatrix}^{-1}
\begin{pmatrix}
1 & 0 & 0 \\
0 & 1 & \frac{q_1}{(t_1-t_2)q_2} \\
0 & 0 & -\frac{q_1}{(t_1-t_2)q_2}
\end{pmatrix}
\tilde{Y}.
\end{gather*}

\noindent
$\boldsymbol{3+1 \to 5/2+1}$

\noindent
$((1)(1))((1)), 21 \to (((1)))_2 1, 21$
\begin{gather*}
\theta^\infty_1=\tilde{\theta}^\infty_1-\varepsilon^{-6},\qquad \theta^\infty_2=\varepsilon^{-6},\qquad
\theta^\infty_3=\tilde{\theta}^\infty_2,\\ t_1=\varepsilon \tilde{t}_1-2\varepsilon^{-3},\qquad t_2=-\varepsilon^{-1}\tilde{t}_2-\varepsilon^{-3},\\
q_1=-\varepsilon \tilde{p}_1,\qquad p_1=\varepsilon^{-1}\tilde{q}_1-\varepsilon^{-3},\qquad
q_2=\varepsilon^{-1} \tilde{q}_2,\qquad p_2=\varepsilon \tilde{p}_2,\\
H_{t_1}=\varepsilon^{-1}\tilde{H}_1,\qquad H_{t_2}=-\varepsilon \tilde{H}_2,\\
x=\varepsilon \tilde{x},\qquad
Y=\exp\left( \frac12 x^2-t_2 x-\varepsilon^{-2}\tilde{t}_1 \right)
\begin{pmatrix}
1 & & \\
 & u & \\
 & & v
\end{pmatrix}^{-1}
\begin{pmatrix}
0 & -1 & 0 \\
1/\varepsilon & -p_1 & 0 \\
0 & 0 & \varepsilon
\end{pmatrix}
\tilde{Y}.
\end{gather*}

\pagebreak  \noindent
$\boldsymbol{2+1+1 \to 3+1}$

\noindent
$(1)_2(1), 21, 21 \to ((1))(1)_2, 21$
\begin{gather*}
\theta^1=\tilde{\theta}^\infty_1+\varepsilon^{-2},\qquad \theta^\infty_1=\tilde{\theta}^\infty_2,\qquad \theta^\infty_2=-\varepsilon^{-2},\qquad
t_1=\varepsilon^{-2}-\varepsilon^{-1}\tilde{t}_2,\qquad t_2=\varepsilon^{-1} \tilde{t}_1,\\
H_{t_1}=-\varepsilon \tilde{H}_2,\qquad H_{t_2}=\varepsilon \left( \tilde{H}_1-\frac{\tilde{p}_1\tilde{q}_1}{\tilde{t}_1} \right),\qquad
q_1=\left( 1+\frac{\tilde{q}_1}{\varepsilon( \tilde{q}_1\tilde{q}_2+\tilde{t}_1 )} \right)^{-1},\\
p_1=-\big(\tilde{p}_1\tilde{q}_1-\tilde{p}_2\tilde{q}_2-\theta^0\big)
\left( \frac{\tilde{q}_1}{\varepsilon(\tilde{q}_1\tilde{q}_2+\tilde{t}_1)}+2+\frac{\varepsilon(\tilde{q}_1\tilde{q}_2+\tilde{t}_1)}{\tilde{q}_1} \right)
-\frac{\tilde{q}_1\tilde{q}_2+\tilde{t}_1}{\varepsilon\tilde{q}_1}-\varepsilon^{-2}, \\
q_2=-\tilde{t}_1\frac{\tilde{q}_1(\tilde{p}_1\tilde{q}_1-\theta^0-\tilde{\theta}^\infty_1)
+\tilde{t}_1\tilde{p}_2}{\tilde{q}_1(\tilde{q}_1\tilde{q}_2+\tilde{t}_1)},\qquad
p_2=1+\frac{\tilde{q}_1\tilde{q}_2}{\tilde{t}_1}, \\
x=\varepsilon \tilde{x},\qquad
Y=f
\begin{pmatrix}
\varepsilon t_2 p_2q_2 & -t_2 & \frac{1}{\varepsilon p_2}((q_1-1)(\varepsilon^{-2}-p_1q_1)+\theta^1)+\frac{q_2}{\varepsilon(q_1-1)} \\
\varepsilon {t_2}^2 p_2 & 0 & \frac{t_2}{\varepsilon(q_1-1)} \\
\varepsilon{t_2}^2 p_2(q_1-1) & 0 & 0
\end{pmatrix}
\tilde{Y},
\end{gather*}
where $f$ satisf\/ies
\begin{gather*}
\frac{1}{f}\frac{\partial f}{\partial \tilde{t}_1}=-\tilde{p}_1\tilde{q}_1-1, \qquad
\frac{1}{f}\frac{\partial f}{\partial \tilde{t}_2}=-\tilde{q}_2-\tilde{t}_2.
\end{gather*}

\noindent
$\boldsymbol{2+1+1 \to 5/3+1+1}$

\noindent
$(1)_2(1),21,21 \to ((1))_3,21,21$
\begin{gather*}
\theta^\infty_1=\tilde{\theta}^\infty_1+\varepsilon^{-2},\qquad \theta^\infty_2=-\varepsilon^{-2},\qquad
t_1=-\varepsilon\tilde{t}_1(1+\varepsilon\tilde{t}_2),\qquad t_2=\varepsilon^{-1}\tilde{t}_1,\\
H_{t_1}=-\frac{1}{\varepsilon^2\tilde{t}_1}\tilde{H}_2,\qquad
H_{t_2}=\varepsilon\tilde{H}_1-\frac{1+\varepsilon\tilde{t}_2}{\tilde{t}_1}\tilde{H}_2
-\frac{\varepsilon\tilde{p}_1\tilde{q}_1}{\tilde{t}_1},\\
q_1=1+\frac{\tilde{q}_1\tilde{q}_2}{\tilde{t}_1}+\frac{\tilde{q}_1}{\varepsilon\tilde{t}_1},\qquad
p_1=\varepsilon\tilde{t}_1\left(
\tilde{p}_1-\frac{\tilde{p}_2\tilde{q}_2+\theta^0}{\tilde{q}_1}\right)+\frac{\tilde{t}_1}{\varepsilon\tilde{q}_1},\\
q_2=\varepsilon\tilde{t}_1\left(\tilde{p}_1-\frac{\tilde{p}_2\tilde{q}_2+\theta^0}{\tilde{q}_1}\right)
+\tilde{t}_1\frac{1-\varepsilon\tilde{p}_2}{\varepsilon\tilde{q}_1},\qquad
p_2=1+\frac{\tilde{q}_1\tilde{q}_2}{\tilde{t}_1},\\
Y=\tilde{t}_1{}^{\varepsilon^{-2}}e^{\varepsilon^{-1}\tilde{t}_2-\frac12 {\tilde{t}_2}{}^2}
\begin{pmatrix}
1 & -1/t_1 & 0 \\
0 & 1 & \frac{q_1-p_2}{\varepsilon^2t_2q_1} \\
0 & 0 & \frac{q_1-p_2}{\varepsilon^2 t_2}
\end{pmatrix}
\tilde{Y}.
\end{gather*}

\noindent
$\boldsymbol{2+1+1 \to 2+2}$

\noindent
$(1)_2(1), 21, 21 \to (1)_2(1), (2)(1)$
\begin{gather*}
\theta^0=\tilde{\theta}^0-\varepsilon^{-1},\qquad \theta^1=\varepsilon^{-1},\qquad
t_1=\varepsilon \tilde{t}_1,\qquad t_2=-\varepsilon \tilde{t}_2,\\
H_{t_1}=\varepsilon^{-1}\tilde{H}_1,\qquad H_{t_2}=-\varepsilon^{-1}\left( \tilde{H}_2-\frac{\tilde{p}_2 \tilde{q}_2}{\tilde{t}_2} \right),\\
q_1=-\frac{1}{\varepsilon \tilde{q}_1},\qquad p_1=\varepsilon \tilde{q}_1\big(\tilde{p}_1\tilde{q}_1+\theta^\infty_2\big),\qquad
q_2=\varepsilon \tilde{t}_2 \tilde{p}_2,\qquad p_2=-\frac{\tilde{q}_2}{\varepsilon \tilde{t}_2},\\
x=\varepsilon^{-1}\tilde{x}, \qquad
Y=
\begin{pmatrix}
1 & 0 & 0 \\
0 & \varepsilon & 0 \\
0 & 0 & \varepsilon q_1
\end{pmatrix}
\tilde{Y}.
\end{gather*}

\noindent
$\boldsymbol{2+2 \to 2+2}$

\noindent
$(1)(1)(1), (2)(1) \to (1)_2(1), (2)(1)$
\begin{gather*}
\theta^\infty_1=\tilde{\theta}^\infty_1-\varepsilon^{-1},\qquad \theta^\infty_3=\varepsilon^{-1},\qquad
t_2=\varepsilon\tilde{t}_2,\qquad H_{t_2}=\varepsilon^{-1}\left( \tilde{H}_2-\frac{\tilde{p}_2\tilde{q}_2}{\tilde{t}_2} \right),\\
q_1=\tilde{q}_1+\frac{\theta^\infty_2}{\tilde{p}_1},\qquad p_1=\tilde{p}_1,\qquad
q_2=\frac{\tilde{t}_2}{\tilde{q}_2},\qquad p_2=-\frac{\tilde{q}_2(\tilde{p}_2\tilde{q}_2-\varepsilon^{-1})}{\tilde{t}_2},\\
Y={\tilde{t}_2}{}^{-\varepsilon^{-1}}\frac{p_2}{t_2}
\cdot\operatorname{diag}(1, \, u, \, v)^{-1}\cdot
\begin{pmatrix}
0 & \varepsilon q_2 & 0 \\
0 & 0 & \varepsilon p_1 q_2 \\
t_2 & \varepsilon p_2q_2 & 0
\end{pmatrix}
\tilde{Y}.
\end{gather*}

\noindent
$\boldsymbol{4 \to 4}$

\noindent
$(((1)(1)))(((1))) \to (((1)))(1)_2$
\begin{gather*}
\theta^\infty_2=\tilde{\theta}^\infty_2+\varepsilon^{-1},\qquad \theta^\infty_3=-\varepsilon^{-1},\qquad
t_1=\tilde{t}_1,\ t_2=\tilde{t}_1+\varepsilon\tilde{t}_2,\\
H_{t_1}=\tilde{H}_1-\varepsilon^{-1}\tilde{H}_2,\qquad
H_{t_2}=\varepsilon^{-1}\tilde{H}_2,\\
q_1=\tilde{q}_1,\qquad p_1=\tilde{p}_1-\varepsilon^{-1}\tilde{p}_2,\qquad
q_2=\tilde{q}_1+\varepsilon\tilde{q}_2,\qquad p_2=\varepsilon^{-1}\tilde{p}_2,\\
x=-\tilde{x},\qquad
Y=\exp\left( -\frac{x^3}{3}-t_2x \right)
\begin{pmatrix}
\varepsilon & 0 & 0 \\
0 & 0 & 1 \\
0 & -\varepsilon & -1
\end{pmatrix}
\tilde{Y}.
\end{gather*}

\noindent
$\boldsymbol{3+1 \to 4}$

\noindent
$((1))(1)_2, 21 \to (((1)))(1)_2$
\begin{gather*}
\theta^0=-\varepsilon^{-6},\qquad \theta^\infty_1=\tilde{\theta}^\infty_1+\varepsilon^{-6},\qquad
t_1=\varepsilon\tilde{t}_2,\qquad t_2=\varepsilon\tilde{t}_1-2\varepsilon^{-3},\\
H_{t_1}=\varepsilon^{-1}\left( \tilde{H}_2-\frac{\tilde{p}_2\tilde{q}_2}{\tilde{t}_2} \right),\qquad
H_{t_2}=\varepsilon^{-1}\tilde{H}_1,\\
 q_1=\varepsilon^2\tilde{t}_2\tilde{p}_2,\qquad p_1=-\frac{\tilde{q}_2}{\varepsilon^2\tilde{t}_2},\qquad
q_2=\varepsilon^{-1}\tilde{q}_1+\varepsilon^{-3},\qquad p_2=\varepsilon\tilde{p}_1,\\
x=\varepsilon^{-1}\tilde{x}+\varepsilon^{-3},\qquad
Y=\exp\big({-}\varepsilon^{-2}\tilde{t}_1\big)
\begin{pmatrix}
1 & 0 & 0 \\
0 & \varepsilon^{-1} & 0 \\
0 & 0 & \tilde{t}_2
\end{pmatrix}
\tilde{Y}.
\end{gather*}

\noindent
$\boldsymbol{3+1 \to 7/3+1}$

\noindent
$((1))(1)_2, 21 \to ((((1))))_3, 21$
\begin{gather*}
\theta^\infty_1=\tilde{\theta}^\infty_1-3\varepsilon^{-4},\qquad \theta^\infty_2=3\varepsilon^{-4},\qquad
t_1=\varepsilon^{-3}\tilde{t}_1+\varepsilon^{-4}\tilde{t}_2+\varepsilon^{-6},\qquad
t_2=\tilde{t}_2-3\varepsilon^{-2},\\
H_{t_1}=\varepsilon^3\tilde{H}_1,\qquad
H_{t_2}=-\varepsilon^{-1}\tilde{H}_1+\tilde{H}_2,\\
q_1=\varepsilon^{-3}\tilde{q}_1+\varepsilon^{-4},\qquad
p_1=\varepsilon^3\tilde{p}_1-\varepsilon^2\tilde{q}_2,\qquad
q_2=\tilde{q}_2,\qquad p_2=\tilde{p}_2-\varepsilon^{-1}\tilde{q}_1-2\varepsilon^{-2},\\
x=\varepsilon\tilde{x},\qquad
Y=\exp(\varepsilon^{-1}\tilde{x}+\varepsilon^{-1}\tilde{t}_1-\varepsilon^{-2}\tilde{t}_2)
\begin{pmatrix}
0 & -\varepsilon & 0 \\
0 & \varepsilon^{-1}-\varepsilon p_2-\varepsilon^3 q_1 & 1 \\
\varepsilon^{-1} & -\varepsilon^{-3} & -\varepsilon^{-2}
\end{pmatrix}
\tilde{Y}.
\end{gather*}

\noindent
$\boldsymbol{5/2+1 \to 7/3+1}$

\noindent
$(((1)))_2 1, 21 \to ((((1))))_3, 21$
\begin{gather*}
\theta^\infty_1=\tilde{\theta}^\infty_1+\varepsilon^{-12},\qquad \theta^\infty_2=-\varepsilon^{-12},\qquad
t_1=-\varepsilon \tilde{t}_1+\frac{\tilde{t}_2}{\varepsilon^2}+\frac{3}{4\varepsilon^8},\qquad
t_2=-\varepsilon^2 \tilde{t}_2+\frac{3}{2\varepsilon^4},\\
H_{t_1}=-\varepsilon^{-1}\tilde{H}_1,\qquad
H_{t_2}=-\varepsilon^{-5}\tilde{H}_1-\varepsilon^{-2}\tilde{H}_2-\varepsilon^{-2}\tilde{q}_2,\\
q_1=-\frac{\tilde{q}_1}{\varepsilon}-\frac{1}{2\varepsilon^4},\qquad p_1=-\varepsilon\tilde{p}_1-\frac{\tilde{q}_2}{\varepsilon^2},\qquad
q_2=\varepsilon^2\tilde{q}_2-\varepsilon^8,\qquad p_2=\frac{\tilde{p}_2-\tilde{t}_2}{\varepsilon^2}+\frac{\tilde{q}_1}{\varepsilon^5}+\frac{1}{\varepsilon^8},\\
x=-\varepsilon \tilde{x},\qquad
Y=\exp\left( \frac{\tilde{x}}{\varepsilon^3}+\frac{\tilde{t}_1}{2\varepsilon^3}+\frac{\tilde{t}_2}{2\varepsilon^6} \right)
\begin{pmatrix}
\varepsilon^3 & \varepsilon q_1+\frac{1}{2\varepsilon^3} & 1 \\
0 & \varepsilon & 0 \\
0 & -\varepsilon p_2-\frac{q_1}{\varepsilon^3}-\frac{1}{2\varepsilon^7} & -\varepsilon^{-4}
\end{pmatrix}
\tilde{Y}.
\end{gather*}

\noindent
$\boldsymbol{5/3+1+1 \to 2+5/3}$

\noindent
$((1))_3, 21, 21 \to ((1))_3, (2)(1)$
\begin{gather*}
\theta^0=-\varepsilon^{-1},\qquad \theta^1=\tilde{\theta}^0+\varepsilon^{-1},\qquad t_1=\varepsilon\tilde{t}_1,\\
H_{t_1}=\frac{1}{\varepsilon}\left( \tilde{H}_1-\frac{\tilde{p}_1\tilde{q}_1+\tilde{p}_2\tilde{q}_2}{\tilde{t}_1} \right),\qquad
H_{t_2}=\tilde{H}_2+\tilde{p}_2\tilde{q}_1,\\
q_1=\frac{\tilde{t}_1}{\tilde{q}_1},\qquad
p_1=-\frac{\tilde{q}_1}{\varepsilon \tilde{t}_1}(\varepsilon \tilde{p}_1\tilde{q}_1+\varepsilon \tilde{p}_2\tilde{q}_2+1),\qquad
q_2=-\tilde{p}_2\tilde{q}_1,\qquad
p_2=-\tilde{p}_2\tilde{q}_1+\frac{\tilde{q}_2}{\tilde{q}_1}+\tilde{t}_2,\\
x=\frac{\tilde{x}}{\varepsilon\tilde{t}_1},\qquad
Y={\tilde{t}_1}{}^{-\tilde{\theta}^0}
\begin{pmatrix}
1 & 0 & 0 \\
0 & \varepsilon\tilde{t}_1 & 0 \\
0 & 0 & \varepsilon^2{\tilde{t}_1}{}^2
\end{pmatrix}
\tilde{Y}.
\end{gather*}

\noindent
$\boldsymbol{2+2 \to 4}$

\noindent
$(1)_2(1), (2)(1) \to (((1)))(1)_2$
\begin{gather*}
\theta^0=\tilde{\theta}^\infty_1+2\varepsilon^{-3},\qquad \theta^\infty_1=\tilde{\theta}^\infty_2,\qquad \theta^\infty_2=-2\varepsilon^{-3},\qquad
t_1=\varepsilon^{-4}\tilde{t}_1+\varepsilon^{-6},\qquad t_2=\varepsilon^{-4}\tilde{t}_2,\\
H_{t_1}=\varepsilon^4 \tilde{H}_1,\qquad
H_{t_2}=\varepsilon^4\left( \tilde{H}_2+\frac{\tilde{p}_1}{\tilde{t}_2\tilde{p}_2}-\frac{\tilde{p}_2\tilde{q}_2}{\tilde{t}_2} \right),\qquad
q_1=\varepsilon^{-3}-\varepsilon^{-2}\left( \tilde{q}_1+\frac{1}{\tilde{p}_2} \right),\\
p_1=1-\varepsilon^2 \tilde{p}_1,\qquad
q_2=-\varepsilon^{-2}\tilde{t}_2\tilde{p}_2,\qquad
p_2=\frac{\varepsilon^2}{\tilde{t}_2}
\left( \tilde{q}_2-\frac{\tilde{\theta}^\infty_2}{\tilde{p}_2}-\frac{\tilde{p}_1}{\tilde{p}_2{}^2} \right),\\
x=-\varepsilon^4(\tilde{x}+\varepsilon^{-1}),\qquad
Y=f
\begin{pmatrix}
0 & -1/t_2 & \frac{p_1-1}{q_2} \\
0 & 0 & 1 \\
-\frac{q_2}{\varepsilon t_2} & 0 & 0
\end{pmatrix}
\tilde{Y}.
\end{gather*}
Here $f$ satisf\/ies
\begin{gather*}
\frac{1}{f}\frac{\partial f}{\partial \tilde{t}_1}=-\tilde{q}_1, \qquad
\frac{1}{f}\frac{\partial f}{\partial \tilde{t}_2}=\frac{\tilde{p}_2\tilde{q}_2+1}{\tilde{t}_2}.
\end{gather*}

\noindent
$\boldsymbol{2+2 \to 2+5/3}$

\noindent
$(1)_2(1), (2)(1) \to ((1))_3, (2)(1)$
\begin{gather*}
\theta^\infty_1=\tilde{\theta}^\infty_1+\varepsilon^{-2},\qquad \theta^\infty_2=-\varepsilon^{-2},\qquad
t_1=\varepsilon\tilde{t}_1,\qquad t_2=\tilde{t}_1\tilde{t}_2+\varepsilon^{-1}\tilde{t}_1,\\
H_{t_1}=\varepsilon^{-1}\left( \tilde{H}_1-\frac{1}{\tilde{t}_1}\big(\tilde{t}_2+\varepsilon^{-1}\big)\tilde{H}_2
-\frac{\tilde{p}_1\tilde{q}_1}{\tilde{t}_1}-\frac{\tilde{p}_2\tilde{q}_1}{\varepsilon\tilde{t}_1} \right),\qquad
H_{t_2}=\frac{1}{\tilde{t}_1}\tilde{H}_2,\\
q_1=\varepsilon(\tilde{q}_1-\varepsilon \tilde{q}_2),\qquad
p_1=-\frac{\tilde{p}_2}{\varepsilon^2},\qquad
q_2=\frac{\tilde{t}_1}{\varepsilon}(\tilde{p}_2+\varepsilon\tilde{p}_1),\qquad
p_2=-\frac{\tilde{q}_1}{\tilde{t}_1},\\
x=\frac{\tilde{x}}{\tilde{t}_1},\qquad
Y=\tilde{t}_1{}^{\varepsilon^{-2}+1-\theta^0}
\begin{pmatrix}
-\frac{\varepsilon}{\tilde{t}_1} & 1/\tilde{t}_1 & 0 \\
0 & -\varepsilon & 1 \\
0 & 0 & 1
\end{pmatrix}
\tilde{Y}.
\end{gather*}

\noindent
$\boldsymbol{2+1+1 \to 3/2+1+1}$

\noindent
$(11)(1), 21, 111 \to (1)_2 1, 21, 111$
\begin{gather*}
\theta^0_1=\tilde{\theta}^0-\varepsilon^{-1}, \qquad \theta^0_2=\varepsilon^{-1}, \qquad t=-\varepsilon\tilde{t}, \qquad H=-\varepsilon^{-1}\tilde{H}, \\
q_1=1-\tilde{p}_1-\tilde{p}_2, \qquad p_1=\tilde{q}_2, \qquad q_2=1-\tilde{p}_1, \qquad p_2=\tilde{q}_1-\tilde{q}_2.
\end{gather*}

\noindent
$\boldsymbol{2+2 \to 2+3/2}$

\noindent
$(11)(1), (11)(1) \to (1)_2 1, (11)(1)$
\begin{gather*}
\theta^\infty_1=\tilde{\theta}^\infty_1+\varepsilon^{-1}, \qquad \theta^\infty_2=-\varepsilon^{-1},\qquad \theta^\infty_3=\tilde{\theta}^\infty_2, \qquad
t=\varepsilon\tilde{t}, \qquad H=\varepsilon^{-1}\tilde{H}-\frac{\tilde{p}_1\tilde{q}_1}{\varepsilon\tilde{t}}, \\
q_1=\frac{\tilde{t}}{\tilde{q}_1}, \qquad p_1=-\frac{\tilde{q}_1}{\tilde{t}}\big(\tilde{p}_1\tilde{q}_1+\tilde{\theta}^\infty_2\big), \qquad
q_2=-\frac{1}{\tilde{p}_2}, \qquad p_2=-\tilde{p}_2\big(\tilde{p}_2\tilde{q}_2-\varepsilon^{-1}+\theta^0_2\big), \\
Y=t^{\varepsilon^{-1}}U^{-1}P^{-1}
\begin{pmatrix}
0 & 1 & 0 \\
-1/q_2 & 0 & -1/q_2 \\
\varepsilon q_1 & 0 & 0
\end{pmatrix}\tilde{Y}.
\end{gather*}

\noindent
$\boldsymbol{3/2+1+1 \to 2+3/2}$

\noindent
$(1)_2 1, 21, 111 \to (1)_2 1, (11)(1)$
\begin{gather*}
\theta^0_1=\tilde{\theta}^0_1-\varepsilon^{-1},\qquad \theta^1=\varepsilon^{-1},\qquad
t=-\varepsilon\tilde{t},\qquad
H=-\varepsilon^{-1}\tilde{H}+\frac{\tilde{p}_1\tilde{q}_1}{\varepsilon\tilde{t}},\\
q_1=\frac{\tilde{t}\tilde{p}_1}{\tilde{p}_1\tilde{q}_1+\theta^\infty_2},\qquad p_1=-\frac{\tilde{q}_1}{\tilde{t}}\big(\tilde{p}_1\tilde{q}_1+\theta^\infty_2\big),\qquad
q_2=-\varepsilon\tilde{q}_2,\qquad p_2=-\varepsilon^{-1}\tilde{p}_2,\\
x=\varepsilon^{-1}\tilde{x},\qquad
Y=
\begin{pmatrix}
1 & 0 & 0 \\
0 & -1 & 0 \\
0 & 0 & -1/t
\end{pmatrix}\tilde{Y}.
\end{gather*}

\noindent
$\boldsymbol{3/2+1+1 \to 4/3+1+1}$

\noindent
$(1)_2 1, 21, 111 \to (1)_3, 21, 111$
\begin{gather*}
\theta^\infty_1=\tilde{\theta}^\infty_1-\varepsilon^{-1}, \qquad \theta^\infty_2=\varepsilon^{-1}, \qquad
t=\varepsilon\tilde{t},\qquad H=\varepsilon^{-1}\tilde{H}-\frac{\tilde{p}_1\tilde{q}_1}{\varepsilon\tilde{t}}, \qquad
q_1=\varepsilon \tilde{q}_2, \\ p_1=\varepsilon^{-1}\tilde{p}_2, \qquad q_2=\frac{\tilde{t}}{\tilde{q}_1}, \qquad
p_2=-\frac{\tilde{q}_1}{\tilde{t}}\big(\tilde{p}_1\tilde{q}_1-\theta^0_1+\theta^0_2\big),\qquad
Y=
\begin{pmatrix}
t & 0 & 0 \\
0 & 1 & 0 \\
0 & 0 & -q_2/t
\end{pmatrix}\tilde{Y}.
\end{gather*}

\noindent
$\boldsymbol{2+3/2 \to 3/2+3/2}$

\noindent
$(1)_2 1, (11)(1) \to (1)_2 1, (1)_2 1$
\begin{gather*}
\theta^0_1=\tilde{\theta}^0+\varepsilon^{-1}, \qquad \theta^0_2=-\varepsilon^{-1},\qquad
t=\varepsilon\tilde{t}, \qquad H=\varepsilon^{-1}\tilde{H}-\frac{\tilde{p}_2\tilde{q}_2}{\varepsilon\tilde{t}},\\
q_1=-\varepsilon\tilde{t}\tilde{p}_2, \qquad p_1=\frac{\tilde{q}_2}{\varepsilon\tilde{t}}-\frac{1}{\varepsilon\tilde{q}_1}, \qquad
q_2=\tilde{q}_1, \qquad p_2=\tilde{p}_1+\frac{-\theta^\infty_2+\varepsilon^{-1}}{\tilde{q}_1}-\frac{\tilde{t}\tilde{p}_2}{\tilde{q}_1^2},\\
x=\varepsilon^{-1}\tilde{x}, \qquad
Y=t^{-\theta^\infty_1}
\begin{pmatrix}
1 & 0 & 0 \\
0 & 1 & 0 \\
0 & 0 & q_2
\end{pmatrix}\tilde{Y}.
\end{gather*}

\noindent
$\boldsymbol{2+3/2 \to 2+4/3}$

\noindent
$(1)_2 1, (11)(1) \to (1)_3, (11)(1)$
\begin{gather*}
\theta^\infty_1=\tilde{\theta}^\infty_1-\varepsilon^{-1},\qquad \theta^\infty_2=\varepsilon^{-1},\qquad
\theta^0_1=\tilde{\theta}^0_2,\qquad \theta^0_2=\tilde{\theta}^0_1,\qquad q_1=-\varepsilon\tilde{t}\tilde{p}_2,\qquad p_1=\frac{\tilde{q}_2}{\varepsilon\tilde{t}},\\
q_2=\tilde{q}_1,\qquad p_2=\tilde{p}_1,\qquad
t=\varepsilon\tilde{t},\qquad H=\varepsilon^{-1}\tilde{H}-\frac{\tilde{p}_2\tilde{q}_2}{\varepsilon\tilde{t}},\qquad
Y=
\begin{pmatrix}
t & 0 & 0 \\
0 & 1 & 0 \\
0 & 0 & 1
\end{pmatrix}\tilde{Y}.
\end{gather*}

\noindent
$\boldsymbol{4/3+1+1 \to 2+4/3}$

\noindent
$(1)_3, 21, 111 \to (1)_3, (11)(1)$
\begin{gather*}
 \theta^1=-\varepsilon^{-1}, \qquad \theta^0_2=\tilde{\theta}^0_2+\varepsilon^{-1},\qquad
 q_1=-\varepsilon \tilde{t}\tilde{p}_1, \qquad p_1=\frac{\tilde{q}_1}{\varepsilon\tilde{t}}, \qquad t=\varepsilon\tilde{t},\\
H=\varepsilon^{-1}\tilde{H}-\frac{\tilde{p}_1\tilde{q}_1}{\varepsilon\tilde{t}},\qquad x=-\varepsilon^{-1}\tilde{x},\qquad
Y=
\begin{pmatrix}
1 & 0 & 0 \\
0 & -\varepsilon & 0 \\
0 & -\varepsilon & \frac{\varepsilon q_1}{t}
\end{pmatrix}\tilde{Y}.
\end{gather*}

\noindent
$\boldsymbol{3/2+3/2 \to 3/2+4/3}$

\noindent
$(1)_2 1, (1)_2 1 \to (1)_3, (1)_2 1$
\begin{gather*}
\theta^\infty_1=-\varepsilon^{-1},\qquad \theta^\infty_2=\tilde{\theta}^\infty_1+\varepsilon,\qquad
t=\varepsilon\tilde{t},\qquad H=\varepsilon^{-1}\tilde{H},\\
q_1=\tilde{q}_1,\qquad p_1=\tilde{p}_1+\frac{\tilde{\theta}^\infty_1+\varepsilon^{-1}}{\tilde{q}_1},\qquad
Y=t^{-\varepsilon^{-1}}
\begin{pmatrix}
t & 0 & 0 \\
0 & 1 & 0 \\
0 & 0 & 1
\end{pmatrix}\tilde{Y}.
\end{gather*}

\noindent
$\boldsymbol{2+4/3 \to 3/2+4/3}$

\noindent
$(1)_3, (11)(1) \to (1)_3, (1)_2 1$
\begin{gather*}
 \theta^0_1=\tilde{\theta}^0-\varepsilon^{-1},\qquad \theta^0_2=\varepsilon^{-1},\qquad
t=\varepsilon\tilde{t},\qquad H=\varepsilon^{-1}\tilde{H},\\
 q_1=\tilde{q}_1,\qquad p_1=\tilde{p}_1+\frac{-\tilde{\theta}^0+\varepsilon^{-1}}{\tilde{q}_1},\qquad
 x=\varepsilon^{-1}\tilde{x},\qquad
Y=
\begin{pmatrix}
1 & 0 & 0 \\
0 & \varepsilon & 0 \\
0 & 0 & \varepsilon q_1
\end{pmatrix}\tilde{Y}.
\end{gather*}

\noindent
$\boldsymbol{3/2+4/3 \to 4/3+4/3}$

\noindent
$(1)_3, (1)_2 1 \to (1)_3, (1)_3$
\begin{gather*}
\theta^0=3\varepsilon^{-1},\qquad \theta^\infty_1=-3\varepsilon^{-1},\qquad
t=\varepsilon\tilde{t},\qquad H=\varepsilon\tilde{H}-\frac{\tilde{p}_1\tilde{q}_1+\tilde{p}_2\tilde{q}_2}{\varepsilon\tilde{t}},\\
q_1=-\frac{\tilde{t}}{\tilde{q}_1},\qquad p_1=\frac{\tilde{q}_1}{\tilde{t}}\left(\tilde{p}_1\tilde{q}_1+\frac13-\frac{1}{\varepsilon}\right),\qquad
q_2=-\frac{\tilde{t}}{\tilde{q}_2},\qquad p_2=\frac{\tilde{q}_2}{\tilde{t}}\left(\tilde{p}_2\tilde{q}_2+\frac23+\frac{1}{\varepsilon}\right),\\
x=\varepsilon^{-1}\tilde{x},\qquad
Y=t^{\frac{1}{\varepsilon}-\frac13}x^{\frac{1}{\varepsilon}-\frac13}
\begin{pmatrix}
1 & 0 & 0 \\
0 & \varepsilon & 0 \\
0 & 0 & t/q_1
\end{pmatrix}\tilde{Y}.
\end{gather*}

\subsection*{Acknowledgements}

The author would like to thank Professor Kazuki Hiroe for his helpful suggestions.
The author would also like to thank Professors Hidetaka Sakai and Akane Nakamura
who provided invaluable comments and continual encouragement.

\pdfbookmark[1]{References}{ref}
\LastPageEnding

\end{document}